\documentclass[preprint,1p,fleqn]{elsarticle}

\usepackage{pgfplots}
\pgfplotsset{compat=1.5}
\usepackage[abs]{overpic} 

\usepackage[utf8]{inputenc}
\usepackage{graphicx,color}
\definecolor{orange}{RGB}{242,121,0}
\definecolor{darkgreen}{RGB}{34,139,34}
\definecolor{purple}{rgb}{0.7,0,0.7}

\usepackage{hyperref}
\usepackage{amsthm}
\usepackage{amsmath}
\usepackage{amsfonts}
\usepackage{amssymb}
\usepackage{subfigure}
\usepackage{algorithm}
\usepackage{algorithmic}
\newcommand{\supp}{\mathop{\mathrm{supp}}}

\newcommand{\trunc}{\mathop{\mathrm{trunc}}}
\newcommand{\Trunc}{\mathop{\mathrm{Trunc}}}

\newcommand{\RR}{{\mathbb{R}}}
\newcommand{\NN}{{\mathbb{N}}}
\newcommand{\VV}{{\mathbb{V}}}
\newcommand{\WW}{{\mathbb{W}}}
\newcommand{\ZZ}{{\mathbb{Z}}}

\newcommand{\cH}{{\mathcal{H}}}
\newcommand{\cT}{{\mathcal{T}}}
\newcommand{\cV}{{\mathcal{V}}}
\newcommand{\cW}{{\mathcal{W}}}
\newcommand{\bzero}{{\boldsymbol{0}}}
\newcommand{\ba}{{\boldsymbol{a}}}

\newcommand{\ff}{{\boldsymbol{f}}}
\newcommand{\bg}{{\boldsymbol{g}}}
\newcommand{\bi}{{\boldsymbol{i}}}

\newcommand{\bn}{{\boldsymbol{n}}}
\newcommand{\bs}{{\boldsymbol{s}}}
\newcommand{\bv}{{\boldsymbol{v}}}
\newcommand{\bx}{{\boldsymbol{x}}}

\newcommand{\bA}{{\boldsymbol{A}}}
\newcommand{\bF}{{\boldsymbol{F}}}
\newcommand{\bG}{{\boldsymbol{G}}}

\newcommand{\bK}{{\boldsymbol{K}}}

\newcommand{\bS}{{\boldsymbol{S}}}
\newcommand{\bX}{{\boldsymbol{X}}}

\newcommand{\bU}{{\boldsymbol{U}}}

\newcommand{\bsigma}{{\boldsymbol{\sigma}}}
\newcommand{\btau}{{\boldsymbol{\tau}}}

\newcommand{\bSigma}{{\boldsymbol{\Sigma}}}
\newcommand{\bXi}{{\boldsymbol{\Xi}}}
\newcommand{\bOmega}{{\boldsymbol{\Omega}}}
\newcommand{\dd}{{\rm{d}}}
\newcommand{\rank}{{\mbox{rank}}}
\newtheorem{theorem}{Theorem}[section]

\newtheorem{definition}[theorem]{Definition}
\newtheorem{configuration}{Strip}

\theoremstyle{definition}

\def\ACM{Adv.\ Comput.\ Math.}
\def\BIT{BIT\ Numer.\ Math.}
\def\CA{Constr.\ Approx.}
\def\CAD{Comput.\ Aided Design}
\def\CAGD{Comput.\ Aided Geom.\ Design}
\def\CG{Comput.\ Graphics}
\def\CMAME{Comput.\ Methods\ Appl.\ Mech.\ Engrg.}

\def\GM{Graph.\ Models}
\def\IJNME{Int.\ J.\ Numer.\ Methods\ Engrg.}
\def\JAMC{J.\ Appl.\ Math.\ Comput.}
\def\JCAM{J.\ Comput.\ Appl.\ Math.}

\def\M3AS{Math.\ Models\ Methods\ Appl.\ Sci.}
\def\NM{Numer.\ Math.}

\usepackage[tracking=true]{microtype}

\begin{document}

\title{Adaptive isogeometric analysis with hierarchical box splines}

\author[indam]{Tadej Kandu\v{c}}
\author[unifi]{Carlotta Giannelli}
\author[unitv]{Francesca Pelosi}
\author[unitv]{Hendrik Speleers}

 \address[indam]{Istituto Nazionale di Alta Matematica, 
 Unit\`a di Ricerca di Firenze c/o DiMaI ``U.~Dini'', 
 Universit\`a di Firenze, Italy}
\address[unifi]{Dipartimento di Matematica e Informatica ``U.~Dini'',
Universit\`a degli Studi di Firenze, Italy}
\address[unitv]{Dipartimento di Matematica, 
Universit\`a degli Studi di Roma ``Tor Vergata'', Italy}

\begin{abstract}
Isogeometric analysis is a recently developed framework based on finite element analysis, where the simple building blocks in geometry and solution space are replaced by more complex and geometrically-oriented compounds.
Box splines are an established tool to model complex geometry, and form an intermediate approach between classical tensor-product B-splines and splines over triangulations. Local refinement can be achieved by considering hierarchically nested sequences of box spline spaces.
Since box splines do not offer special elements to impose boundary conditions for the numerical solution of partial differential equations (PDEs), we discuss a weak treatment of such boundary conditions. Along the domain boundary, an appropriate domain strip is introduced to enforce the boundary conditions in a weak sense. 
The thickness of the strip is adaptively defined in order to avoid
unnecessary computations. Numerical examples show the optimal convergence
rate of box splines and their hierarchical variants for the solution of PDEs.
\end{abstract}


\maketitle


\section{Introduction}
Isogeometric analysis (IgA) is a recently established paradigm based on finite element analysis (FEA) that replaces standard simple building blocks in geometry and solution space with more complex and geometrically-oriented compounds coming from computer-aided design (CAD), see \cite{cottrell2009,hughes2005}. 
The research community has made important steps towards the ambitious goal of incorporating CAD geometries directly into the FEA environment. Even though substantial contributions concerning both theoretical and computational fundamentals have been obtained in the isogeometric context over the last years, critical geometrical challenges still need to be addressed to effectively handle complex multivariate geometries.

Standard CAD geometries are usually represented in terms of tensor-product B-splines and their rational version NURBS \cite{rogers2001,schumaker2007}. Such a tensor-product representation is attractive as it is computationally very easy to work with. Unfortunately, the rigid grid-like structure prevents to model geometries with complicated shapes and to do local mesh refinements.
On the other hand, splines over triangulations do not suffer from such disadvantages
but the geometry of the triangulation and smoothness conditions strongly affect the dimension of the spline spaces and basis properties; hence, all these parameters need to be chosen carefully \cite{lai2007}.
Using triangulations with macro-structures helps to localize and unify the construction of stable spaces and good bases. For example, 
{in \cite{dierckx1997,speleers2013},} 
B-spline-like bases were developed on Powell--Sabin triangulations (these are triangulations endowed with a particular 6-split macro-structure) 
and they turn out to be an interesting ingredient for (adaptive) isogeometric methods \cite{beirao2015,smp2013,smps2012}.

Box splines are an attractive alternative that combine several advantages from tensor-product B-splines and splines over triangulations, somehow constituting an intermediate approach between the two concepts \cite{deboor1993, lai2007}. 
Thanks to the possibility of being defined over special types of regular triangulations, they can handle more complex domains than tensor-product counterparts. The basis elements are shifts of a single box spline, and, consequently, the basis is uniform on the entire domain, an important computational advantage.

Adaptivity can be built on top of a regular structure by considering several (local) layers of different resolutions in a hierarchy. Hierarchically nested spaces of tensor-product B-splines were initially explored in \cite{forsey1988}, and a hierarchical B-spline (HB-spline) basis was developed in \cite{kraft1997,vuong2011}.
An alternative basis for the same hierarchical space has been proposed in \cite{giannelli2012} and is called truncated hierarchical B-spline (THB-spline) basis. The truncated basis possesses several enhanced properties \cite{giannelli2014,speleersM2016} which explain its increasing attention in the design and analysis of adaptive isogeometric methods \cite{buchegger2016,buffa2016c,giannelli2016,hennig2016}.
The (truncated) hierarchical framework has been extended to more general multilevel spline spaces built from any kind of B-spline-like bases \cite{giannelli2014} and generating systems \cite{zore2014}.
The characterization of (T)HB-spline spaces was addressed for bivariate \cite{giannelli2013} and more general multivariate configurations \cite{berdinsky2013,mokris2014a}. This study was also repeated for adaptive spline spaces spanned by $C^2$ quartic box splines \cite{villamizar2016} and $C^1$ quadratic box splines \cite{juttler2016}. Hierarchically nested triangular spline spaces were constructed and analyzed on Powell--Sabin triangulations \cite{speleers2009} and on regular triangulations~\cite{kang2014}.

In this paper we present a hierarchical framework for box spline constructions in the spirit of the multilevel approach proposed in \cite{giannelli2014}, with the aim of solving partial differential equations.
Box spline basis functions have a uniform structure on the entire domain. This means that specific boundary basis functions need to be designed if boundary conditions are imposed in a strong form. However, the corresponding anisotropy in the assembly process can be avoided by considering a weak treatment of the boundary conditions.
Therefore, the hierarchical box spline model is integrated in the immersed boundary method proposed in \cite{baiges2012,kollm2015,pgmss}. One of the main difficulties of immersed boundary methods is accurate
numerical integration over possibly very small mesh cells cut by the boundary; it is a subtle source of ill-conditioning and loss of accuracy. In order to avoid this, we combine the method with a geometry map describing the physical domain, according to the isogeometric philosophy. In this way we combine the benefits offered by immersed boundary and isogeometric methods in the hierarchical box spline context.
Special attention is paid to the construction of a suitable domain strip along the domain boundary, which is needed for the weak imposition of the boundary conditions. The thickness of the boundary strip is adaptively defined in order to avoid unnecessary computations.

The remainder of this paper is organized as follows. Sections~\ref{sec:hbox} describes the general construction of adaptively refined box spline spaces and corresponding (truncated) hierarchical box spline -- (T)HBox-spline -- bases. Fundamental properties for the choice of the underlying box spline spaces 
in the multilevel construction are considered. In Section~\ref{sec:weak} we introduce the weak formulation of the isogeometric method based on hierarchical box splines, together with several solutions for the design of the required domain boundary strip. Section~\ref{sec:exm} presents a selection of examples that demonstrate the use of hierarchical box splines in isogeometric analysis for solving advection-diffusion problems. Finally, Section~\ref{sec:conclusions} concludes the paper.

\section{Hierarchical box splines} \label{sec:hbox}
We consider \emph{(truncated) hierarchical box splines} defined over nested sequences of box spline spaces by following the general approach for the construction of adaptively refined multilevel spline spaces presented in \cite{giannelli2014}. We start by summarizing the definition of box splines and some of their main properties.

\subsection{Box splines and their main properties}

A $d$-variate \emph{box spline} is determined by a set of (possibly repeated) direction vectors $\bv_k\in\ZZ^d\setminus\bzero$, $k=1,\ldots,n$. It is usually denoted by $M_\bXi$ where $\Xi:=[\bv_1 \cdots \bv_n]$. For simplicity, we assume that $n\geq d$ and the submatrix $\Xi_d:=[\bv_1 \cdots \bv_d]$ is non-singular.

\begin{definition}
 The box spline $M_\bXi$ is defined by successive convolutions as follows:
\begin{equation} \label{eq:def-box}
 M_\bXi(\bx):=\int_0^1 M_{\bXi\setminus\bv_n}(\bx-t\bv_n)\, \dd t, \quad n>d,
\end{equation}
starting from
\begin{equation}\label{eq:def-box-start}
 M_{\bXi_d}(\bx):=\left\{\begin{array}{ll}
    1/|\det(\bXi_d)|, \quad & \mbox{if } \bx \in \bXi_d\, [0,1)^d,\\
    0, & \mbox{otherwise}.
  \end{array}\right.
\end{equation}
Here, $\bXi\setminus\bv$ stands for the matrix obtained from $\bXi$ by omitting the vector $\bv$ once,
and $\bXi_d\, [0,1)^n$ is the set of all points in $\RR^d$ obtained after multiplication of $\bXi_d$ with any point in $[0,1)^n$.
\end{definition}

\begin{figure}[t!]
\begin{center}
\subfigure[two-directional mesh]{
\includegraphics[height=4cm]{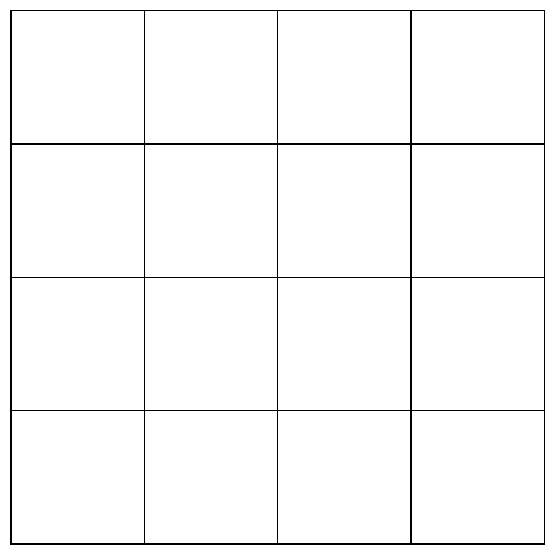}
} \hspace*{0.5cm}
\subfigure[three-directional mesh]{
\includegraphics[height=4cm]{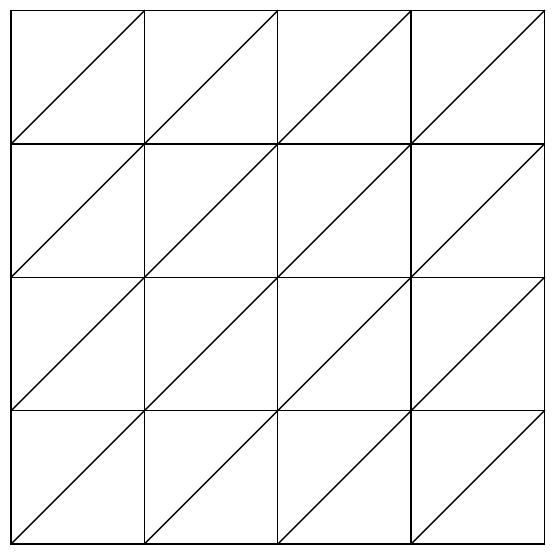}
} \hspace*{0.5cm}
\subfigure[four-directional mesh]{
\includegraphics[height=4cm]{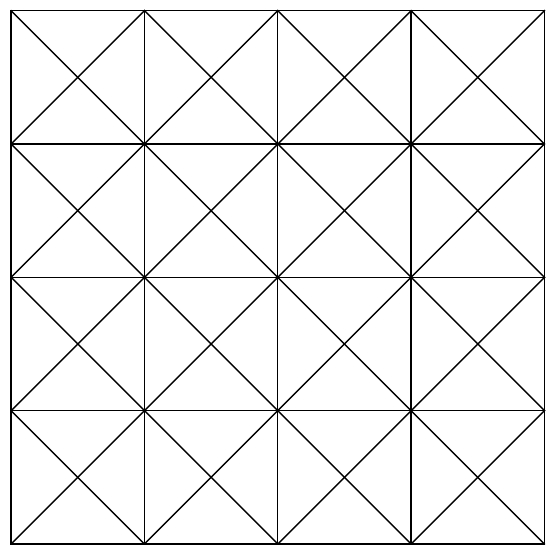}
}
\caption{Three different types of box spline meshes $\Delta(\bXi)$.}
\label{fig:grids}
\end{center}
\end{figure}

\begin{figure}[t!]
\begin{center}
\subfigure[cubic B-spline support (16 quads)]{
 \begin{overpic}[height=3.25cm]{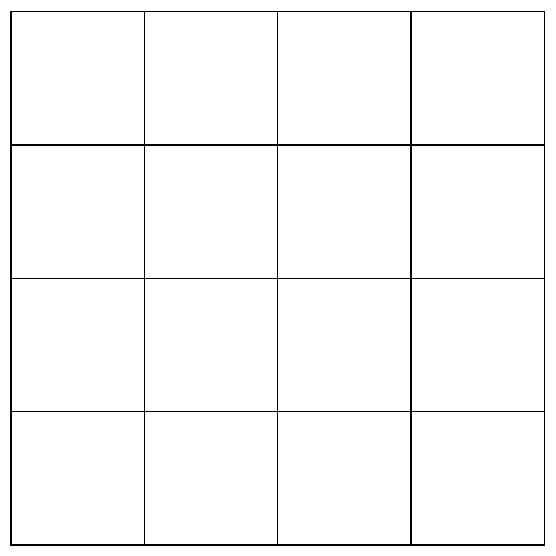}
 \put(41.8,41.5){\includegraphics[ height=.3cm]{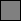}}
 \end{overpic}
} 
\subfigure[shape of cubic B-spline]{
 \includegraphics[height=3.25cm]{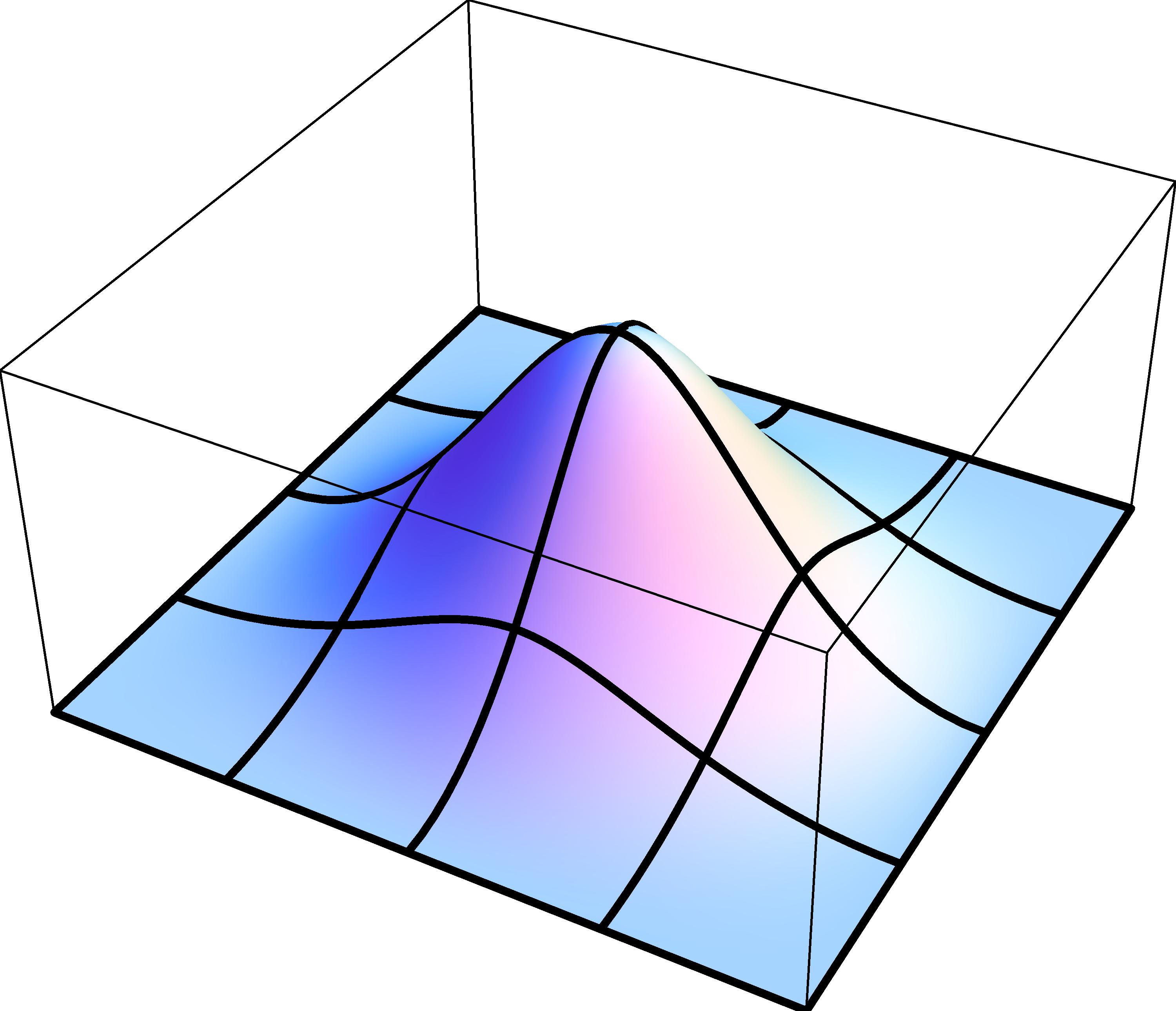}
}
\subfigure[quartic box spline support (24 triangles)]{
 \begin{overpic}[height=3.25cm]{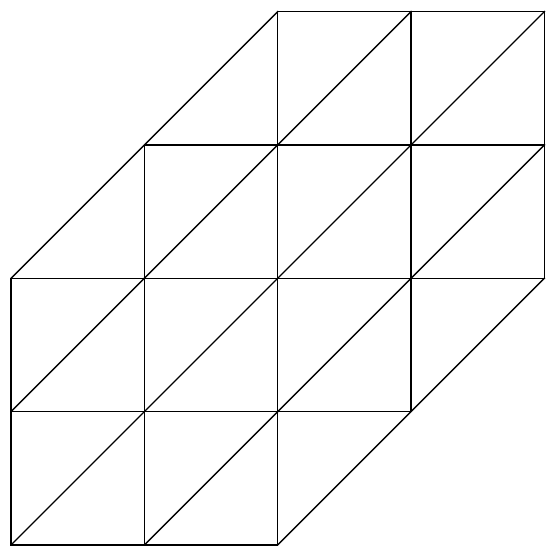}
  \put(41.8,41.5){\includegraphics[height=.3cm]{legend_a0_mod.png}}
 \end{overpic}
} 
\subfigure[shape of quartic box spline]{
 \includegraphics[height=3.25cm]{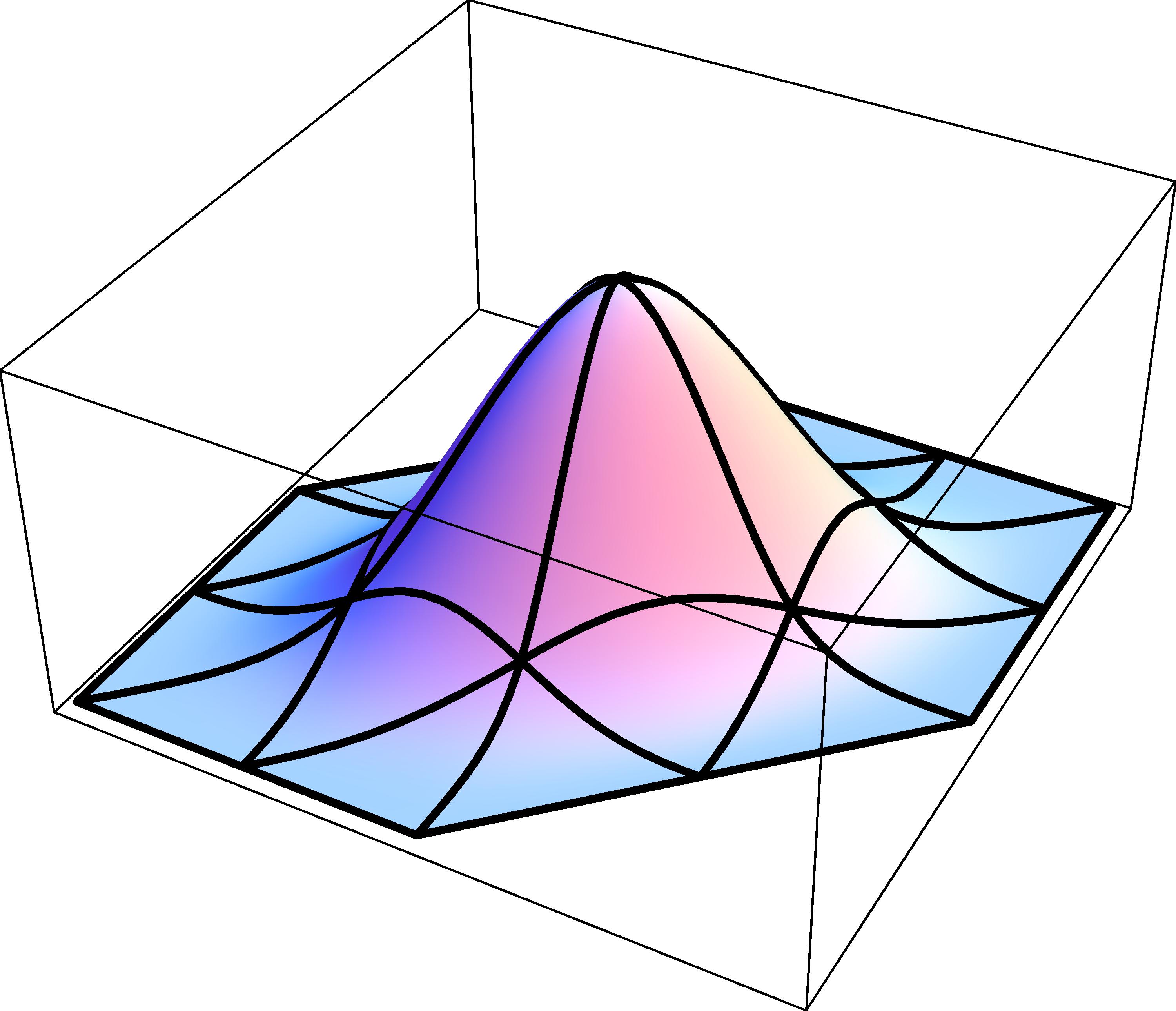}
}
\caption{$C^2$ cubic B-spline on a two-directional mesh, and $C^2$ quartic box spline on a three-directional mesh. The anchors of both splines are represented as {squares} in the center of the supports.}
\label{fig:C2quartic}
\end{center}
\end{figure}

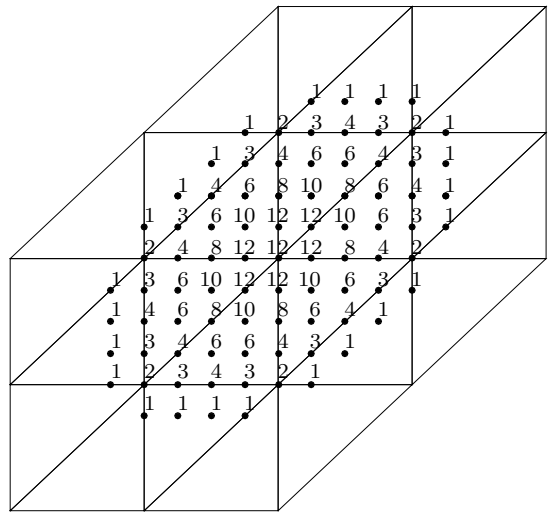
\begin{figure}[t!]
\begin{center}
\begin{overpic}[trim= 1.3cm .9cm .9cm .7cm, clip=true, height=7cm]{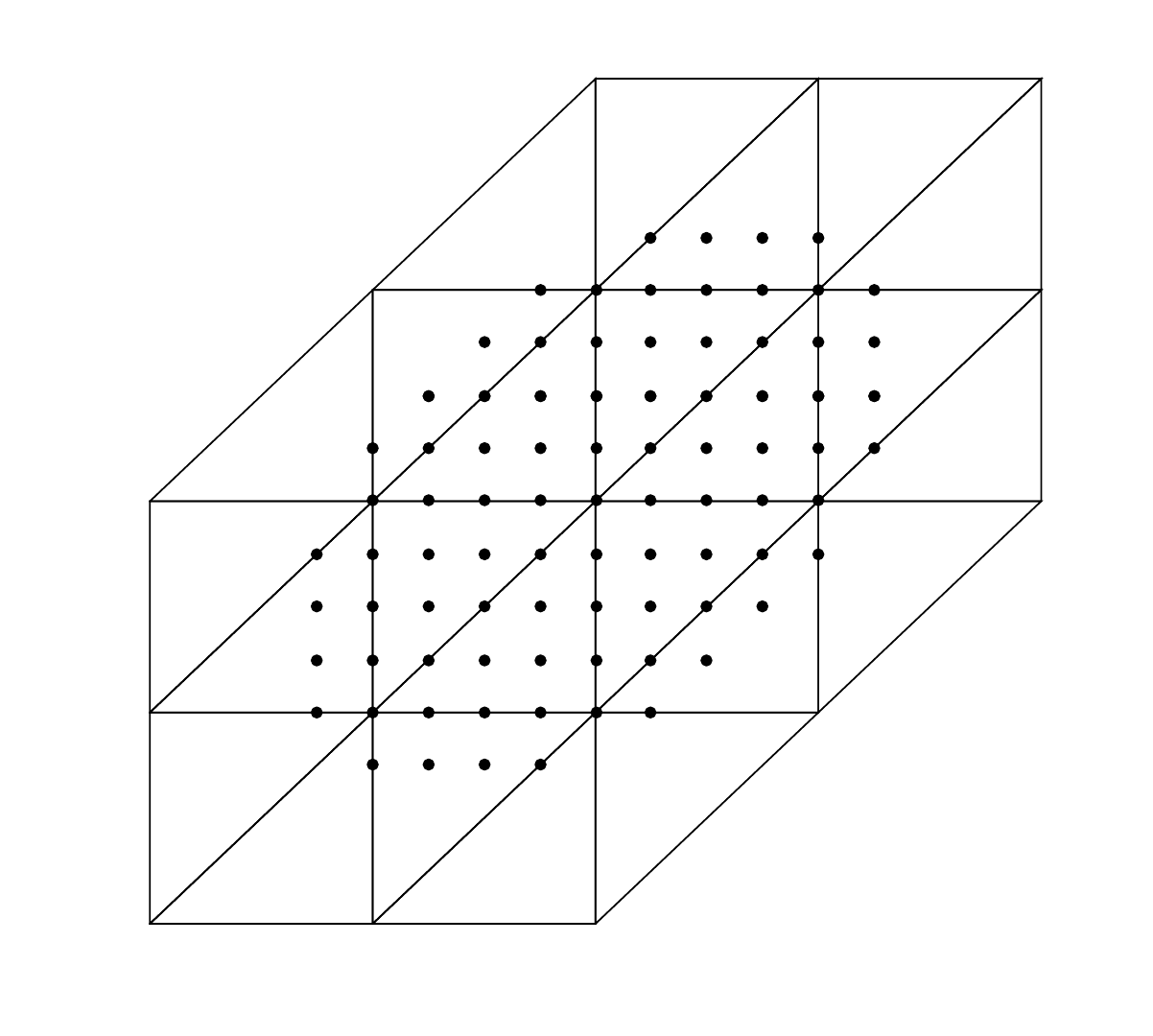}
\put(52,44){\footnotesize$\phantom{0}1$}
\put(64.5,44){\footnotesize$\phantom{0}1$}
\put(77,44){\footnotesize$\phantom{0}1$}
\put(89.5,44){\footnotesize$\phantom{0}1$}

\put(39.5,56){\footnotesize$\phantom{0}1$}
\put(52,56){\footnotesize$\phantom{0}2$}
\put(64.5,56){\footnotesize$\phantom{0}3$}
\put(77,56){\footnotesize$\phantom{0}4$}
\put(89.5,56){\footnotesize$\phantom{0}3$}
\put(102,56){\footnotesize$\phantom{0}2$}
\put(114,56){\footnotesize$\phantom{0}1$}

\put(39.5,67.5){\footnotesize$\phantom{0}1$}
\put(52,67.5){\footnotesize$\phantom{0}3$}
\put(64.5,67.5){\footnotesize$\phantom{0}4$}
\put(77,67.5){\footnotesize$\phantom{0}6$}
\put(89.5,67.5){\footnotesize$\phantom{0}6$}
\put(102,67.5){\footnotesize$\phantom{0}4$}
\put(114,67.5){\footnotesize$\phantom{0}3$}
\put(127,67.5){\footnotesize$\phantom{0}1$}

\put(39.5,79.5){\footnotesize$\phantom{0}1$}
\put(52,79.5){\footnotesize$\phantom{0}4$}
\put(64.5,79.5){\footnotesize$\phantom{0}6$}
\put(77,79.5){\footnotesize$\phantom{0}8$}
\put(89.5,79.5){\footnotesize$10$}
\put(102,79.5){\footnotesize$\phantom{0}8$}
\put(114,79.5){\footnotesize$\phantom{0}6$}
\put(127,79.5){\footnotesize$\phantom{0}4$}
\put(139.5,79.5){\footnotesize$\phantom{0}1$}

\put(39.5,91){\footnotesize$\phantom{0}1$}
\put(52,91){\footnotesize$\phantom{0}3$}
\put(64.5,91){\footnotesize$\phantom{0}6$}
\put(77,91){\footnotesize$10$}
\put(89.5,91){\footnotesize$12$}
\put(102,91){\footnotesize$12$}
\put(114,91){\footnotesize$10$}
\put(127,91){\footnotesize$\phantom{0}6$}
\put(139.5,91){\footnotesize$\phantom{0}3$}
\put(152,91){\footnotesize$\phantom{0}1$}

\put(52,103){\footnotesize$\phantom{0}2$}
\put(64.5,103){\footnotesize$\phantom{0}4$}
\put(77,103){\footnotesize$\phantom{0}8$}
\put(89.5,103){\footnotesize$12$}
\put(102,103){\footnotesize$12$}
\put(114.5,103){\footnotesize$12$}
\put(127,103){\footnotesize$\phantom{0}8$}
\put(139.5,103){\footnotesize$\phantom{0}4$}
\put(152,103){\footnotesize$\phantom{0}2$}

\put(52,115){\footnotesize$\phantom{0}1$}
\put(64.5,115){\footnotesize$\phantom{0}3$}
\put(77,115){\footnotesize$\phantom{0}6$}
\put(89.5,115){\footnotesize$10$}
\put(102,115){\footnotesize$12$}
\put(114.5,115){\footnotesize$12$}
\put(127,115){\footnotesize$10$}
\put(139.5,115){\footnotesize$\phantom{0}6$}
\put(152,115){\footnotesize$\phantom{0}3$}
\put(164.5,115){\footnotesize$\phantom{0}1$}

\put(64.5,126.5){\footnotesize$\phantom{0}1$}
\put(77,126.5){\footnotesize$\phantom{0}4$}
\put(89.5,126.5){\footnotesize$\phantom{0}6$}
\put(102,126.5){\footnotesize$\phantom{0}8$}
\put(114.5,126.5){\footnotesize$10$}
\put(127,126.5){\footnotesize$\phantom{0}8$}
\put(139.5,126.5){\footnotesize$\phantom{0}6$}
\put(152,126.5){\footnotesize$\phantom{0}4$}
\put(164.5,126.5){\footnotesize$\phantom{0}1$}

\put(77,138.5){\footnotesize$\phantom{0}1$}
\put(89.5,138.5){\footnotesize$\phantom{0}3$}
\put(102,138.5){\footnotesize$\phantom{0}4$}
\put(114.5,138.5){\footnotesize$\phantom{0}6$}
\put(127,138.5){\footnotesize$\phantom{0}6$}
\put(139.5,138.5){\footnotesize$\phantom{0}4$}
\put(152,138.5){\footnotesize$\phantom{0}3$}
\put(164.5,138.5){\footnotesize$\phantom{0}1$}

\put(89.5,150){\footnotesize$\phantom{0}1$}
\put(102,150){\footnotesize$\phantom{0}2$}
\put(114.5,150){\footnotesize$\phantom{0}3$}
\put(127,150){\footnotesize$\phantom{0}4$}
\put(139.5,150){\footnotesize$\phantom{0}3$}
\put(152,150){\footnotesize$\phantom{0}2$}
\put(164.5,150){\footnotesize$\phantom{0}1$}

\put(114.5,162){\footnotesize$\phantom{0}1$}
\put(127,162){\footnotesize$\phantom{0}1$}
\put(139.5,162){\footnotesize$\phantom{0}1$}
\put(152,162){\footnotesize$\phantom{0}1$}
 \end{overpic}

\caption{{Schematic representation of the Bernstein--B\'ezier form of the three-directional $C^2$ quartic box spline. Only the non-zero B\'ezier coefficients are shown and are multiplied by 24 for a better visualization.}}
\label{fig:bezierCoeffs}
\end{center}
\end{figure}

The construction of the box spline $M_\bXi$ does not depend on the ordering of the direction vectors $\bv_i$, $i=1,\ldots,n$. Moreover, it has the following properties (see, e.g., \cite{deboor1993}):
\begin{itemize}
 \item it is non-negative and its support is given by  $\Xi\, [0,1]^n$;
 \item it is $\rho-2$ times continuously differentiable, where $\rho$ is the minimal number of columns that need to be removed from $\bXi$ to obtain a matrix whose columns do not span $\RR^d$;
 \item it is a $d$-variate piecewise polynomial of total degree $n-d$ over the mesh
  \begin{equation} \label{eq:mesh}
    \Delta(\bXi) := \mathbb{H}(\bXi)+\mathbb{Z}^d,
  \end{equation}
  where $\mathbb{H}(\bXi)$ is the collection of all hyperplanes spanned by columns of any submatrix $\bX$ of $\bXi$ with $\rank(\bX)=d-1$.
\end{itemize}
Box splines and their derivatives can be evaluated through simple recurrence relations,
and elegant expressions are known for their inner products \cite{Speleers2015}.
There is a well-established theory for spaces spanned by integer translates of box splines, namely
\begin{equation} \label{eq:set-box}
\left\{M_\bXi(\cdot-\bi),\ \bi\in \ZZ^d \right\}.
\end{equation}
The elements in the set (\ref{eq:set-box}) have the following properties (see, e.g., \cite{deboor1993}):
\begin{itemize}
  \item they form a partition of unity;
  \item they are linearly independent if and only if
  \begin{equation}\label{eq:li}
  \text{$\det(\bX)\in \{-1,0,1\}$ for each $d\times d$ submatrix $\bX$ of $\bXi$;}
  \end{equation}
 \item local linear independence is equivalent to (global) linear independence;
  \item they reproduce polynomials of total degree $\rho-1$.
\end{itemize}
In the bivariate case, the most popular choices of box splines are defined on two-directional (tensor-product) meshes, three-directional (type I) meshes, and four-directional (type II) meshes, see Figure~\ref{fig:grids}. Box splines on two-directional meshes are nothing else than uniform tensor-product B-splines of maximal smoothness.
For example, $C^2$ cubic tensor-product B-splines can be seen as scaled integer translates of the box spline $M_\bXi$ generated by the matrix
\begin{equation}
  \bXi= \left[
\begin{array}{cccccccc}
1 & 0 & 1 & 0 & 1 & 0 & 1 & 0\\
0 & 1 & 0 & 1 & 0 & 1 & 0 & 1
\end{array}\right],\nonumber
\end{equation}
see Figure~\ref{fig:C2quartic}(a)--(b). The (local) linear
independence condition (\ref{eq:li}) is satisfied on two- and
three-directional meshes, but not on four-directional meshes. This
property is required for the construction of adaptive box spline
bases using the hierarchical approach \cite{giannelli2014}. In our
numerical experiments (see Section~\ref{sec:exm}) we will solely
focus on $C^2$ quartic splines defined over three-directional
meshes ($d=2$), i.e., scaled integer translates of the box spline
$M_\bXi$ generated by the matrix
\begin{equation}\label{eq:C2-quartic}  \bXi= \left[
\begin{array}{cccccc}
1 & 0 & 1 & 1 & 0 & 1\\
0 & 1 & 1 & 0 & 1 & 1
\end{array}\right],
\end{equation}
see Figure~\ref{fig:C2quartic}(c)--(d). {A common
way to manipulate box splines on a three-directional mesh is using
their Bernstein--B\'ezier form on the triangles in their support.
Figure~\ref{fig:bezierCoeffs} shows the non-zero coefficients of
the $C^2$ quartic box spline generated by \eqref{eq:C2-quartic} in
its Bernstein--B\'ezier form. On each of the 24 triangles in its
support, the box spline is a quartic bivariate polynomial which
can be represented in terms of the Bernstein polynomials
\begin{equation}
B_{ijk}^4(u,v,w):=\frac{4!}{i!j!k!}u^i v^j w^k,\nonumber
\end{equation}
for $i+j+k=4$ and $i,j,k \in \NN_0$, and $u,v,w$ are barycentric
coordinates defined on the corresponding triangle.} More details
on the computation of $C^2$ quartic box splines can be found
in~\cite{pgmss}.

\subsection{Hierarchical basis constructions}
Let $\hat{\Omega}^0$ be a given domain in $\RR^d$. We consider a nested sequence of $d$-variate spline spaces  defined on $\hat{\Omega}^0$,
\begin{equation} \label{eq:spaces}
\hat{\VV}^0\subset \hat{\VV}^1 \subset \hat{\VV}^2\subset \cdots,
\end{equation}
where any $\hat{\VV}^\ell$ is a space spanned by scaled integer translates of a certain box spline $M_\bXi$ using a proper scaling factor $h_\ell$ according to the level $\ell$ (we assume $h_\ell > h_{\ell+1}$), i.e.,
\begin{equation} \label{eq:set-box-ell}
 \hat{{\cal B}}^\ell := \left\{M_\bXi\Bigl(\frac{\cdot}{h_\ell}-\bi\Bigr),\ \bi\in \ZZ^d \right\}.
\end{equation}
For simplicity of notation, we denote an element of $\hat{{\cal B}}^\ell$ by $\hat{\beta}^\ell$.
In addition, we consider a nested sequence of subsets of $\hat{\Omega}^0$,
\begin{equation*}
\hat{\Omega}^0\supseteq \hat{\Omega}^1 \supseteq \hat{\Omega}^2\supseteq \cdots,
\end{equation*}
where each $\hat{\Omega}^\ell$ is chosen to be aligned with the scaled mesh $\Delta^\ell(\bXi):=h_\ell\Delta(\bXi)$, recalling $\Delta(\bXi)$ from (\ref{eq:mesh}).
We assume that $\hat{\Omega}^N=\emptyset$ for some $N\in\NN$, and
we denote the corresponding finite sequence by
$\hat{\bOmega}:=\{\hat{\Omega}^{0},\hat{\Omega}^{1},\ldots,\hat{\Omega}^{N-1}\}$.
This sequence represents the regions to be refined at different levels of resolution.

\begin{figure}[t!]
\begin{center}
\subfigure[uniform mesh (1 level)]{
\includegraphics[trim= 3.25cm 2.5cm 3.5cm 3cm, clip=true, height=4cm]{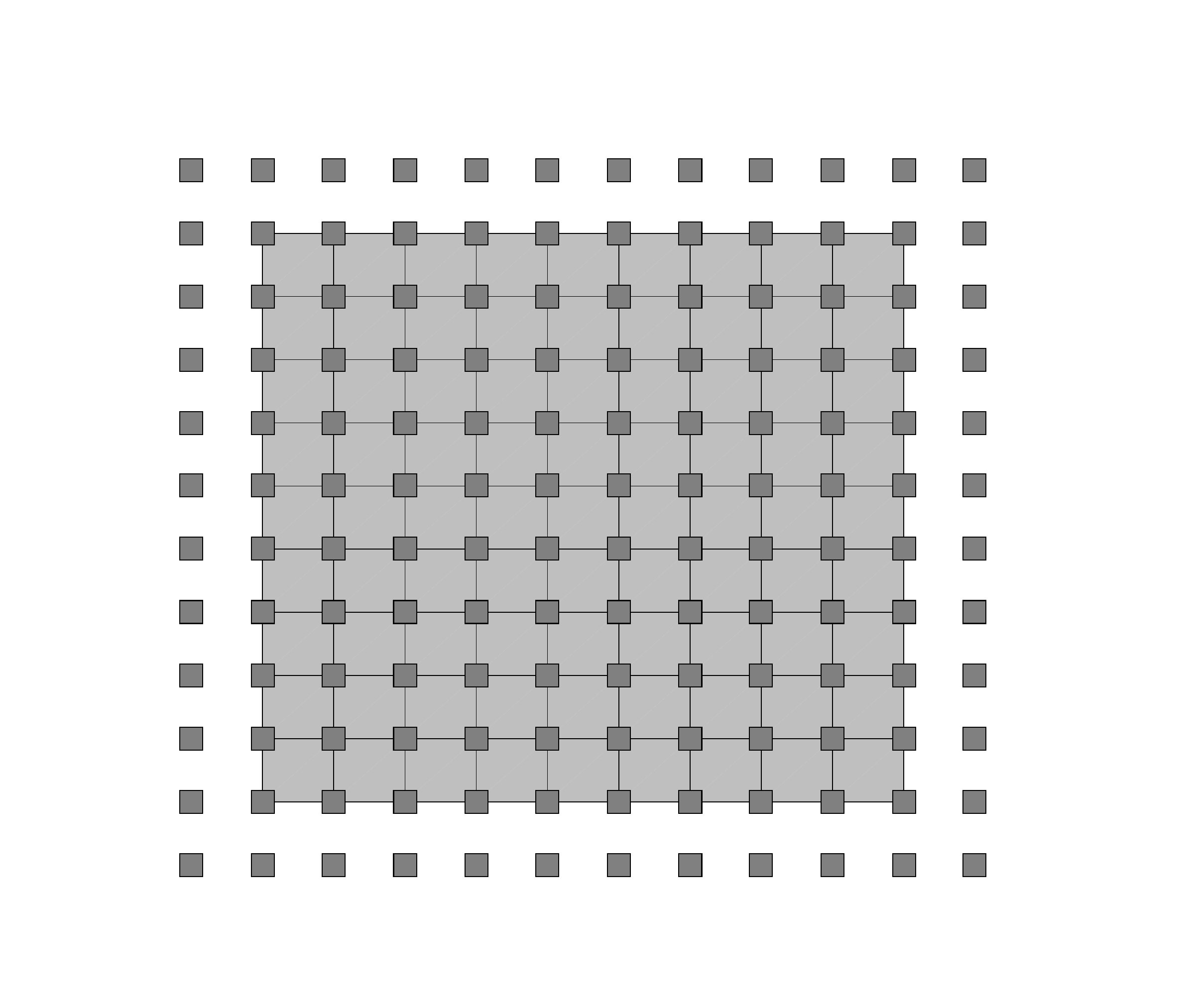}
}
\subfigure[hierarchical mesh with 2 levels]{
\includegraphics[trim= 3.25cm 2.5cm 3.5cm 3cm, clip=true, height=4cm]{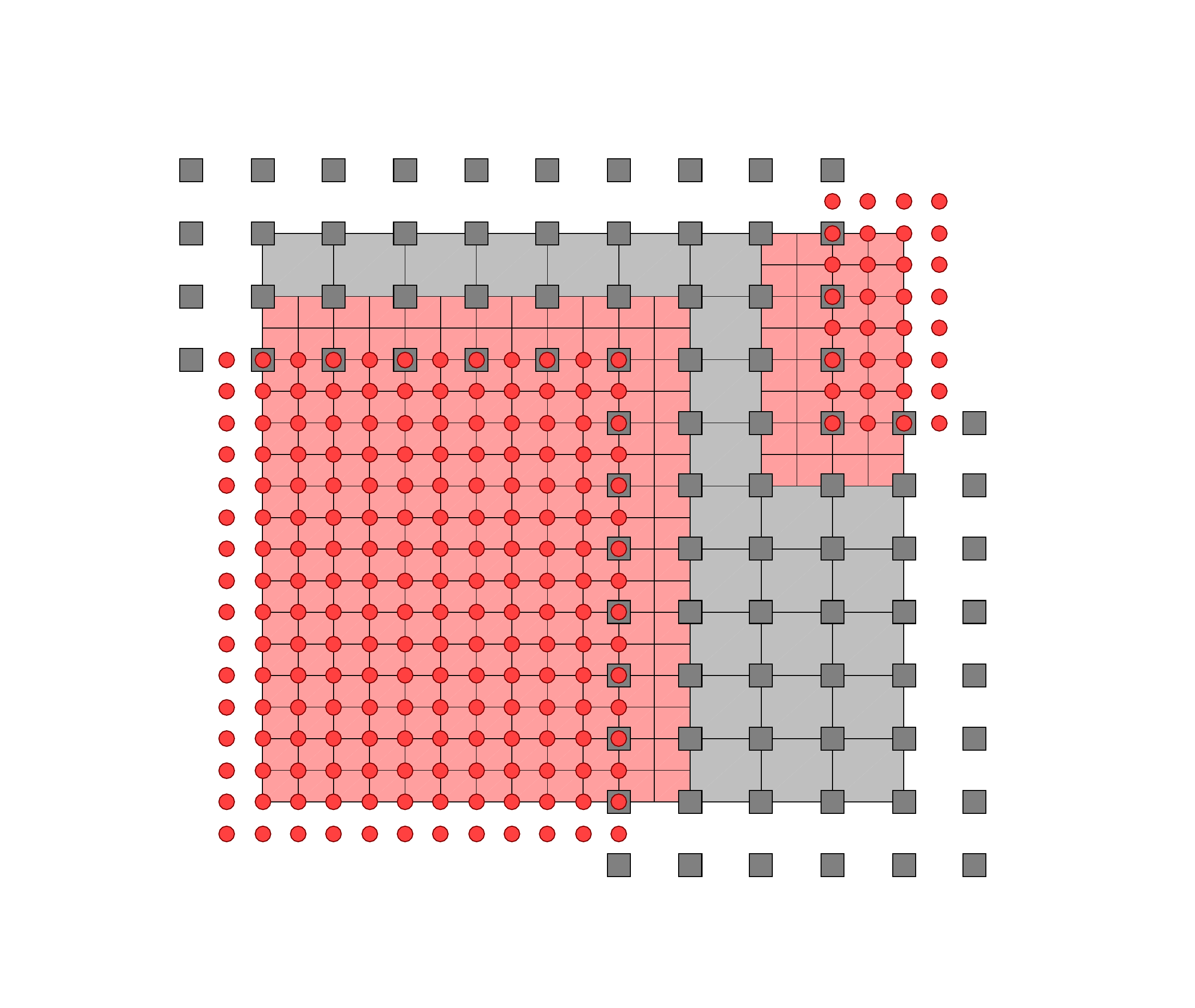}
}
\subfigure[hierarchical mesh with 3 levels]{
\includegraphics[trim= 3.25cm 2.5cm 3.5cm 3cm, clip=true, height=4cm]{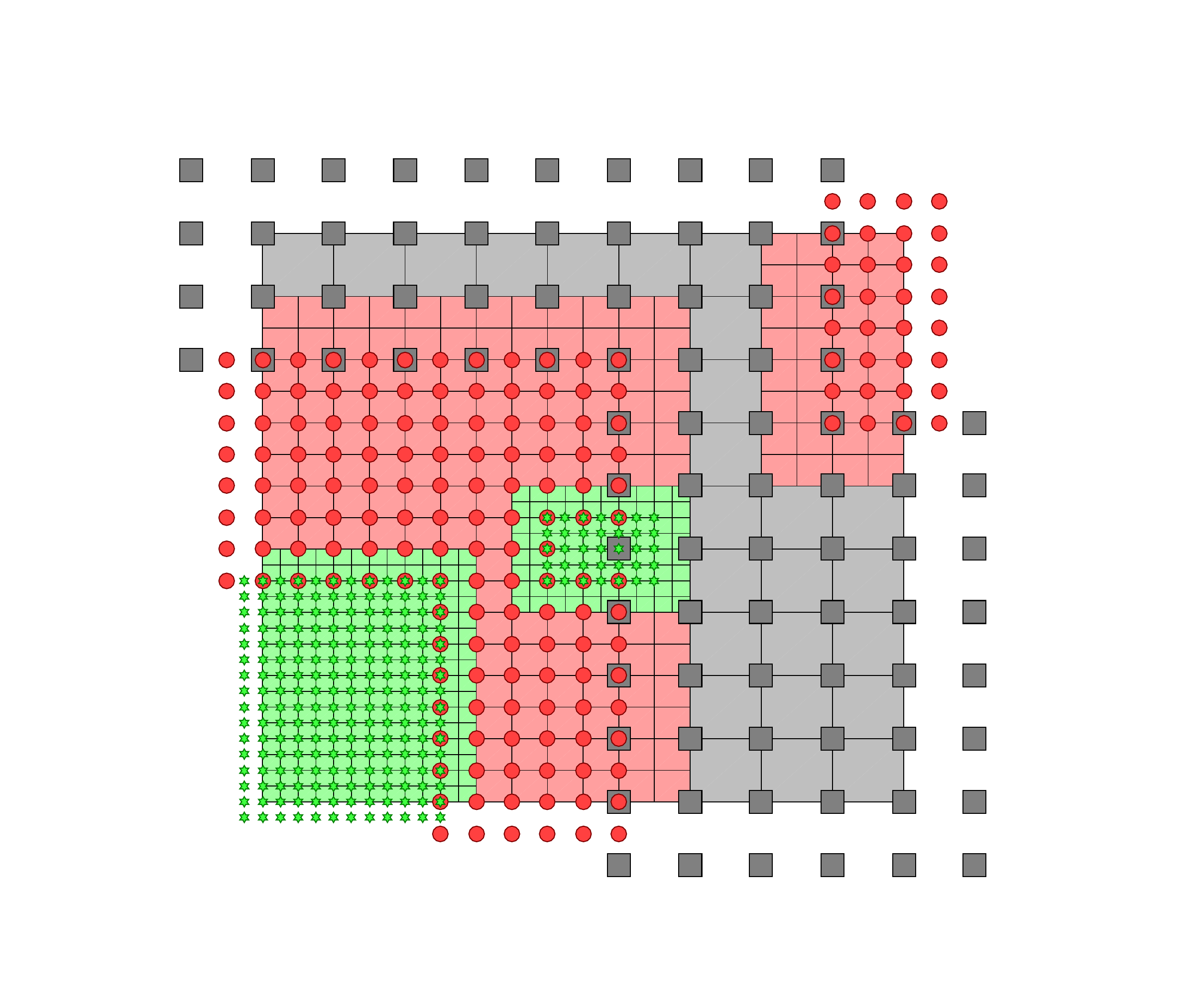}
}
\caption{A sequence of locally refined hierarchical tensor-product meshes.
The anchors of $C^2$ cubic (T)HB-splines are indicated for the different levels.}
\label{fig:hmesh-tp}
\end{center}
\end{figure}

\begin{figure}[t!]
\begin{center}
\subfigure[uniform mesh (1 level)]{
\includegraphics[trim= 3.25cm 2.5cm 3.5cm 3cm, clip=true, height=4cm]{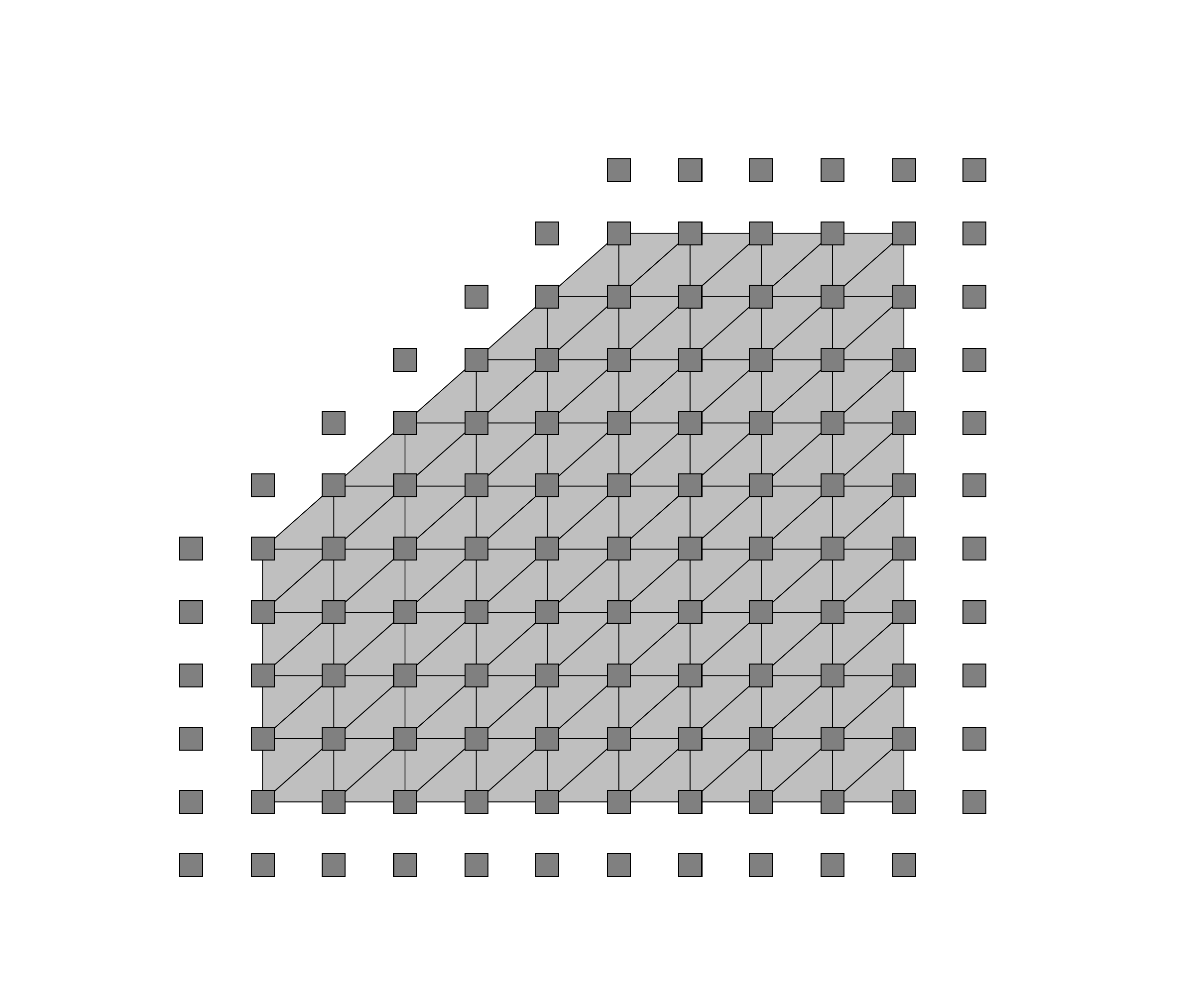}
}
\subfigure[hierarchical mesh with 2 levels]{
\includegraphics[trim= 3.25cm 2.5cm 3.5cm 3cm, clip=true, height=4cm]{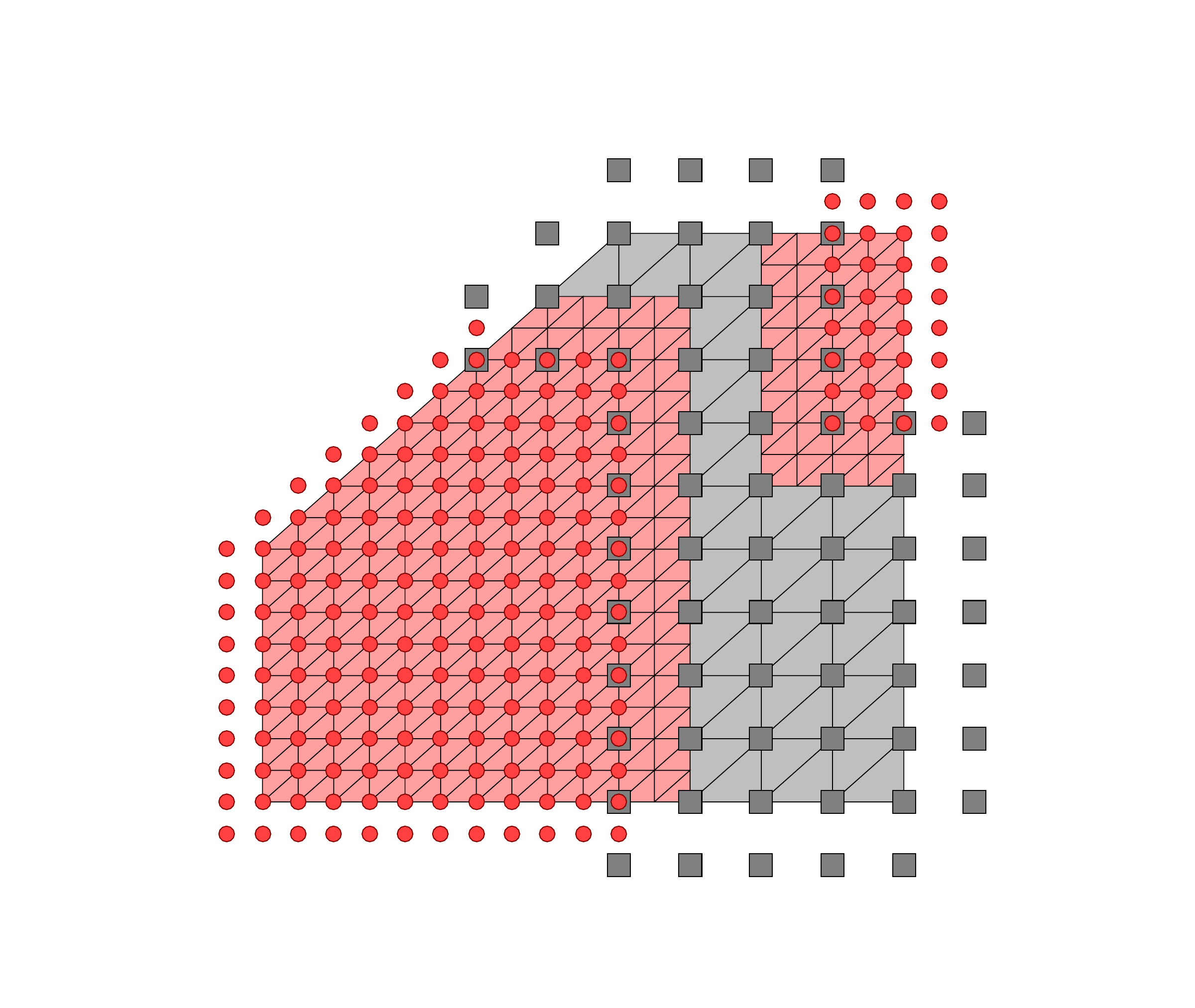}
}
\subfigure[hierarchical mesh with 3 levels]{
\includegraphics[trim= 3.25cm 2.5cm 3.5cm 3cm, clip=true, height=4cm]{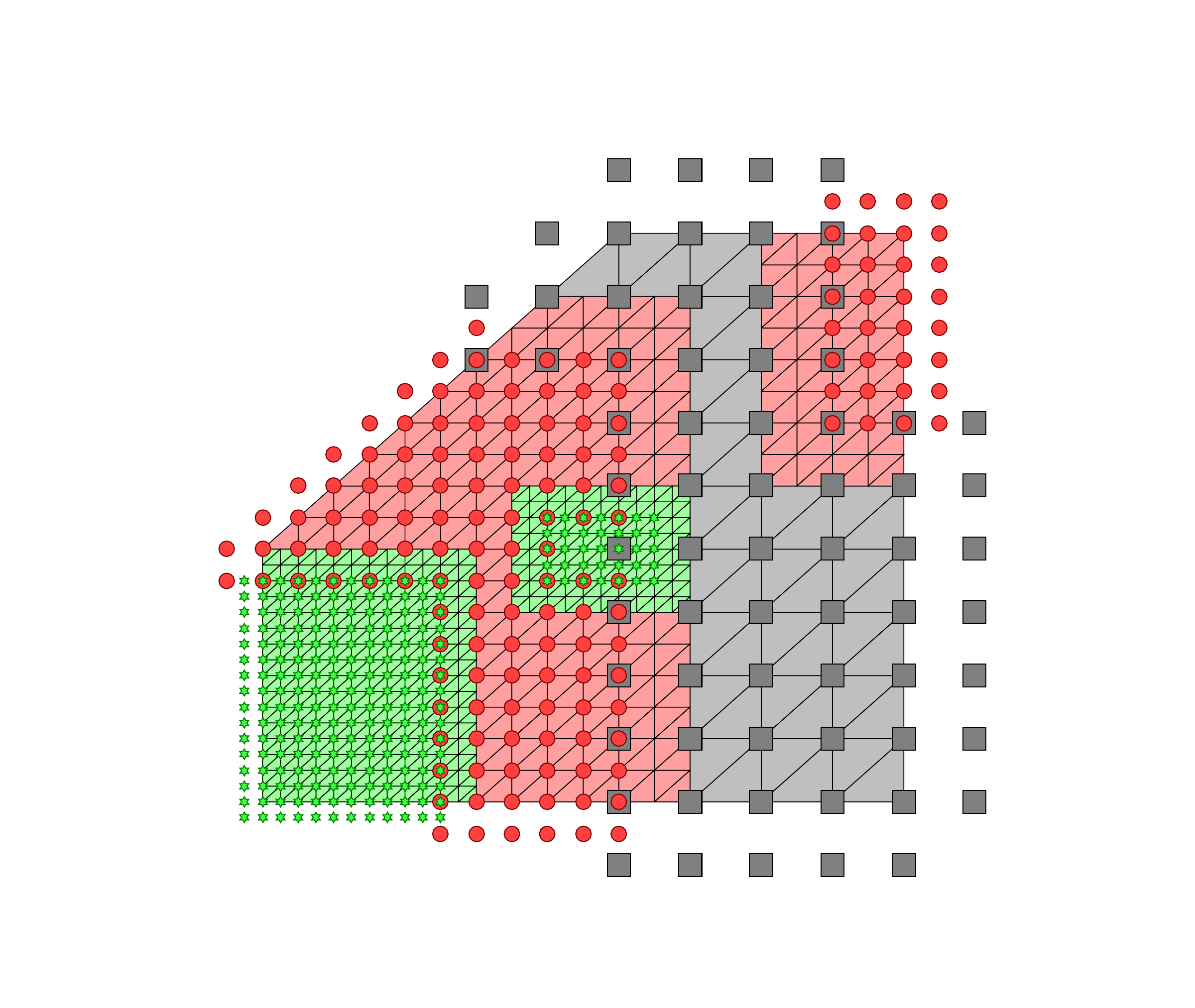}
}
\caption{A sequence of locally refined hierarchical three-directional meshes.
The anchors of $C^2$ quartic (T)HBox-splines are indicated for the different levels.}
\label{fig:hmesh-box}
\end{center}
\end{figure}

Throughout the paper, we graphically represent mesh cells and basis functions at different levels as shown in Table~\ref{tab:legend}.
An example of a hierarchical refinement for a two- and three-directional mesh is shown in Figures~\ref{fig:hmesh-tp} and \ref{fig:hmesh-box}, respectively.

%
%

\begin{table}[t!]
\small
\begin{center}
\begin{tabular}{l|ccccc}
 \rule{0cm}{.6cm} & level 0 & level 1 & level 2 & level 3 & level 4\\
\hline
cell \rule{0cm}{.6cm}
& \includegraphics[trim= .5cm .5cm .2cm .2cm, clip=true, height=.5cm]{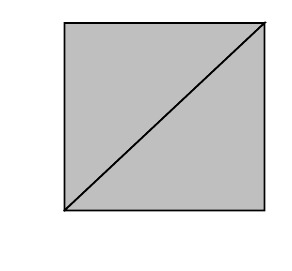}
& \includegraphics[trim= .5cm .5cm .2cm .2cm, clip=true, height=.5cm]{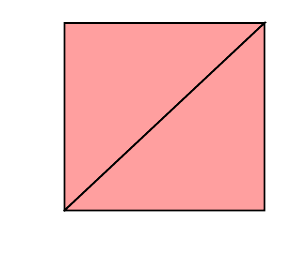}
& \includegraphics[trim= .5cm .5cm .2cm .2cm, clip=true, height=.5cm]{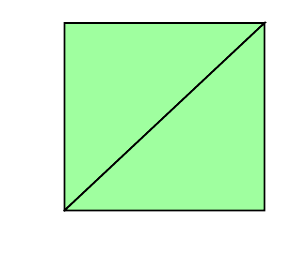}
& \includegraphics[trim= .5cm .5cm .2cm .2cm, clip=true, height=.5cm]{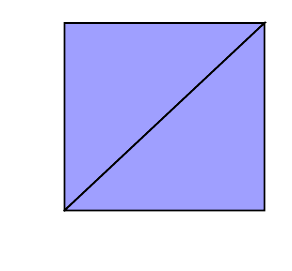}
& \includegraphics[trim= .5cm .5cm .2cm .2cm, clip=true, height=.5cm]{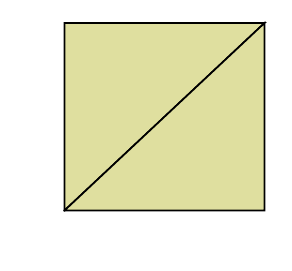}\\
anchor \rule{0cm}{.6cm}
& \includegraphics[trim= 2.63cm 2.55cm 2.37cm 2.35cm, clip=true, height=.5cm]{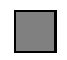}
& \includegraphics[trim= 2.63cm 2.55cm 2.37cm 2.35cm, clip=true, height=.5cm]{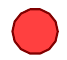}
& \includegraphics[trim= 2.63cm 2.55cm 2.37cm 2.35cm, clip=true, height=.5cm]{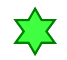}
& \includegraphics[trim= 1.5cm .9cm .8cm 1cm, clip=true, height=.5cm]{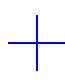}
& \includegraphics[trim= 1.5cm .9cm .8cm 1cm, clip=true, height=.5cm]{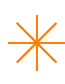}
\end{tabular}
\caption{Legend for graphical elements representing mesh cells and anchors (= centers) of basis functions at different hierarchical levels.}
\label{tab:legend}
\end{center}
\end{table}

\begin{figure}[t!]
\begin{center}
\subfigure[{box spline before refinement}]{
\begin{overpic}[trim= 1.3cm .9cm .9cm .7cm, clip=true, height=4cm]{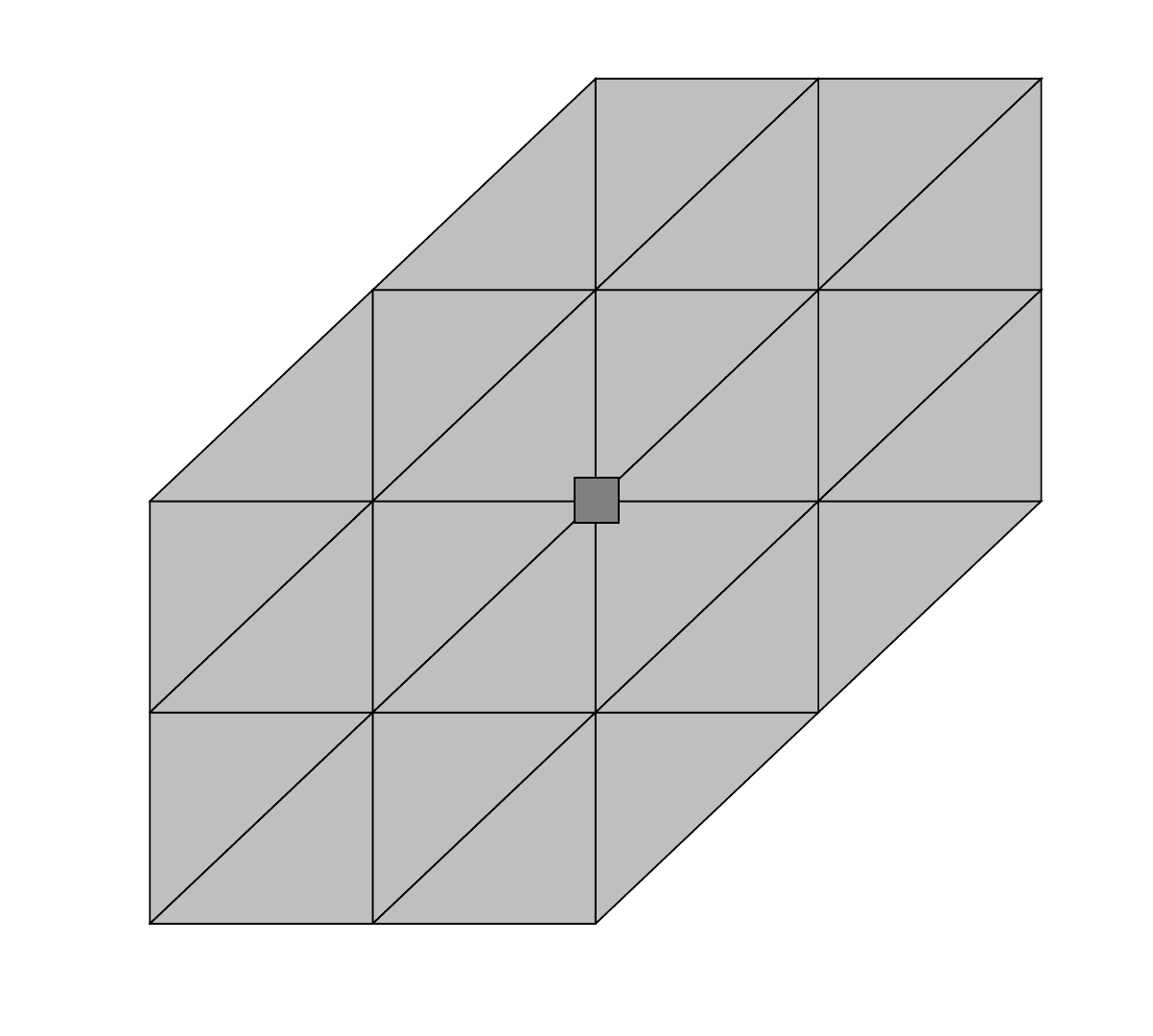}
\put(51,62){\footnotesize$16$}
 \end{overpic}
}\hspace{2cm} \subfigure[{refinement of the box
spline}]{
\begin{overpic}[trim= 1.3cm .9cm .9cm .7cm, clip=true, height=4cm]{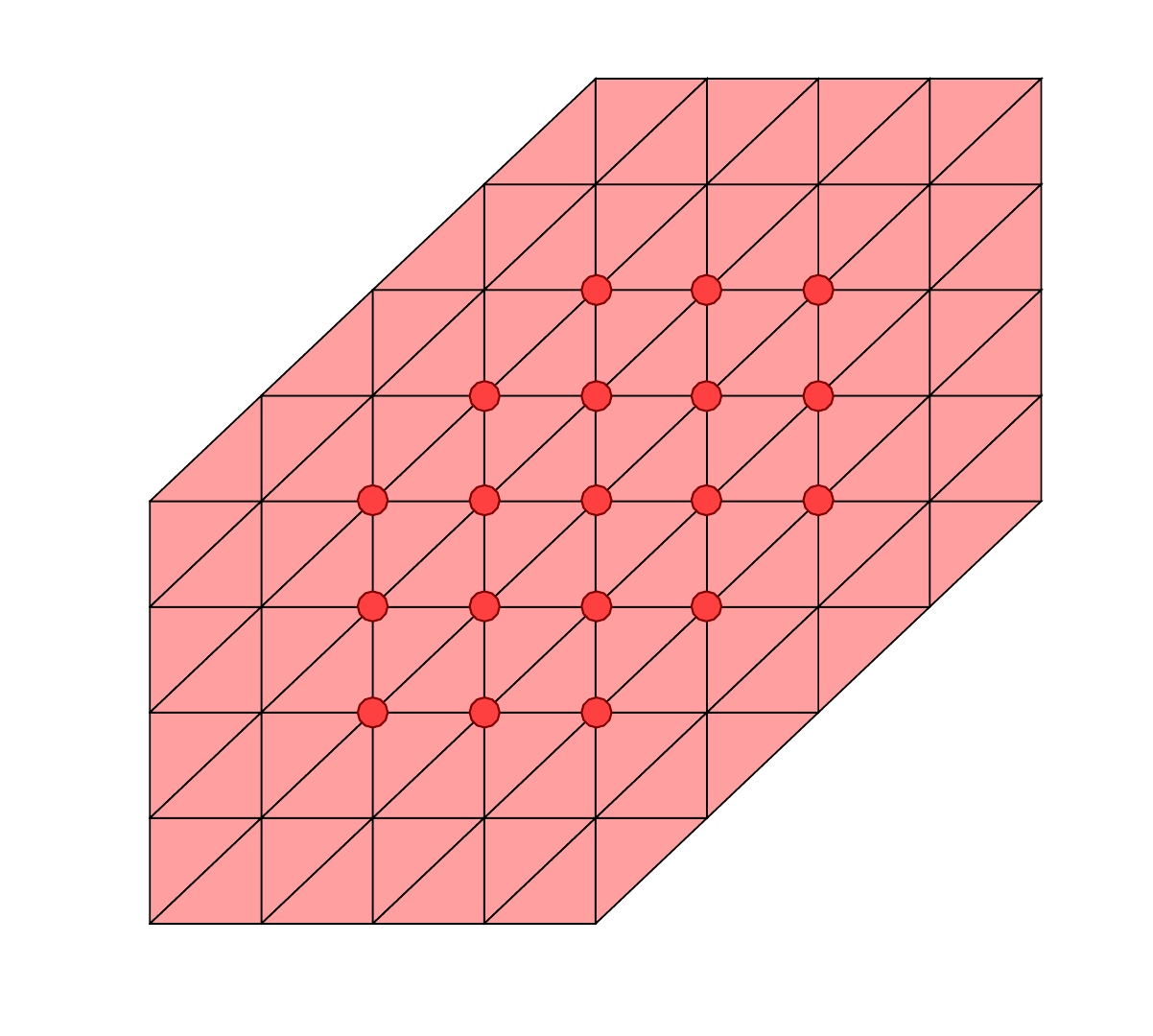}
\put(25.5,32.5){\footnotesize$1$}
\put(39.75,32.5){\footnotesize$2$}
\put(54,32.5){\footnotesize$1$}

\put(25.5,46.25){\footnotesize$2$}
\put(39.75,46.25){\footnotesize$6$}
\put(54,46.25){\footnotesize$6$}
\put(68.25,46.25){\footnotesize$2$}

\put(25.5,60){\footnotesize$1$}
\put(39.75,60){\footnotesize$6$}
\put(54,60){\footnotesize$10$}
\put(68.25,60){\footnotesize$6$}
\put(82.5,60){\footnotesize$1$}

\put(39.75,73.5){\footnotesize$2$}
\put(54,73.5){\footnotesize$6$}
\put(68.25,73.5){\footnotesize$6$}
\put(82.5,73.5){\footnotesize$2$}

\put(54,87){\footnotesize$1$}
\put(68.25,87){\footnotesize$2$}
\put(82.5,87){\footnotesize$1$}
 \end{overpic}

}
\caption{{Schematic representation of the two-scale (dyadic) refinement relation of three-directional $C^2$ quartic box splines. Only the non-zero box spline coefficients are shown and are multiplied by 16 for a better visualization.}}
\label{fig:boxSplineRefine}
\end{center}
\end{figure}
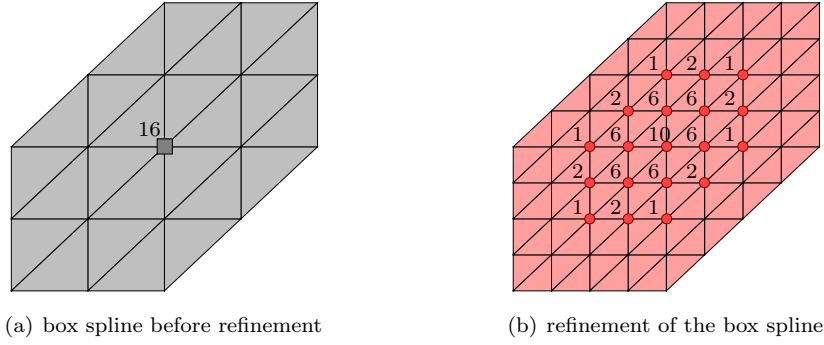

Any set of box splines $\hat{{\cal B}}^\ell$ forms a partition of unity, and its elements are non-negative with local compact support. Moreover, we assume that the chosen sequence of $\hat{{\cal B}}^\ell$, $\ell=0,\ldots,N-1$, possesses the following additional properties:
\begin{itemize}
\item each $\hat{{\cal B}}^\ell$ is locally linearly independent, so $\bXi$ satisfies (\ref{eq:li});
\item two-scale relations between $\hat{{\cal B}}^\ell$ and $\hat{{\cal B}}^{\ell+1}$ have only non-negative
coefficients, so
\begin{equation}\label{eq:2scale}
\hat{\beta}^\ell=\sum_{{\hat{\beta}}^{\ell+1}\in \hat{{\cal B}}^{\ell+1}}
c_{{\hat{\beta}}^{\ell+1}}(\hat{\beta}^\ell)\, \hat{\beta}^{\ell+1},
\quad c_{{\hat{\beta}}^{\ell+1}}(\hat{\beta}^\ell)\geq0.
\end{equation}
\end{itemize}
These properties are satisfied for {the popular families of} box splines defined over nested sequences of two- or three-directional meshes constructed by dyadic refinement.
{For example, the coefficients in the two-scale (dyadic) refinement relation of three-directional $C^2$ quartic box splines are shown in Figure~\ref{fig:boxSplineRefine}.}

According to the hierarchical approach \cite{giannelli2014}, we define the hierarchical box spline basis as follows.
\begin{definition}\label{dfn:hb}
The hierarchical box spline (HBox-spline) basis $\hat{\cH}(\hat\bOmega)$ related to the hierarchy $\hat\bOmega$
is defined as
\begin{equation*}
\hat{\cH}(\hat\bOmega) :=
\left\{
\hat{\beta}^\ell\in\hat{{\cal B}}^\ell :
{\supp}^{\hat \Omega^0}(\hat{\beta}^\ell) \subseteq\hat{\Omega}^\ell \wedge
{\supp}^{\hat \Omega^0}(\hat{\beta}^\ell) \not\subseteq \hat{\Omega}^{\ell+1},
\, \ell =0,\ldots,N-1
\right\},
\end{equation*}
where ${\supp}^{\hat \Omega^0}(\hat{\beta}^\ell)$ stands for the intersection of the support of $\hat{\beta}^\ell$ with $\hat{\Omega}^0$, i.e., the largest set in $\hat \bOmega$.
\end{definition}
From \cite{giannelli2014} we know that the aforementioned box spline properties ensure that the elements in $\hat{\cH}(\hat\bOmega)$ are non-negative and linearly independent.

By exploiting the two-scale relations (\ref{eq:2scale}), we define the truncation of a function $\hat{s}\in \hat{\VV}^\ell\subset \hat{\VV}^{\ell+1}$ expressed with respect to $\hat{\cal B}^{\ell+1}$
as the linear combination of only those basis functions in $\hat{\cal B}^{\ell+1}$ that are not included in $\hat{\cH}(\hat\bOmega)$,
\begin{equation}\label{eq:trunc}
{\trunc}^{\ell+1} \hat{s} := \sum_{\hat{\beta}^{\ell+1} \in
\hat{\cal B}^{\ell+1}, \,
\supp^{\hat\Omega^0}(\hat{\beta}^{\ell+1})\not\subseteq\hat{\Omega}^{\ell+1}}
c_{{\hat{\beta}}^{\ell+1}}(\hat{s}) \,\hat{\beta}^{\ell+1}.
\end{equation}
By iterating the truncation over the spline hierarchy, we define the truncated basis for hierarchical box splines as follows.
\begin{definition}\label{dfn:thb}
The truncated hierarchical box spline (THBox-spline) basis $\hat{\cT}(\hat\bOmega)$ related to the hierarchy $\hat\bOmega$ is defined as
\begin{equation*}
\hat{\cT}(\hat\bOmega) := 
\left\{
 {\Trunc}^{\ell+1}(\hat{\beta}^{\ell}):\hat{\beta}^{\ell}\in\hat{{\cal B}}^\ell
\cap\hat{\cH}(\hat\bOmega),\,  \ell =0,\ldots,N-1\right\},
\end{equation*}
where $ {\Trunc}^{\ell+1}(\hat{\beta}^{\ell}) := {\trunc}^{N-1}({\trunc}^{N-2}(\cdots ({\trunc}^{\ell+1}(\hat{\beta}^{\ell}))\cdots))$.
\end{definition}
The non-negativity and linear independence of HBox-splines are preserved by
the truncation mechanism. In addition, the truncated basis $\hat{\cT}(\hat\bOmega)$ spans the same space as $\hat{\cH}(\hat\bOmega)$ and forms a partition of unity, see~\cite{giannelli2014}. Recently, a general and very simple procedure has been developed for the construction of quasi-interpolants in such hierarchical spaces \cite{speleersM2016}.

Since uniform tensor-product splines are a special instance of box splines, uniform (T)HB-splines\footnote{(T)HB-splines can be defined on any kind of nested (non-uniform) knot sequences, see \cite{giannelli2012,giannelli2014}.} defined in \cite{giannelli2012} are a special instance of (T)HBox-splines.
The anchors (i.e., centers of untruncated supports) of $C^2$ cubic (T)HB-splines are shown in
Figure~\ref{fig:hmesh-tp} for a sequence of hierarchical tensor-product meshes, and the anchors of $C^2$ quartic (T)HBox-splines are depicted in Figure~\ref{fig:hmesh-box} for a similar sequence of hierarchical three-directional meshes.

\section{Isogeometric methods using weakly imposed boundary conditions} \label{sec:weak}
When Dirichlet boundary conditions are imposed explicitly in the classical Galerkin method, the basis functions considered in the discretization space need to satisfy certain requirements along the boundary. 
In our hierarchical box spline setting, this would require special basis functions at the boundary that diverge from standard box splines used in the interior of the domain. 
In order to exploit the full potential of the uniform structure of box splines on the entire domain, we prefer to impose Dirichlet boundary conditions in a weak sense in combination with a sufficiently smooth invertible geometry map $\bF:\hat \Omega \to \Omega$, so uniting features of the immersed boundary and isogeometric approach. Throughout the paper, all quantities related to the parametric domain $\hat \Omega$ will be denoted by symbols with a hat. We explain and illustrate the approach for $d=2$, but all the ideas extend to higher dimensionality.

\subsection{A weak boundary formulation} \label{sec-WBC-method}

We briefly summarize the weak boundary formulation proposed in \cite{baiges2012} which will be considered in our numerical examples.
We explain the method by means of the following advection-diffusion problem defined on the domain $\Omega\subset\RR^2$:
\begin{align}
\label{problem0}
\left\{
\begin{array}{rll}
  -\kappa \Delta u + \ba\cdot\nabla u &=\ f, & \hbox{in }\ \Omega, \\
 u &=\ g, & \hbox{on } \Gamma:=\partial \Omega, \\
\end{array}
\right.
\end{align}
where $0<\kappa\leq1$, $\ba$ is the advection velocity, $f$ is the given forcing function, and $g$ is the given Dirichlet boundary function. 

\begin{figure}[t!]
\begin{center}
\subfigure[$\Omega$ and $\Omega_h$]{
\begin{overpic}[trim= 2.5cm 1cm 2cm .5cm, clip=true, height=4cm]{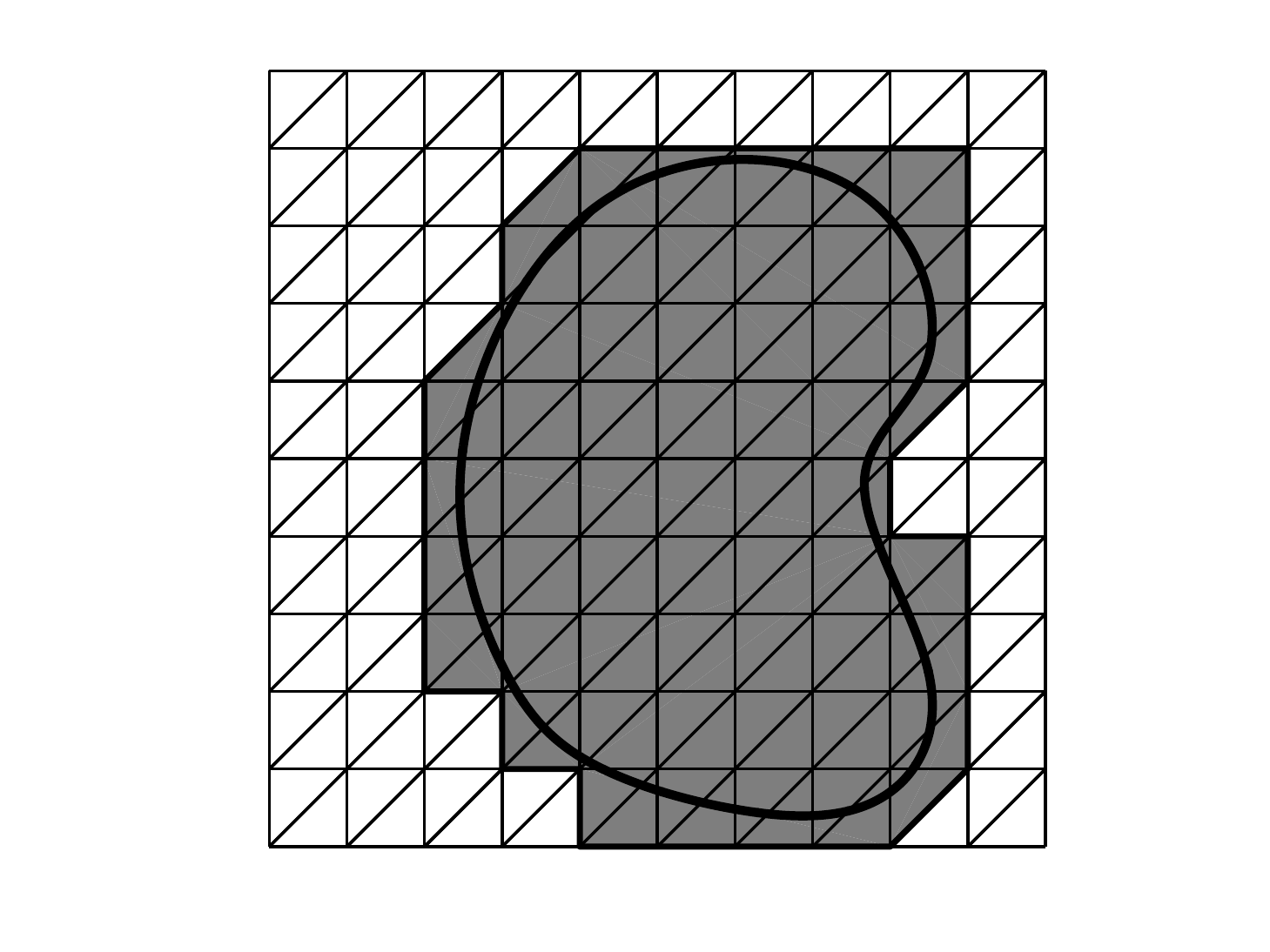}
 \put(52,47.3){{\textcolor{white}{\small$\Omega\subseteq\Omega_h$}}}
 \end{overpic}
 }
\subfigure[$\Omega_\Gamma$ and $\Omega_\textrm{in}$]{
\begin{overpic}[trim= 2.5cm 1cm 2cm .5cm, clip=true, height=4cm]{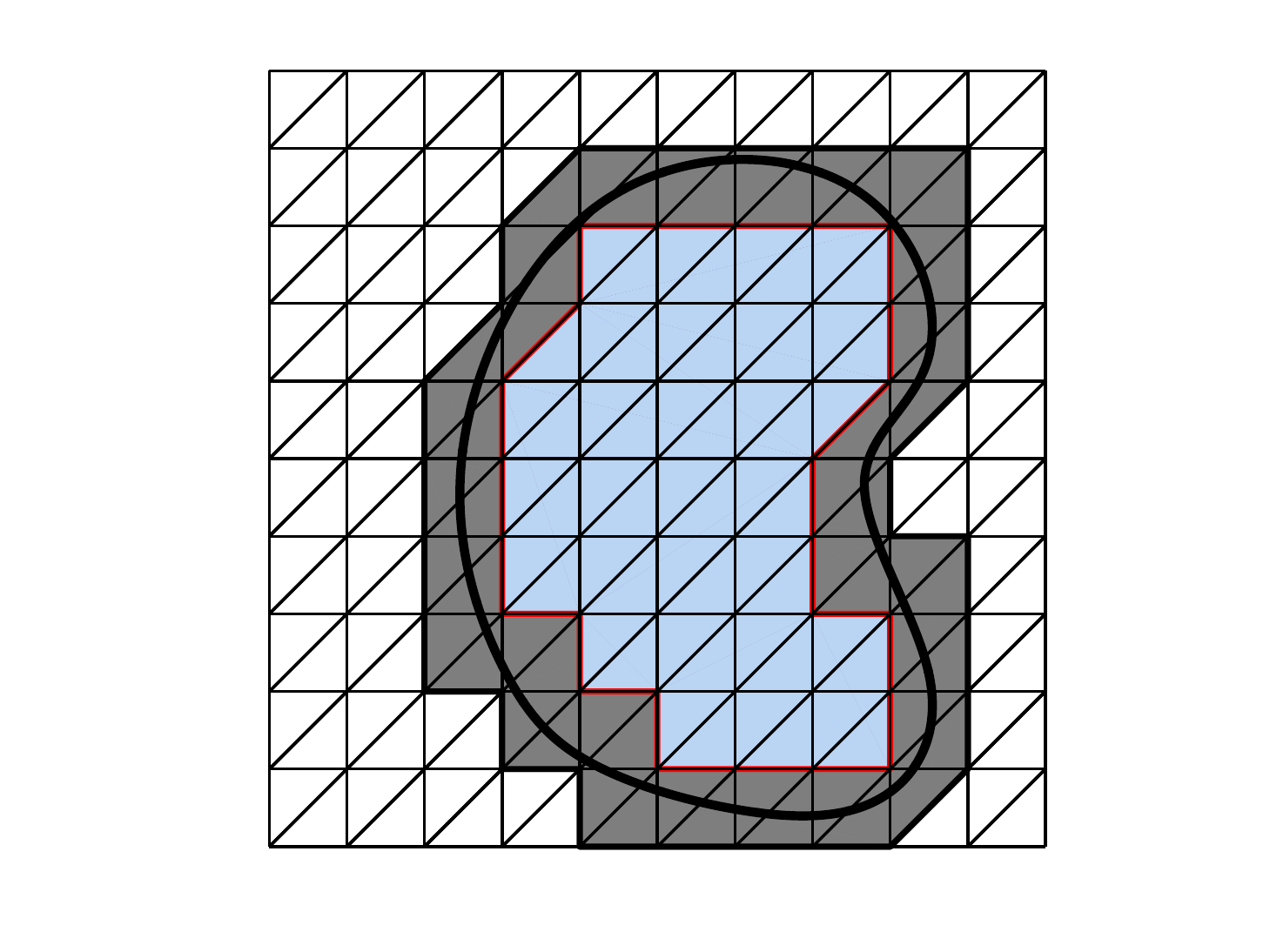}
\put(52,47.3){{\textcolor{black}{\small$\Omega_{\textrm{in}}$}}}
 \put(73,90.3){{\textcolor{white}{\small$\Omega_\Gamma$}}}
 \end{overpic}
}
\subfigure[$\Omega_{\Gamma, \textrm{in}}$ and $\Omega_{\Gamma, \textrm{out}}$]{
\begin{overpic}[trim= 2.5cm 1cm 2cm .5cm, clip=true, height=4cm]{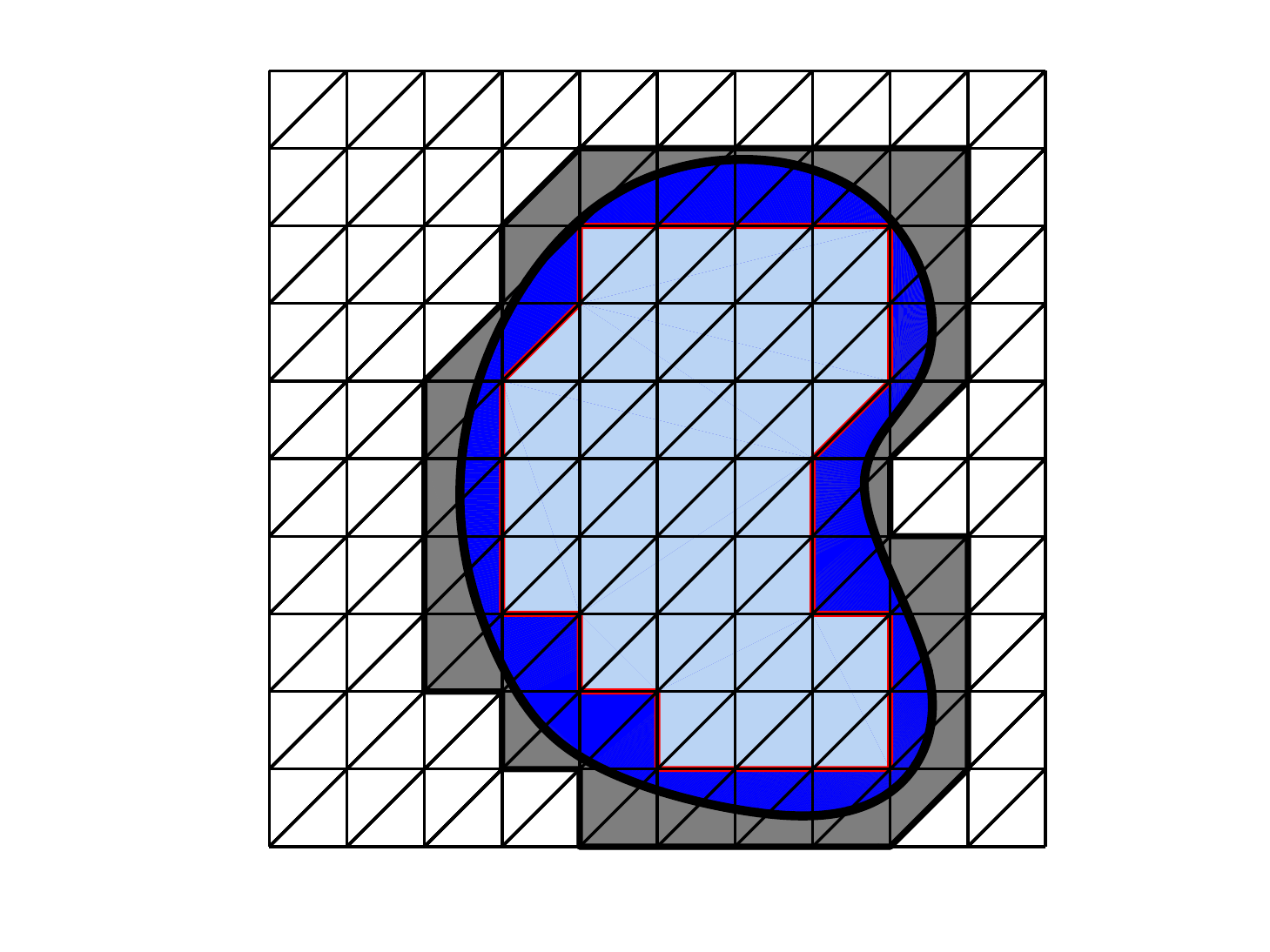}
\put(65,90.3){{\textcolor{white}{\small$\Omega_{\Gamma, \textrm{in}}$}}}
 \put(51.5,4){{\textcolor{white}{\small$\Omega_{\Gamma, \textrm{out}}$}}}
 \end{overpic}
}
\caption{A domain immersed into a three-directional mesh.}
\label{fig:weakBC_mesh}
\end{center}
\end{figure}

Let $\Omega_h$ be the domain of the mesh that covers the initial domain $\Omega$, so $\Omega \subseteq \Omega_h$. The boundary $\Gamma$ does not {necessarily} coincide with the boundary of the mesh.
More precisely, we set
\begin{itemize}
\item $\Omega_h=\Omega_\text{in}\cup\Omega_{\Gamma}$, where
$\Omega_\text{in}\subset\Omega$ is formed by a set of mesh cells 
in the interior of $\Omega$ and $\Omega_{\Gamma}$ is formed by a disjoint set of mesh cells that contains at least all cells cut by $\Gamma$;
 \item $\Omega_{\Gamma} = {\Omega_{\Gamma,\text{in}}}\cup\Omega_{\Gamma,\text{out}}$, where ${\Omega_{\Gamma,\text{in}}}:=\Omega\cap\Omega_{\Gamma}$ and $\Omega_{\Gamma,\text{out}}:=\Omega_h\setminus {\Omega}$.
\end{itemize}
This is illustrated in Figure~\ref{fig:weakBC_mesh} for a three-directional mesh. The choice of the domain $\Omega_{\Gamma}$, called the \emph{boundary strip}, will be discussed in Section~\ref{sec:strip}. 

In order to enforce the boundary conditions weakly, we reformulate problem (\ref{problem0}) by introducing an additional flux unknown $\bsigma$:
\begin{equation}
\label{problem1}
\left\{
\begin{array}{rll}
  -\kappa \Delta u + \ba\cdot\nabla u &=\ f, & \hbox{in }\ \Omega, \\
  \frac{1}{\kappa} \bsigma &=\ {\nabla u}, & {\hbox{in }\ \Omega_{\Gamma,\text{in}}}, \\
  u &=\ g, & \hbox{on } \Gamma. \\
\end{array}
\right.
\end{equation}
{Then, we consider two approximation spaces, $\VV_h \subset H^1(\Omega_h)$ and $\WW_h \subset (L^2(\Omega_{\Gamma}))^2$, for the discretization of the problem.
When the Peclet number is high ($\|\ba\|/\kappa \gg 1$), the advection is the dominating term and the diffusion is only important in very small layers. In such case we consider the standard {streamline upwind} Petrov--Galerkin (SUPG) approach (see, e.g., \cite{brooks1982}) in combination with the weak boundary approach.
More precisely, we look at the following symmetric discrete variational formulation of (\ref{problem1}):
 find $u_h \in \VV_h$ and $\bsigma_h \in \WW_h$, such that  
\begin{align}
\label{eq:WBC}
\left\{
 \begin{array}{ll}
  \kappa( \nabla u_h,\nabla v_h )_{\Omega} 
  + \langle \ba\cdot\nabla u_h,v_h \rangle_{\Omega}
  + { \frac{1}{{\eta}}} (\bsigma_h - \kappa\nabla u_h, \nabla v_h)_{\Omega_{\Gamma,\text{in}}} \\[1ex]
  + \langle \frac{\alpha}{2} u_h - \bsigma_h \cdot \bn,v_h \rangle_{\Gamma}
  + \delta_h\, \langle (-\kappa \Delta u_h + \ba \cdot \nabla u_h)(\ba \cdot \nabla u_h), v_h \rangle_\Omega\\[1ex]
  \hspace*{3.5cm}=\ \langle f,v_h \rangle_{\Omega} + \frac{1}{2}\langle \alpha g,v_h \rangle_{\Gamma}
  + \delta_h\, \langle f\, (\ba \cdot \nabla u_h), v_h\rangle_\Omega,
    & \forall v_h \in \VV_h,\\[2ex]
  - \frac{1}{{{\eta}}} (\frac{1}{ \kappa}\bsigma_h + \nabla u_h, \btau_h)_{\Omega_{\Gamma,\text{in}}}
  - \langle u_h, \btau_h \cdot \bn \rangle_{\Gamma} \\[1ex]
   \hspace*{3.5cm}=\ - \langle g, \btau_h\cdot \bn \rangle_{\Gamma},
    & \forall \btau_h \in \WW_h,
\end{array}
\right.
\end{align}
where $\bn$ denotes the outward unit normal to~$\Gamma$, 
$(\cdot,\cdot)_\omega$ denotes the $L^2$ inner product for vector functions over $\omega$,
and $\langle \cdot,\cdot\rangle_\omega$ denotes the $L^2$ inner product for scalar functions over $\omega$.
The parameter $\delta_h\geq0$ is the local SUPG stabilization parameter,
which is set to zero if stabilization is not needed.
The free parameters ${\eta}$ and ${\alpha}$ are fixed a priori.
In the following we consider $\eta=2/\kappa$ and }
\begin{equation*}
\alpha=\left\{\begin{array}{ll}
  -\ba\cdot\bn, & \ \mbox{if\ } \ba\cdot\bn<0, \\
  0, & \ \mbox{otherwise}.
\end{array}\right.
\end{equation*}

\begin{figure}[t!]
\begin{center}
\subfigure[domain]{
 \includegraphics[trim = 3.5cm 2.5cm 3.5cm 3cm, clip = true, height=6cm]{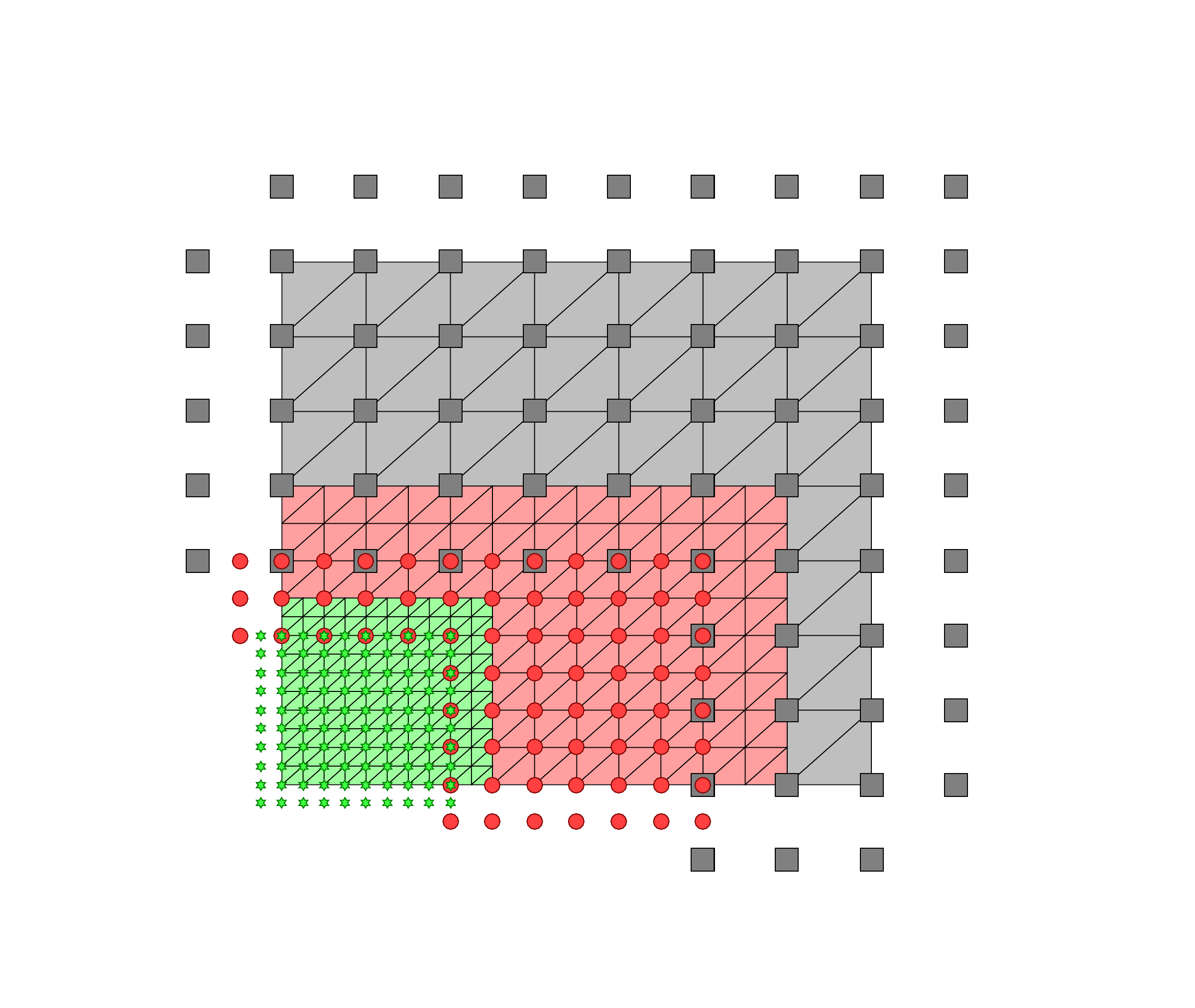}
}  \hspace*{0.5cm}
\subfigure[hierarchical fixed strip]{ 
 \includegraphics[trim = 3.5cm 2.5cm 3.5cm 3cm, clip = true, height=6cm]{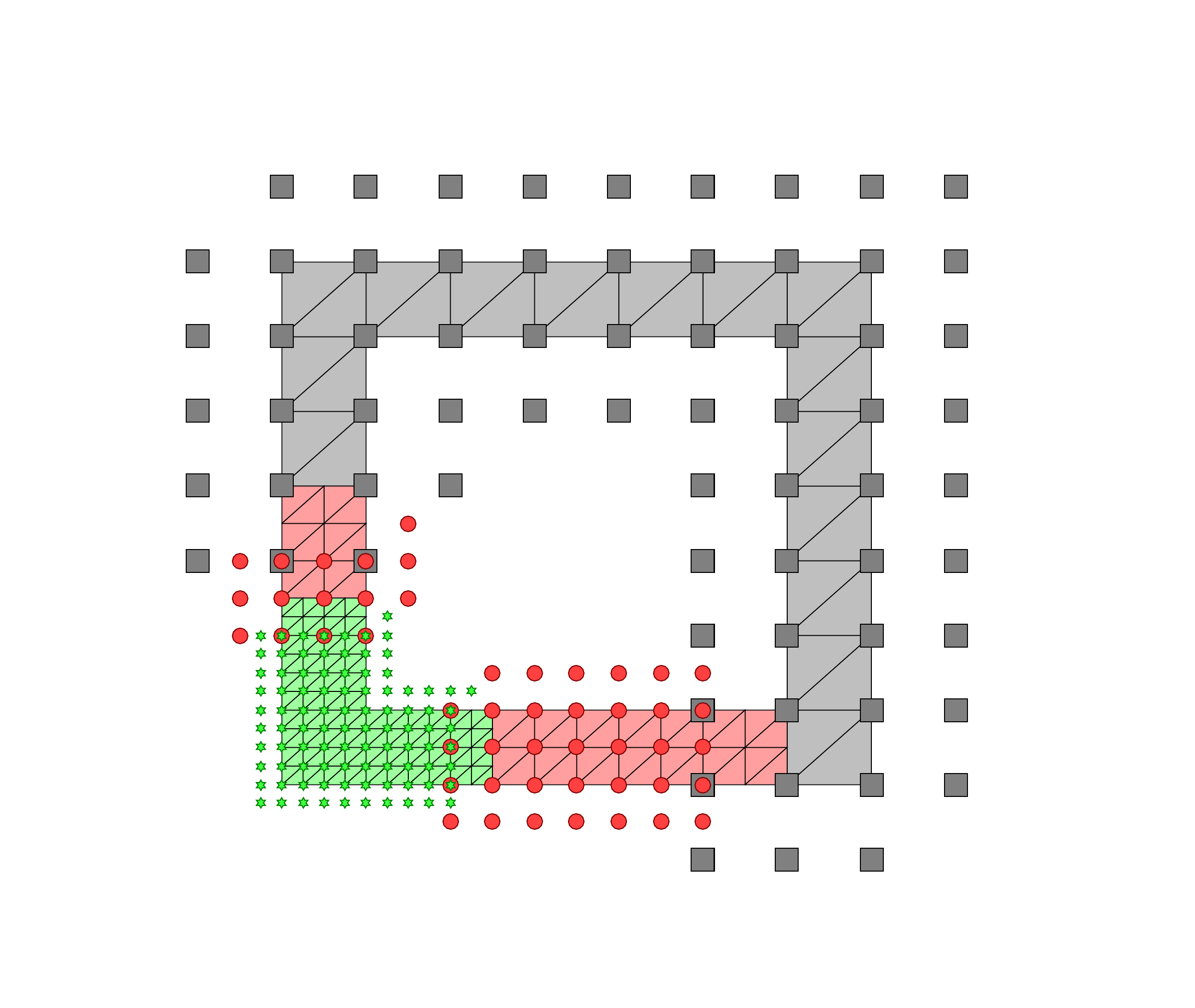}
} \\
\subfigure[HBox-spline strip]{ 
 \includegraphics[trim = 3.5cm 2.5cm 3.5cm 3cm, clip = true, height=6cm]{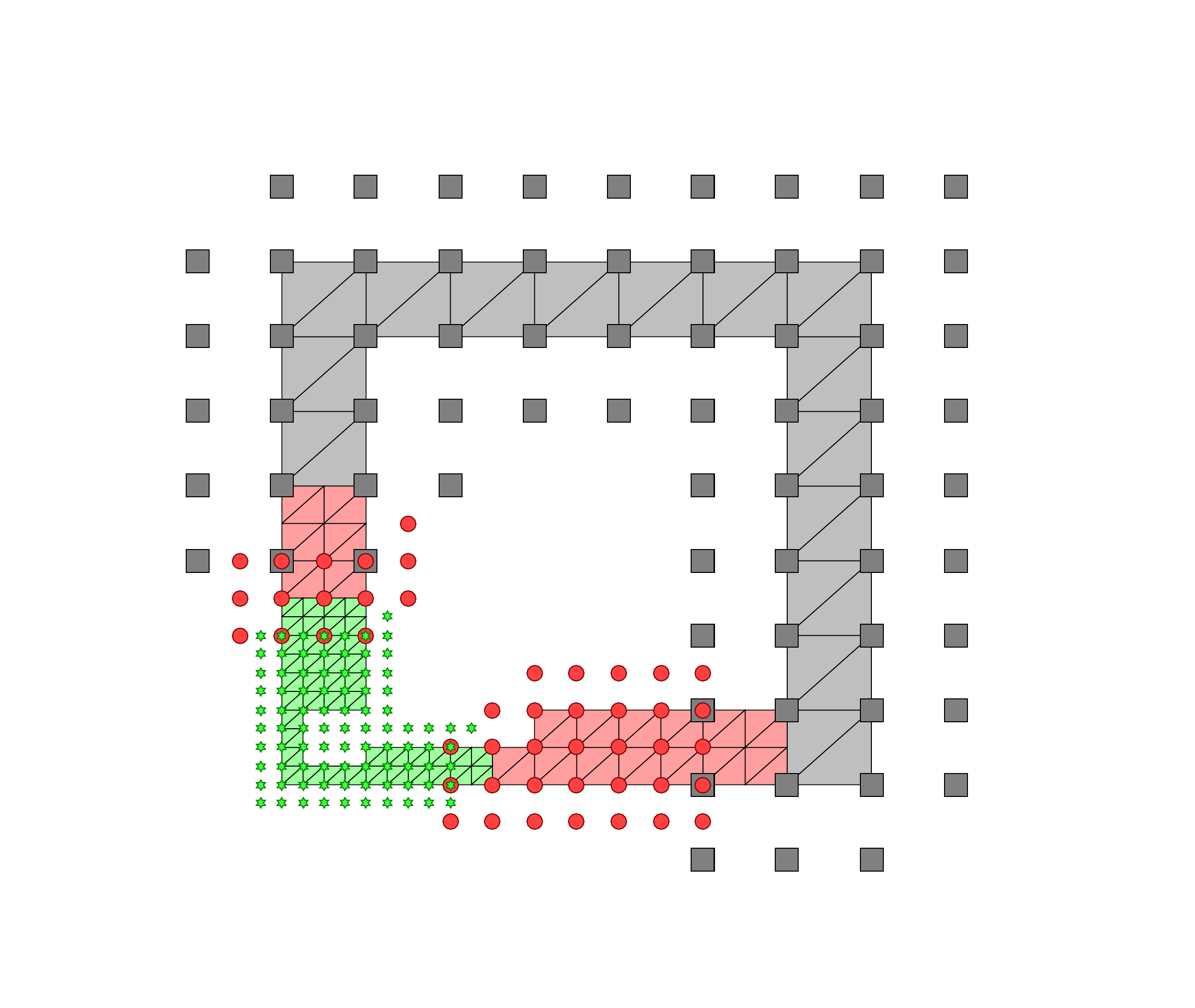}
}  \hspace*{0.5cm}
\subfigure[THBox-spline strip]{ 
 \includegraphics[trim = 3.5cm 2.5cm 3.5cm 3cm, clip = true, height=6cm]{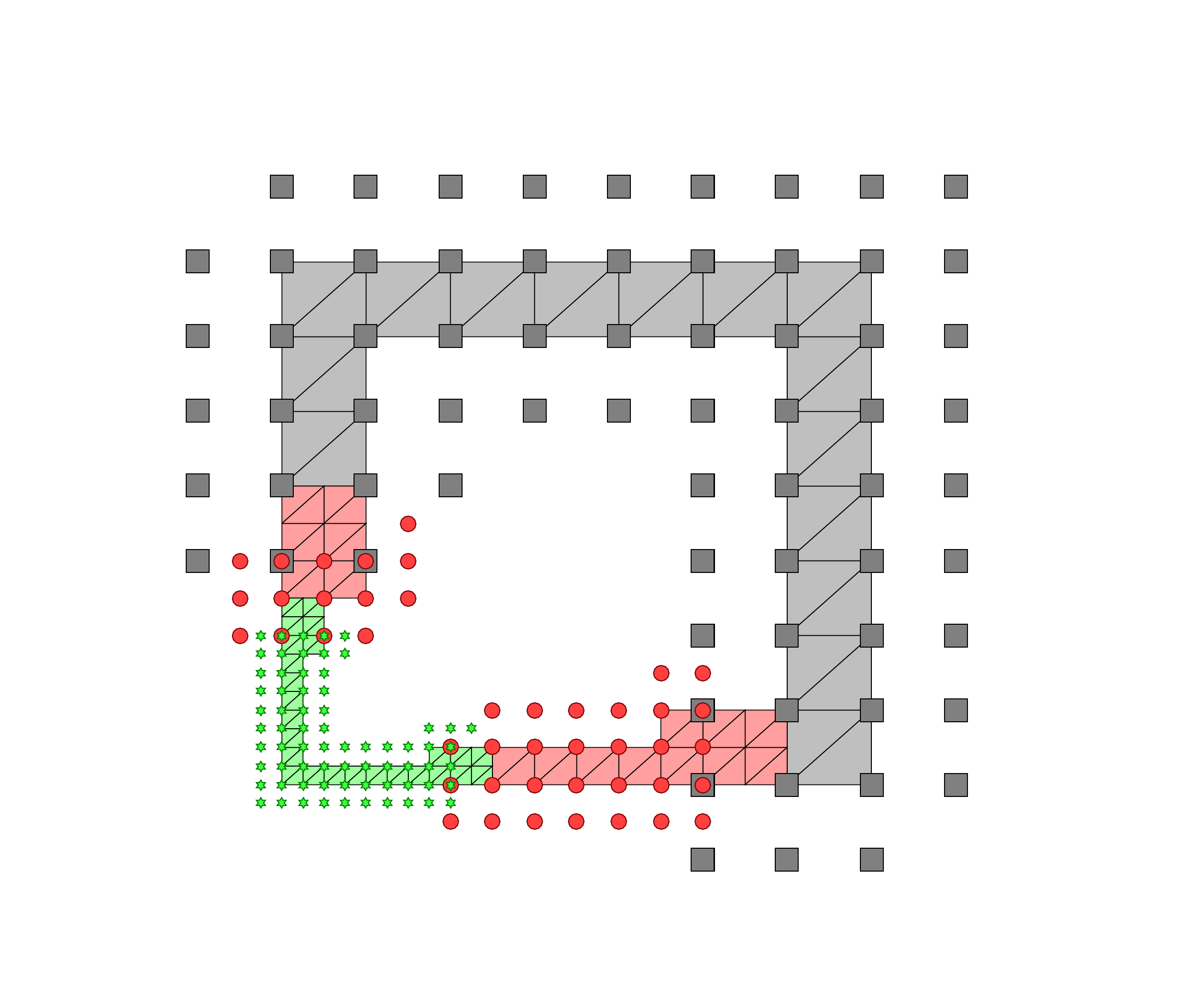}
}
\caption{A hierarchical three-directional mesh of three levels over the domain $\hat{\Omega}^0$, and three different types of the boundary strip $\hat{\Omega}_\Gamma$. The anchors of $C^2$ quartic (T)HBox-spline bases are indicated for different levels for $\hat \Omega$ (a) and $\hat \Omega_\Gamma$ (b)--(d), respectively.}
\label{fig:strips}
\end{center}
\end{figure}

We now specify the approximation spaces $\VV_h$ and $\WW_h$ in terms of hierarchical box splines.
Given a parametric domain $\hat{\Omega}^0$, suppose that the domain $\Omega$ can be described by a sufficiently smooth geometry function $\bF:\hat{\Omega}^0\rightarrow\Omega$, which is invertible and satisfies $\bF(\partial\hat{\Omega}^0)=\partial\Omega$. 
Then, using the (T)HBox-spline basis related to the hierarchy $\hat{\bOmega}:=\{\hat{\Omega}^{0},\ldots,\hat{\Omega}^{N-1}\}$, we specify
\begin{equation} \label{eq:Vh}
\cV_h:=\left\{{\beta}:=\hat{\beta}\circ\bF^{-1}:\hat{\beta}\in\hat{\cH}(\hat\bOmega)\right\},
\quad
\VV_h:=\left\langle \cV_h \right\rangle,
\end{equation}
where $\langle \cdot \rangle$ denotes the linear span of a set of functions.
Similarly, using the hierarchy of restricted subsets 
$\hat\bOmega_\Gamma:=\{\hat{\Omega}^{0}\cap\hat{\Omega}_\Gamma,\ldots,\hat{\Omega}^{N-1}\cap\hat{\Omega}_\Gamma\}$ where $\hat{\Omega}_\Gamma:=\bF^{-1}(\Omega_\Gamma)$,
we specify
\begin{equation} \label{eq:Wh}
\cW_h:=\left\{
{\beta}:=\hat{\beta}\circ\bF^{-1}:\hat{\beta}\in\hat{\cH}(\hat\bOmega_\Gamma)\right\},
\quad
\WW_h:=\left\langle \cW_h \right\rangle^2.
\end{equation}
The hierarchical construction in Definition~\ref{dfn:hb} guarantees that both sets $\cV_h$ and $\cW_h$ form a basis, and that $\left\langle \cV_h \right\rangle$ is a subspace of $\left\langle \cW_h \right\rangle$ if we restrict the domains to $\Omega_\Gamma$. In other words, near $\Omega_\Gamma$ the basis $\cW_h$ can locally consists of a larger number of elements than $\cV_h$. {This phenomenon is seen in Figure~\ref{fig:strips}, where some anchors corresponding to the basis $\cW_h$ are not active in the basis $\cV_h$.}
Since THBox-splines span the same hierarchical spline space as HBox-splines \cite{giannelli2014}, the bases $\hat{\cT}(\hat\bOmega)$ and $\hat{\cT}(\hat\bOmega_\Gamma)$ can also be considered in the definitions of $\cV_h$ and $\cW_h$, respectively.
An example of such hierarchical bases is shown in Figure~\ref{fig:strips} considering three-directional $C^2$ quartic box splines. The dots in the figure represent the anchors of the corresponding basis functions. 

Using the bases (\ref{eq:Vh}) and (\ref{eq:Wh}), 
the discrete weak formulation (\ref{eq:WBC}) leads to a discrete algebraic formulation of our problem. Let $n_h, m_h$ be the dimensions of $\VV_h$ and $\WW_h$, respectively, and
let $\bU\in \mathbb{R}^{n_h}, \bSigma\in \mathbb{R}^{m_h}$ be the unknown coefficient vectors describing $u_h$ and $\bsigma_h$, respectively. Then, we get the linear system
\begin{equation}
\label{eq:linear-sys}
  \left[
    \begin{array}{cc}
      \bK_{uu}+\bA_{uu}+\bK^{\text{in}}_{uu}+\bG_{uu} + \bS_1 + \bS_2 &  \bK_{u\sigma}+\bG_{u\sigma}\\
      \bK_{\sigma u}+\bG_{\sigma u} & \bK_{\sigma\sigma} \\
    \end{array}
  \right]\
  \left[
    \begin{array}{c}
      \bU\\
      \bSigma \\
    \end{array}
  \right]\ =\
  \left[
    \begin{array}{c}
      \ff+\bg_{u g} + \bs_f\\
      \bg_{\sigma g} \\
    \end{array}
  \right],
\end{equation}
    where the meaning of the different blocks is summarized in the diagram below:
\begin{equation*}
\hspace*{-0.17cm}
\begin{array}{c@{\hspace{.75em}}c@{\hspace{.75em}}c@{\hspace{.75em}}c@{\hspace{.75em}}c@{\hspace{.75em}}}
  \kappa( \nabla u_h,\nabla v_h )_{\Omega}
  &  \langle \ba\cdot\nabla u_h,v_h \rangle_{\Omega}
  & -\frac{\kappa}{\eta}( \nabla u_h,\nabla v_h )_{\Omega_{\Gamma,\text{in}}}
  & - \frac{1}{{{\eta}} \kappa}(\bsigma_h,\btau_h)_{\Omega_{\Gamma,\text{in}}}
  & \delta_h\, \langle -\kappa \Delta u_h\, (\ba \cdot \nabla u_h), v_h \rangle_\Omega
  \\
  \downarrow & \downarrow & \downarrow & \downarrow & \downarrow\\
  \bK_{uu}\bU & \bA_{uu}\bU & \bK^{\text{in}}_{uu}\bU& \bK_{\sigma\sigma}\bSigma &
  \bS_1 \bU\\
  &  & \\
  {\frac{1}{{\eta}}} ( \nabla u_h,\btau_h )_{\Omega_{\Gamma,\text{in}}} 
  & { \frac{1}{{\eta}}} (\bsigma_h,\nabla v_h)_{\Omega_{\Gamma,\text{in}}} 
  & \langle f,v_h \rangle_{\Omega}
  & \frac{1}{2}\langle \alpha g,v_h \rangle_{\Gamma}
  & \delta_h\, \langle (\ba \cdot \nabla u_h)^2, v_h \rangle_\Omega
  \\
  \downarrow & \downarrow & \downarrow & \downarrow & \downarrow
  \\
  \bK_{\sigma u}\bU & \bK_{u\sigma }\bSigma & \ff & \bg_{u g} & \bS_2 \bU
  \\
  &  &  \\
  \frac{1}{2}\langle \alpha u_h,v_h \rangle_{\Gamma} 
  & - \langle \bsigma_h \cdot \bn,v_h\rangle_{\Gamma} 
  & - \langle u_h,\btau_h \cdot \bn \rangle_{\Gamma} 
  & - \langle g,\btau_h \cdot \bn \rangle_{\Gamma} 
  &  \delta_h\, \langle f\, (\ba \cdot \nabla u_h), v_h\rangle_\Omega
    \\
 \downarrow & \downarrow & \downarrow & \downarrow & \downarrow
 \\
 \bG_{uu}\bU & \bG_{u\sigma}\bSigma & \bG_{\sigma u}\bU & \bg_{\sigma g}
 & \bs_f
  \end{array}
\end{equation*}
When writing the flux $\bSigma = \bK_{\sigma\sigma}^{-1}(-( \bK_{\sigma u}+\bG_{\sigma u})\bU+\bg_{\sigma g})$, the problem can be reformulated in terms of the original unknown $\bU$ only.

The described method fits in the class of immersed boundary methods. One of the main difficulties of immersed boundary methods is an accurate numerical integration over possibly very small mesh cells in $\Omega_{\Gamma, \text{in}}$ cut by the boundary. This is a subtle source of ill-conditioning and loss of accuracy.
In order to avoid such problems, we have combined the method with a geometry map $\bF$ describing the physical domain, according to the isogeometric philosophy. This allows us to take $\Omega_h=\Omega$ and $\Omega_{\Gamma}=\Omega_{\Gamma,\text{in}}$.
In this way we combine the benefits offered by immersed boundary and isogeometric methods in the hierarchical box spline context.

\subsection{Choice of the boundary strip} \label{sec:strip}

As described before, we assume $\Omega_h=\Omega$ and $\Omega_{\Gamma}=\Omega_{\Gamma,\text{in}}$.
The boundary strip $\Omega_\Gamma$ plays the role of the domain where the equality 
$\frac{1}{\kappa} \bsigma = {\nabla u}$ is enforced in the formulation (\ref{problem1}),
and the hierarchical basis $\cW_h$ in (\ref{eq:Wh}) can be constructed once $\Omega_\Gamma$ has been identified.
In this section we detail several strategies to choose such a boundary strip in the hierarchical setting.

A valid boundary strip must satisfy the following requirements.
First, the strip must contain all the (not further refined) cells of the (mapped) mesh that touch the boundary $\Gamma$ in order to properly enforce the conditions on the function $u$ and the flux $\bsigma$ on $\Gamma$. Second, due to interactions between different levels of the hierarchical basis and interactions between $\VV_h$ and $\WW_h$, it might be necessary to include additional cells into the strip. 
In particular, if the support of a HBox-spline of level $\ell_1$ overlaps a finer region $\hat \Omega^{\ell_2}$ ($\ell_1<\ell_2$) along a certain part of the boundary, then the boundary strip should not just contain cells of level $\ell_2$ at such a place but cells of level $\ell_1$ as well. Nested sequences of meshes guarantee that such coarser cells can always be tessellated with finer cells.

As we are working with a geometry map, it suffices to specify $\hat{\Omega}_\Gamma$ on the parametric domain $\hat{\Omega}^0$.
Three different boundary strips are introduced: 
\begin{enumerate}
  \item \emph{hierarchical fixed boundary strip}, 
  \item \emph{HBox-spline boundary strip}, 
  \item \emph{THBox-spline boundary strip}.
\end{enumerate}
The thickness of the first strip does not change along the boundary, while the last two are adaptively constructed in the spirit of the (truncated) hierarchical framework.
Recall that $h_\ell$ is the scaling factor corresponding to the box splines at level $\ell$, see (\ref{eq:set-box-ell}), and is a measure of the size of cells at level $\ell$ in the hierarchical mesh. 

\begin{configuration}
The {hierarchical fixed boundary strip} is a strip along the boundary $\partial\hat{\Omega}^0$ with a fixed thickness of~$h_0$. 
\end{configuration}
Due to the nested nature of the spaces in (\ref{eq:spaces}), at places where the finer region $\hat{\Omega}^\ell$ touches the boundary, the strip consists of $h_0/h_\ell$ rows of cells of level $\ell$. 
For example, when uniform dyadic refinement is considered on a three-directional mesh, along the part of the boundary refined at level $\ell$, the fixed boundary strip consists of $2^{\ell}$ rows of triangular cells of level $\ell$, as shown in Figure~\ref{fig:strips}(b).

The thickness of the strip effects the degrees-of-freedom of the overall system and its sparsity. If the refinements of the domain are mostly applied around the boundary, solving the system (\ref{eq:linear-sys}) can lead to significant additional computational costs. 
To construct the most compact system, the strip should be as thin as possible, meaning that the optimal choice should have preferably a thickness of $h_\ell$ at places where the space $\hat{\VV}^\ell$ is active at the boundary of the parametric domain $\hat{\Omega}^0$. 
In other words, we would like to shrink the uniform strip at regions that correspond to finer-level basis functions and where there is no influence of basis functions related to coarser levels. 

\begin{figure}[t!]
\begin{center}
\subfigure[$\hat H_{\Gamma}^0$]{ 
\includegraphics[trim = 3.75cm 2.5cm 4.5cm 3cm, clip = true, height=4.25cm]{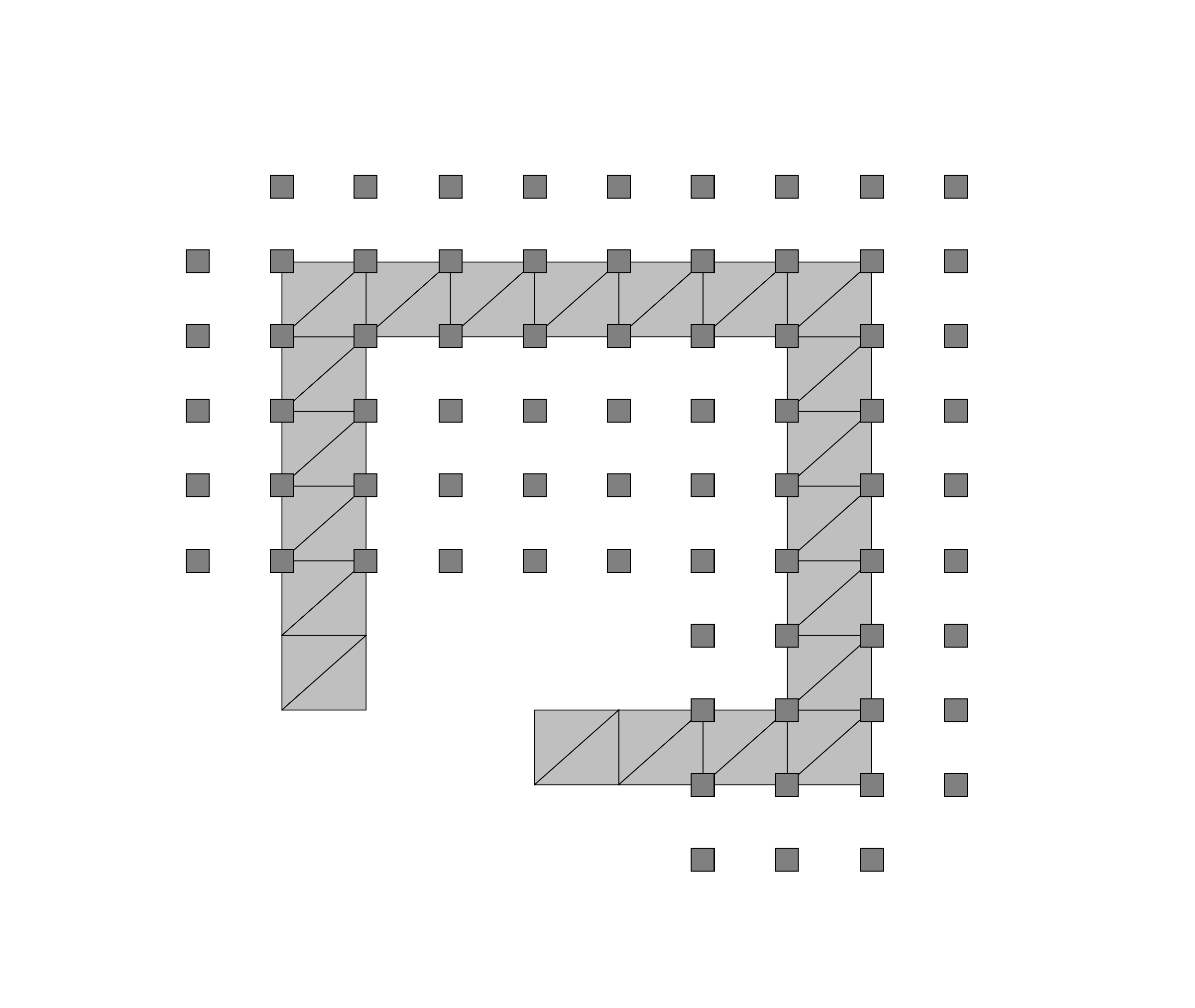}
}
\subfigure[$\hat H_{\Gamma}^1$]{
\includegraphics[trim = 4.75cm 2.5cm 8.25cm 3cm, clip = true, height=4.25cm]{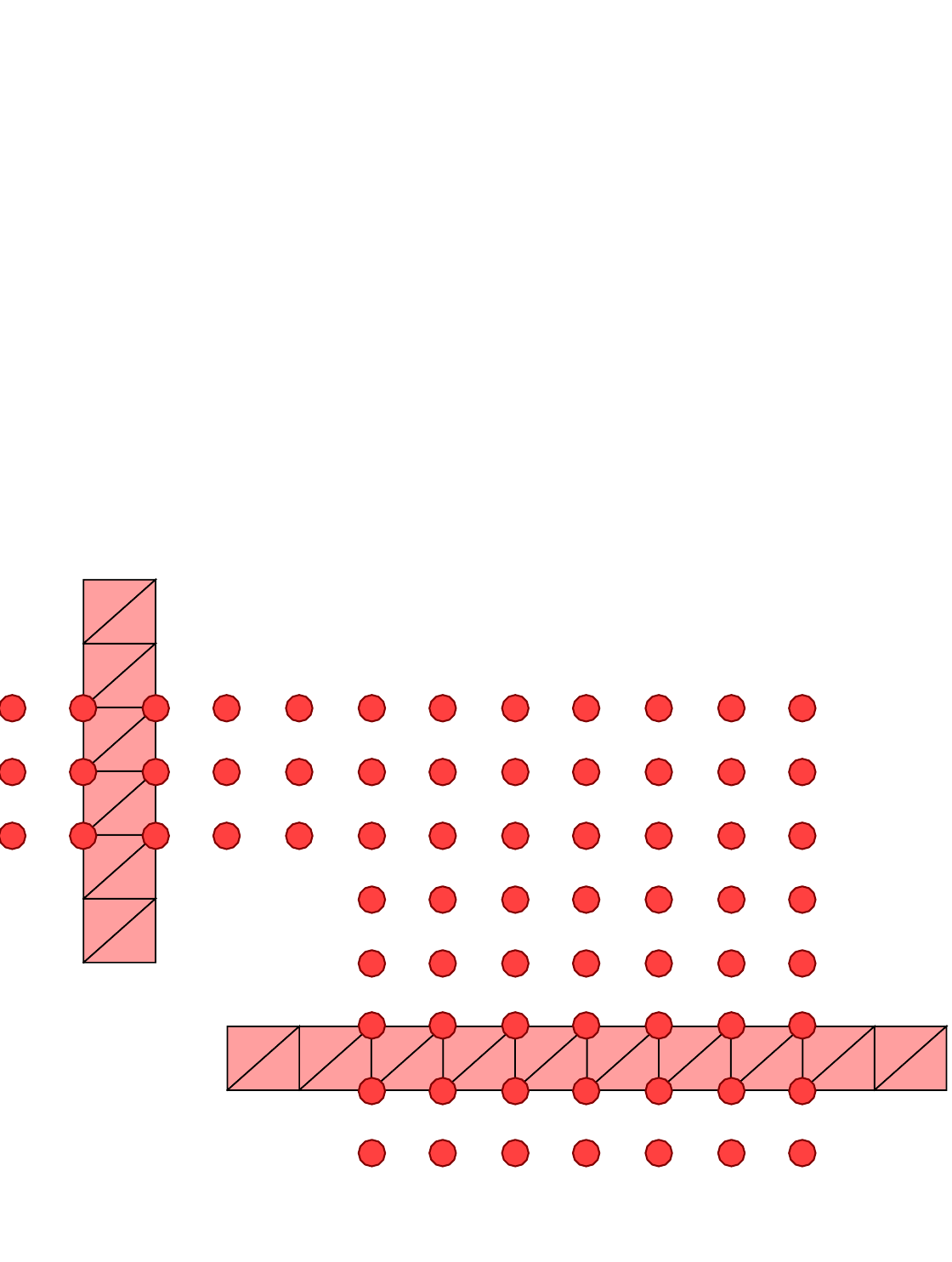}
}
\subfigure[$\hat H_{\Gamma}^2$]{
\includegraphics[trim = 5.25cm 2.5cm 14.25cm 3cm, clip = true, height=4.25cm]{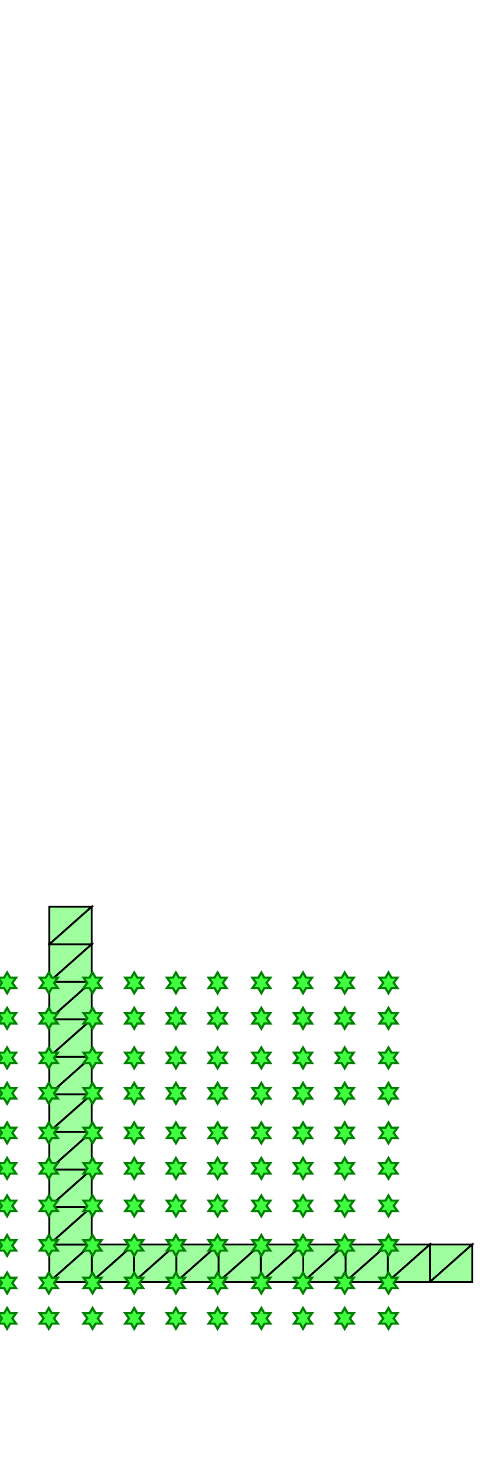}
}
\subfigure[$\hat \Omega_{\Gamma}$]{
\includegraphics[trim = 3.75cm 2.5cm 4.5cm 3cm, clip = true, height=4.25cm]{strip_adv1_mod-eps-converted-to}
}
\caption{Definition of the HBox-spline boundary strip (d) in terms of $\hat H_{\Gamma}^\ell$ with $\ell=0,1,2$ (a)--(c) for a three-directional configuration. The anchors of $C^2$ quartic HBox-splines in $\hat{\cH}(\hat\bOmega) \cap \hat {\mathcal B}^\ell$ are also shown in (a)--(c), and the anchors of the hierarchical basis $\hat{\cH}(\hat\bOmega_\Gamma)$ is shown in (d).}
\label{fig:HBoxstrip}
\end{center}
\end{figure}

\begin{configuration}
Let $\hat H_{\Gamma}^\ell$ be a strip of thickness $h_\ell$ along the part of the boundary $\partial\hat{\Omega}^0$ contained in the union of the supports of $\hat{\beta}^{\ell}\in\hat{\cH}(\hat\bOmega) \cap \hat {\mathcal B}^\ell$.
The HBox-spline boundary strip is defined as the union of the substrips $\hat H_{\Gamma}^\ell$, $\ell=0,\ldots,N-1$,
\begin{equation*}
\hat\Omega_{\Gamma} := \bigcup_{\ell=0}^{N-1} \hat H_{\Gamma}^\ell.
\end{equation*}
\end{configuration}
An example of the HBox-spline boundary strip together with its substrips is shown in Figure~\ref{fig:HBoxstrip}. 
By comparing Figure~\ref{fig:strips}(b) and Figure~\ref{fig:strips}(c), it is clear that this adaptive strip reduces the number of degrees-of-freedom with respect to the fixed boundary strip. 

\begin{figure}[t!]
\begin{center}
\subfigure[$\hat T_{\Gamma}^0$]{ 
\includegraphics[trim = 3.75cm 2.5cm 4.5cm 3cm, clip = true, height=4.25cm]{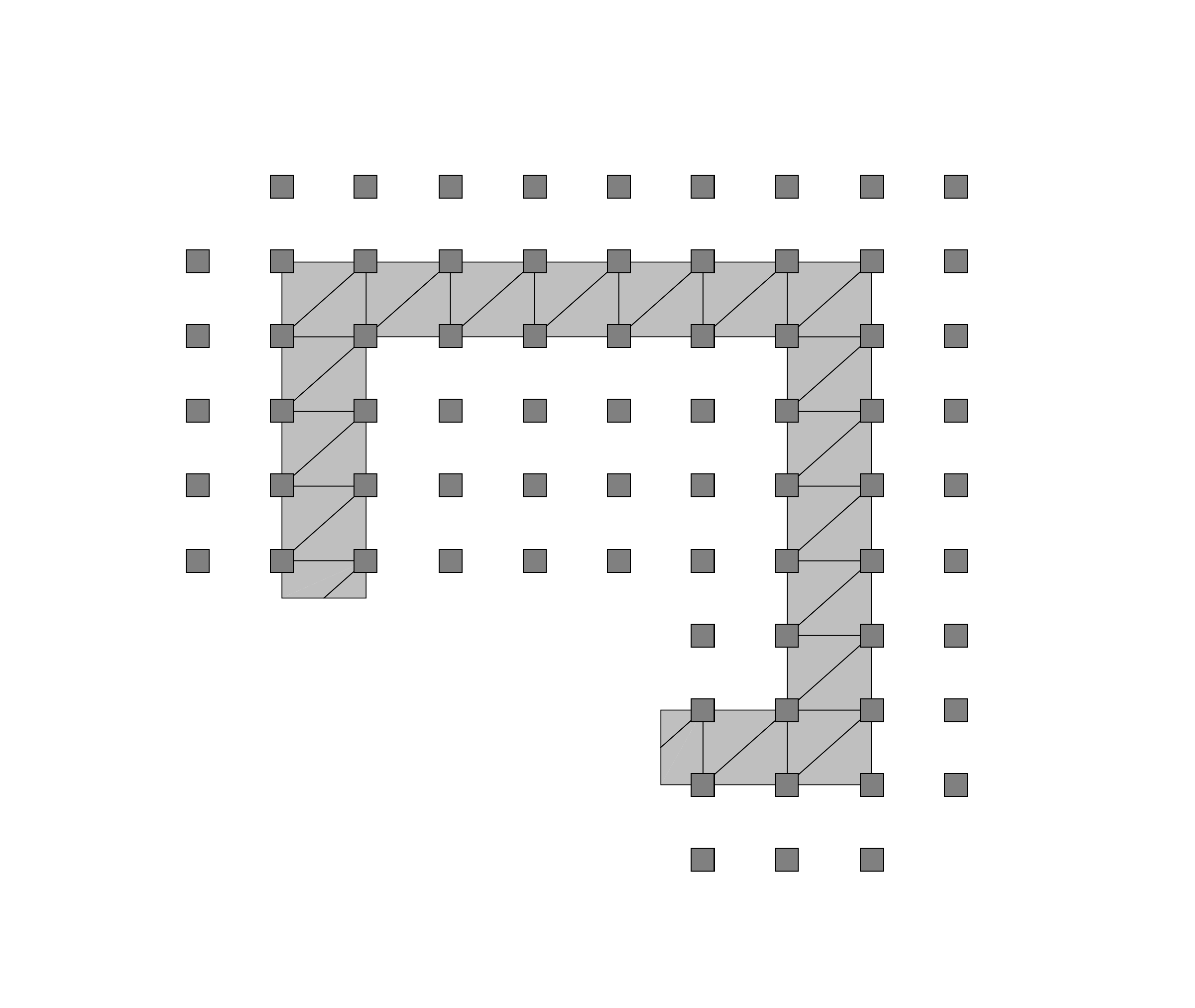}
}
\subfigure[$\hat T_{\Gamma}^1$]{
\includegraphics[trim = 4.75cm 2.5cm 8.25cm 3cm, clip = true, height=4.25cm]{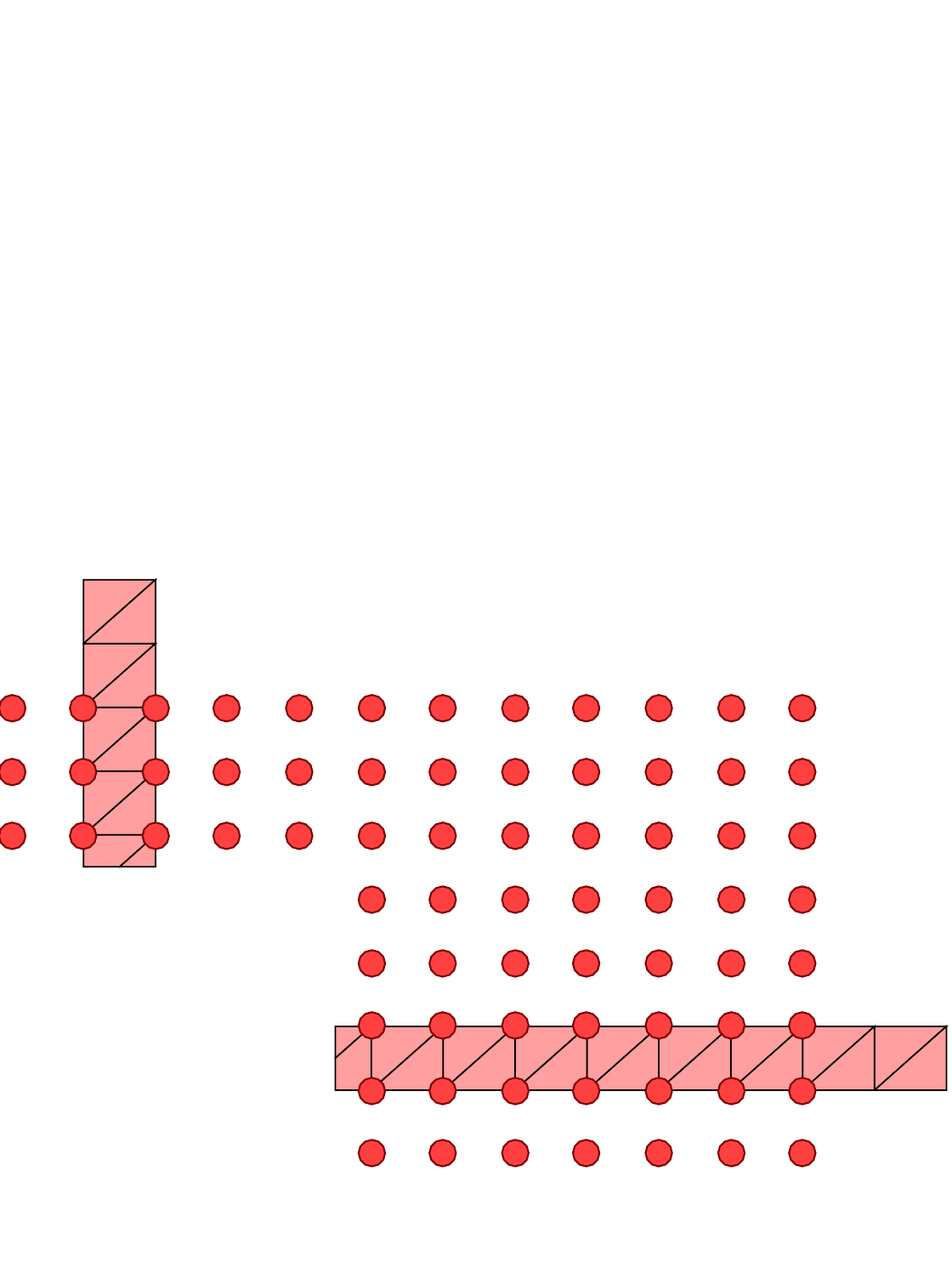}
}
\subfigure[$\hat T_{\Gamma}^2$]{
\includegraphics[trim = 5.25cm 2.5cm 14.25cm 3cm, clip = true, height=4.25cm]{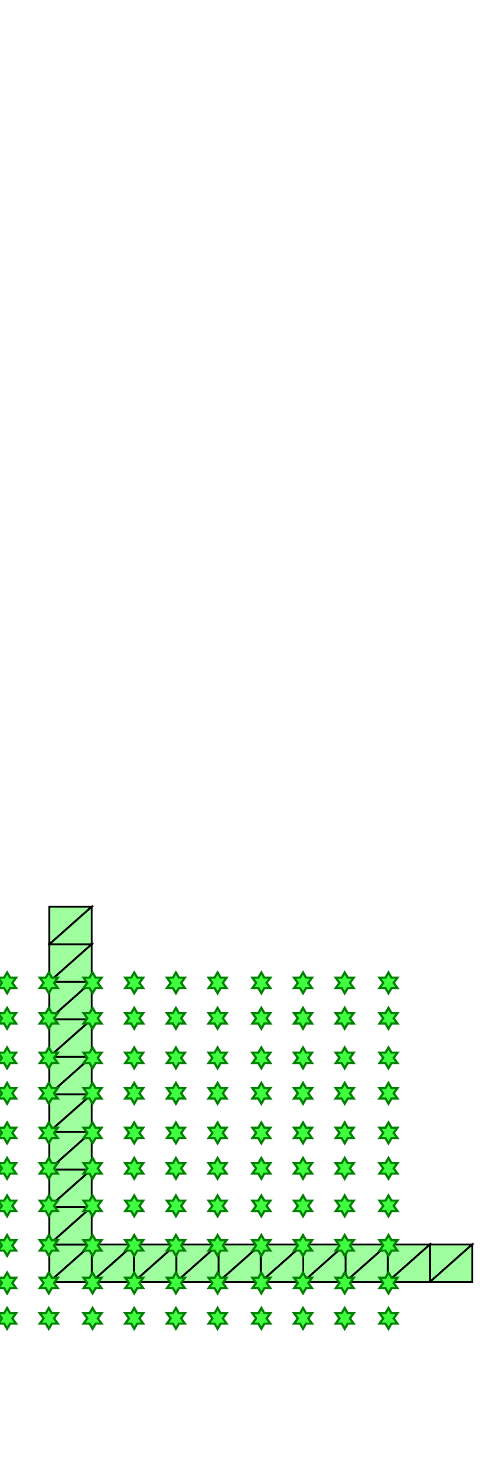}
}
\subfigure[$\hat \Omega_{\Gamma}$]{
\includegraphics[trim = 3.75cm 2.5cm 4.5cm 3cm, clip = true, height=4.25cm]{strip_adv2_mod-eps-converted-to}
}
\caption{Definition of the THBox-spline boundary strip (d) in terms of $\hat T_{\Gamma}^\ell$ with $\ell=0,1,2$ (a)--(c) for a three-directional configuration. The anchors of $C^2$ quartic THBox-splines in $\hat{\cT}(\hat\bOmega) \cap \hat {\mathcal B}^\ell$ are also shown in (a)--(c), and the anchors of the hierarchical basis $\hat{\cT}(\hat\bOmega_\Gamma)$ is shown in (d).}
\label{fig:THBoxstrip}
\end{center}
\end{figure}

For the THBox-spline boundary strip we follow the idea behind the truncation mechanism in Definition~\ref{dfn:thb}. In order to further shrink the boundary strip, we exploit the fact that the THBox-splines have a smaller support than the corresponding HBox-splines.

\begin{configuration}
Let $\hat T_{\Gamma}^\ell$ be a strip of thickness $h_\ell$ along the part of the boundary $\partial\hat{\Omega}^0$ contained in the union of the supports of ${\Trunc}^{\ell+1}(\hat{\beta}^{\ell})$ where $\hat{\beta}^{\ell}\in\hat{\cH}(\hat\bOmega) \cap \hat {\mathcal B}^\ell$.
The THBox-spline boundary strip is defined as the union of the substrips $\hat T_{\Gamma}^\ell$, $\ell=0,\ldots,N-1$,
\begin{equation*}
\hat\Omega_{\Gamma} := \bigcup_{\ell=0}^{N-1} \hat T_{\Gamma}^\ell.
\end{equation*}
\end{configuration}
An example of the THBox-spline boundary strip is shown in Figure~\ref{fig:THBoxstrip}. 
By comparing Figure~\ref{fig:strips}(c) and Figure~\ref{fig:strips}(d), we see that the THBox-spline strip allows us to further reduce the size of the boundary basis as naturally suggested by the hierarchical configuration.

\begin{figure}[t!]
\centering
\includegraphics[trim = 2cm .5cm 1cm 0cm, clip=true, height=2.5cm]{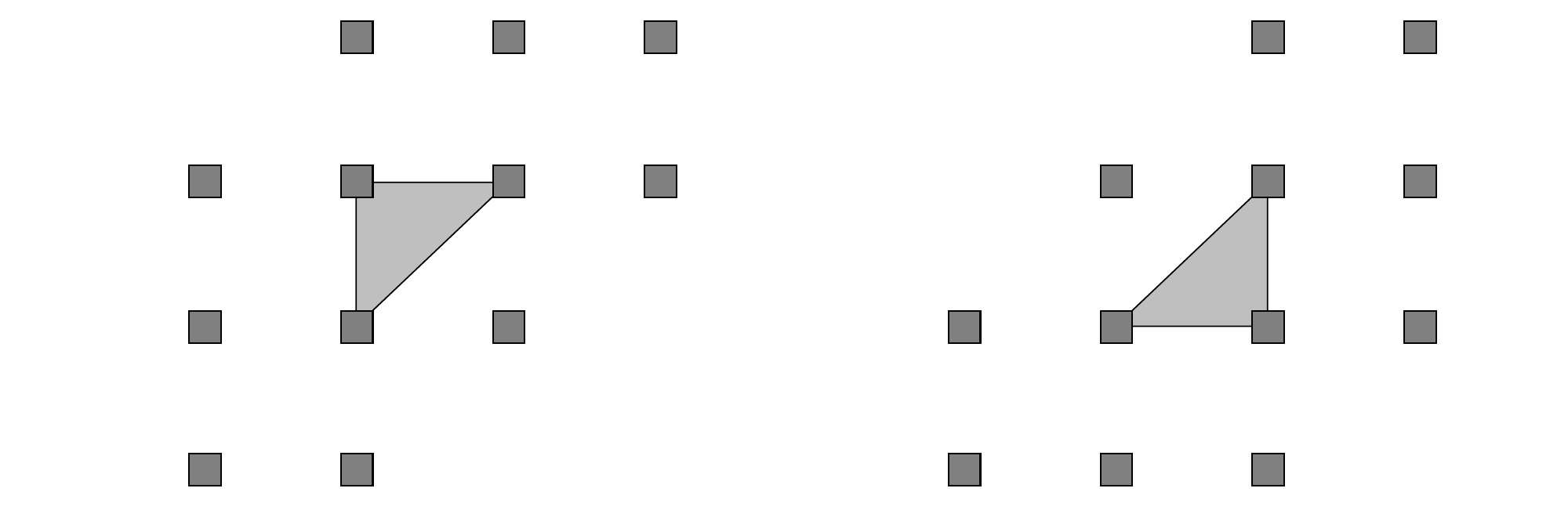}
\caption{All {three-directional $C^2$ quartic box splines} (represented by their anchors) that are non-zero on two types of triangles.}
\label{fig:neighbouringBasis}
\end{figure}

\begin{algorithmic}[ht!]
\begin{algorithm}
\STATE \textbf{Input}: hierarchical fixed strip $\hat \Omega_\Gamma$ 
\FOR {$\ell = 0,1, \dots, N-1$}
    \FOR {each cell $\pi^\ell$ in $\hat \Omega^\ell \cap \hat \Omega_\Gamma$}
        \STATE ${\cal N} := \{ \hat{\beta}^\ell \in \hat{{\cal B}}^\ell: \hat{\beta}^\ell |_{\pi^\ell} \not\equiv 0 \}$ 
        \STATE ${\cal N}_{\rm a} := \{ \hat{\beta}^\ell \in \hat{{\cal B}}^\ell: \hat{\beta}^\ell |_{\pi^\ell} \not \equiv 0 \; \wedge \;   \supp^{\hat \Omega_\Gamma} (\hat\beta^\ell) \subseteq \hat \Omega^\ell \cap \hat \Omega_\Gamma \}$ 
        \IF {${\rm closure}(\pi^\ell) \cap \partial\hat{\Omega}^0 = \emptyset \; \wedge \; {\cal N}={\cal N}_{\rm a}$}
            \STATE remove $\pi^\ell$ from strip $\hat \Omega_\Gamma$
        \ENDIF
    \ENDFOR
\ENDFOR
\caption{Cell removing algorithm for THBox-spline strip}
\label{alg:removeCells}
\end{algorithm}
\end{algorithmic}
The THBox-spline boundary strip can be alternatively constructed from the hierarchical fixed boundary strip by removing unnecessary cells (level by level) by following Algorithm~\ref{alg:removeCells}. 
A cell at level $\ell$ is removed if it does not touch the boundary and all box splines in $\hat{\cal B}^\ell$ 
that are non-zero on the cell are active on $\hat \Omega_\Gamma$. For example, for $C^2$ quartic THBox-splines a cell is removed if all 12 basis elements, whose supports overlap the cell, are active (see Figure~\ref{fig:neighbouringBasis}). 

\section{Numerical examples} \label{sec:exm}
In this section we present numerical experiments where we solve our
model problem (\ref{problem1}) 
using (truncated) hierarchical bivariate box splines over different domains. As previously mentioned, we
focus on $C^2$ quartic box splines defined over nested sequences of three-directional meshes
constructed by dyadic refinement, hence $h_{\ell+1}=h_{\ell}/2$.
The hierarchical approach is applied using the boundary strips described in Section~\ref{sec:strip}, and it is compared against the uniform case in terms of degrees-of-freedom and accuracy with respect to the exact solution (when available).
{In the first four examples predefined hierarchical refinements are considered (Section~\ref{sec:num-poisson}). In the last advection-diffusion problem, instead, the refinement is automatically performed using a simple gradient based a posteriori error estimator (Section~\ref{sec:num-advection}).}

\subsection{Poisson problem} \label{sec:num-poisson}

We start by considering four examples where (\ref{problem1}) reduces to the Poisson problem, i.e.,
\begin{equation} \label{problem1-poisson}
\left\{
\begin{array}{rll}
  -\Delta u &=\ f, & \hbox{in }\ \Omega, \\
  \bsigma &=\ {\nabla u}, & {\hbox{in }\ \Omega_{\Gamma,\text{in}}}, \\
  u &=\ g, & \hbox{on } \Gamma, 
\end{array}
\right.
\end{equation}
defined over different kinds of domains $\Omega\subset\RR^2$. 
{Of course, no SUPG stabilization is required here ($\delta_h=0$).}

\subsubsection{Hexagonal domain}

\begin{figure}[t!]
\begin{center}
\includegraphics[trim = 0.5cm 0.25cm 0.5cm 0cm, clip = true, height = 4cm]
{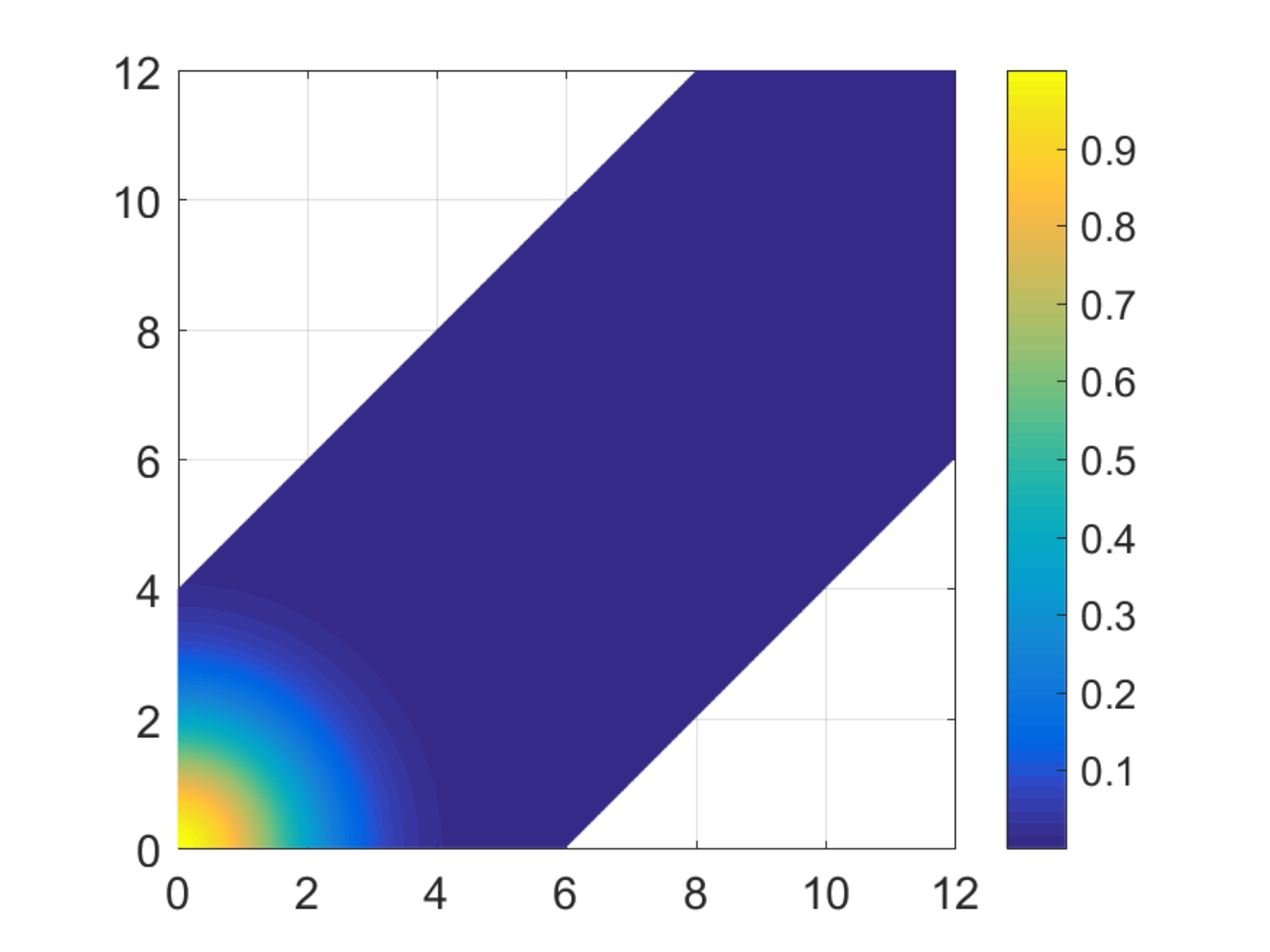}
\caption{The hexagonal domain problem: Exact solution.}
\label{fig:hexagon-sol}
\end{center}
\end{figure}

\begin{figure}
\begin{center}
\subfigure[$\Omega$, uniform]
{\includegraphics[trim = 0.5cm 0.5cm 0.5cm 0cm, clip = true, height = 4cm]
{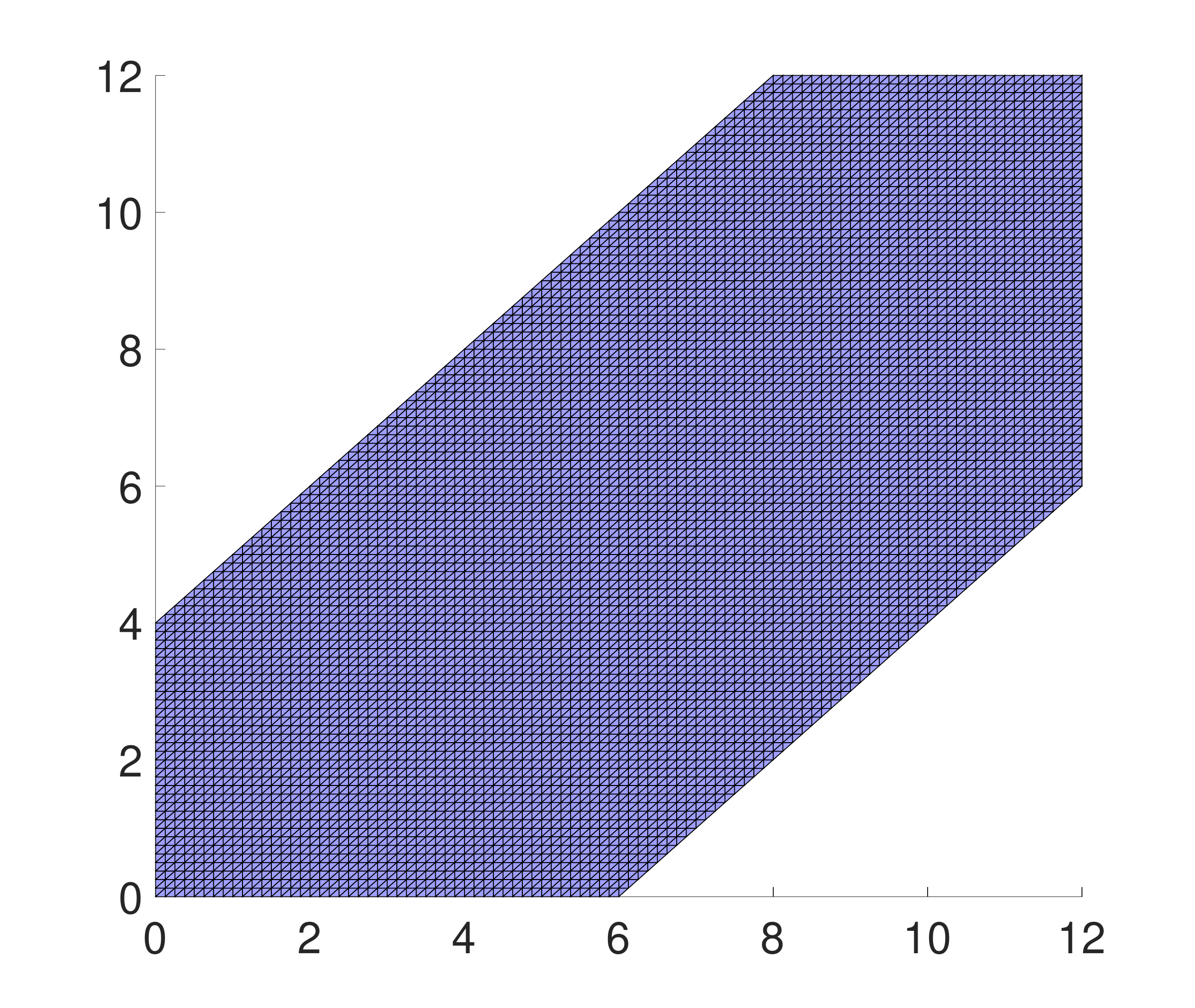}}
\subfigure[$\Omega_\Gamma$, uniform]
{\includegraphics[trim = 0.5cm 0.5cm 0.5cm 0cm, clip = true, height = 4cm]
{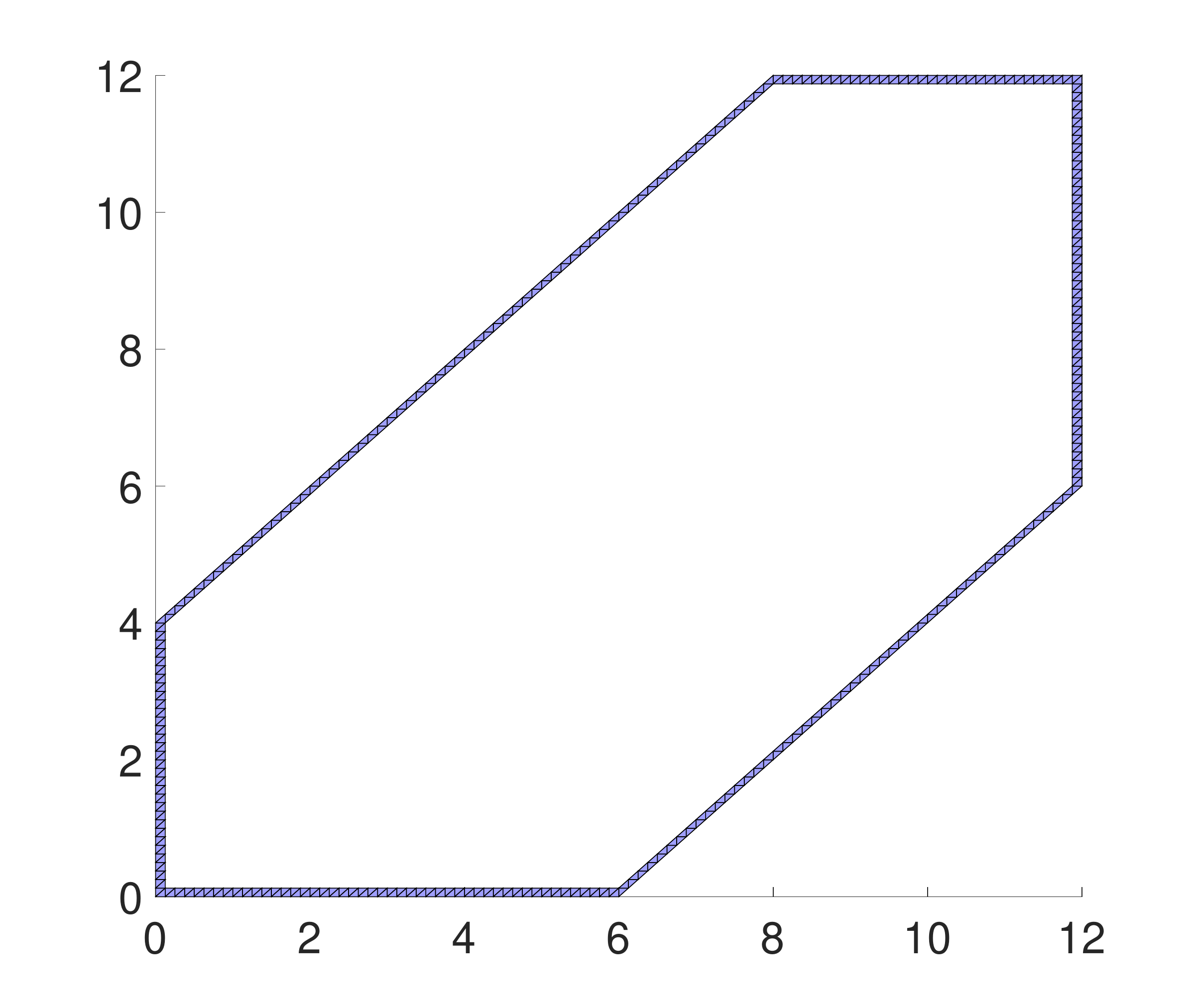}}
\subfigure[error, uniform]
{\includegraphics[trim = 0.5cm 0.25cm 0.5cm 0cm, clip = true, height = 4cm]
{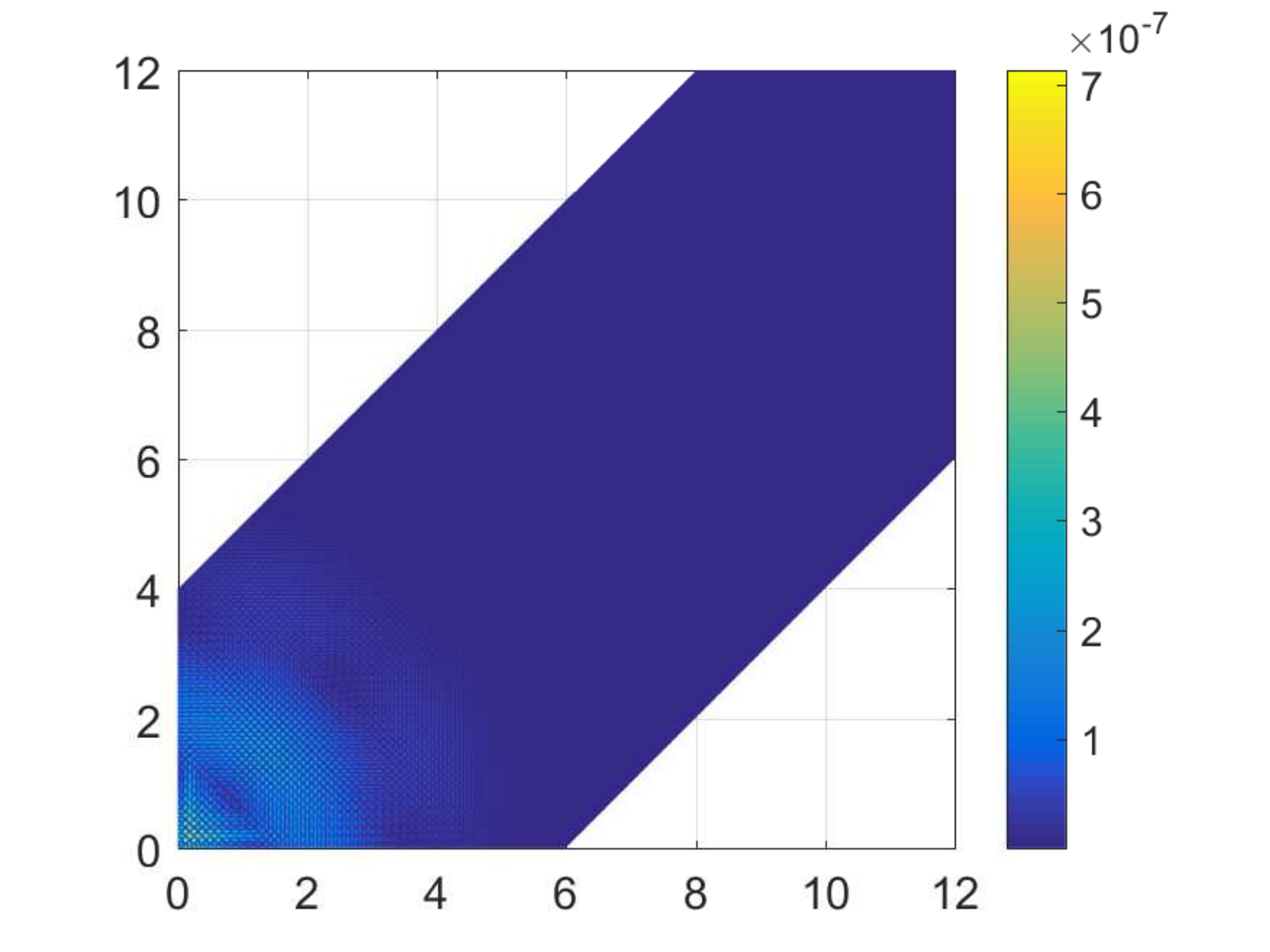}} \\
\subfigure[$\Omega$, hierarchical fixed]
{\includegraphics[trim = 0.5cm 0.5cm 0.5cm 0cm, clip = true, height = 4cm]
{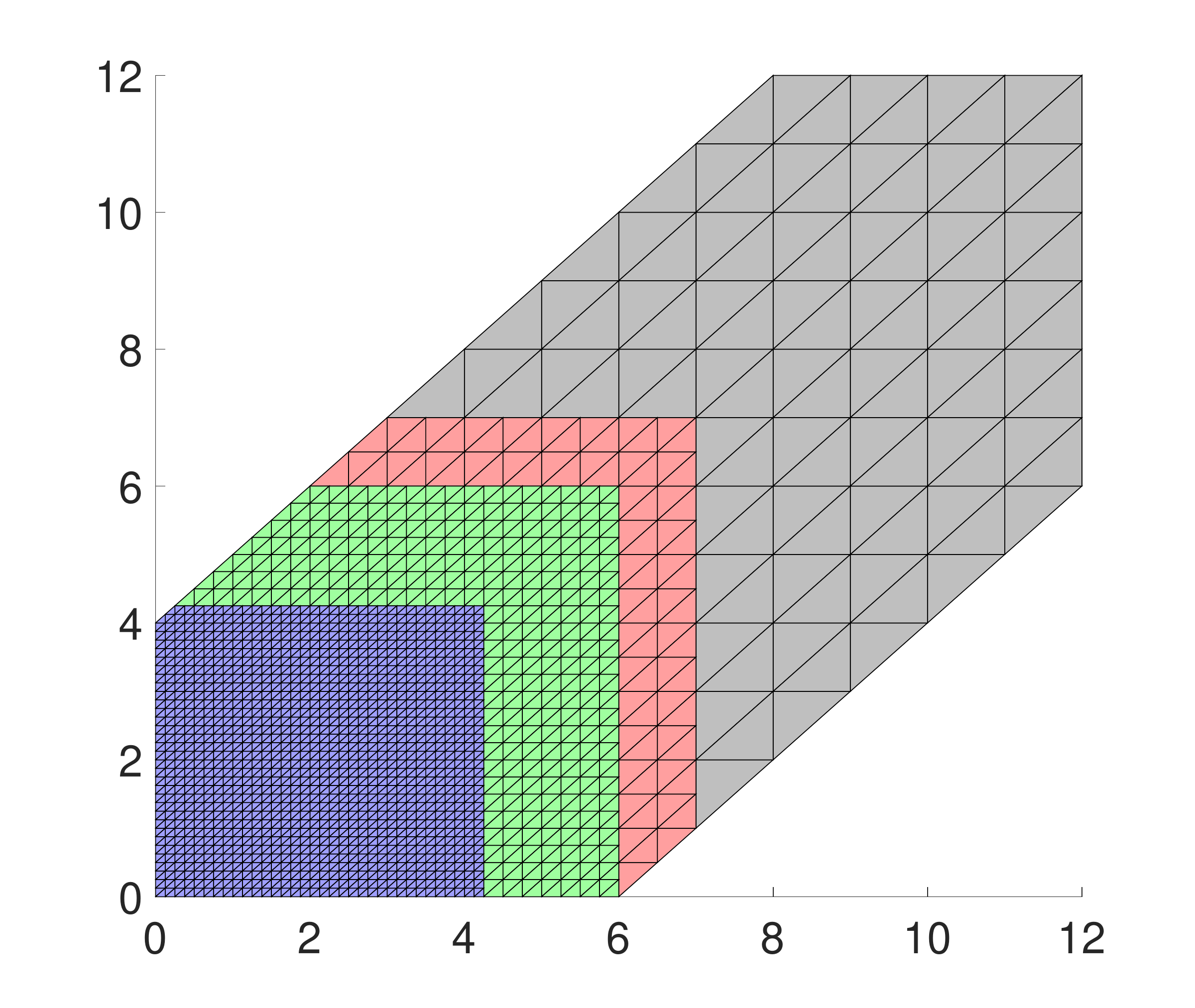}}
\subfigure[$\Omega_\Gamma$, hierarchical fixed]
{\includegraphics[trim = 0.5cm 0.5cm 0.5cm 0cm, clip = true, height = 4cm]
{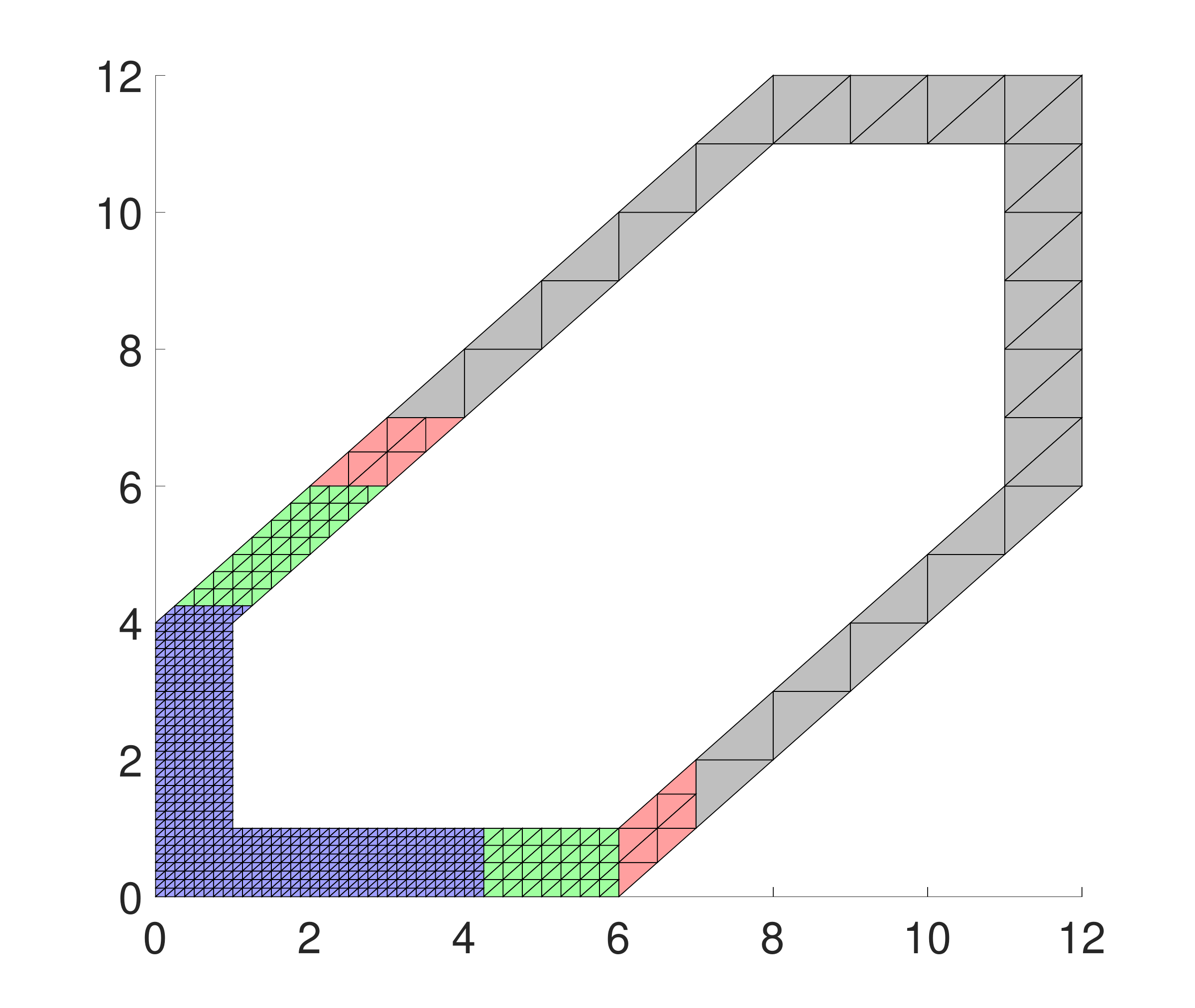}}
\subfigure[error, hierarchical fixed]
{\includegraphics[trim = 0.5cm 0.25cm 0.5cm 0cm, clip = true, height = 4cm]
{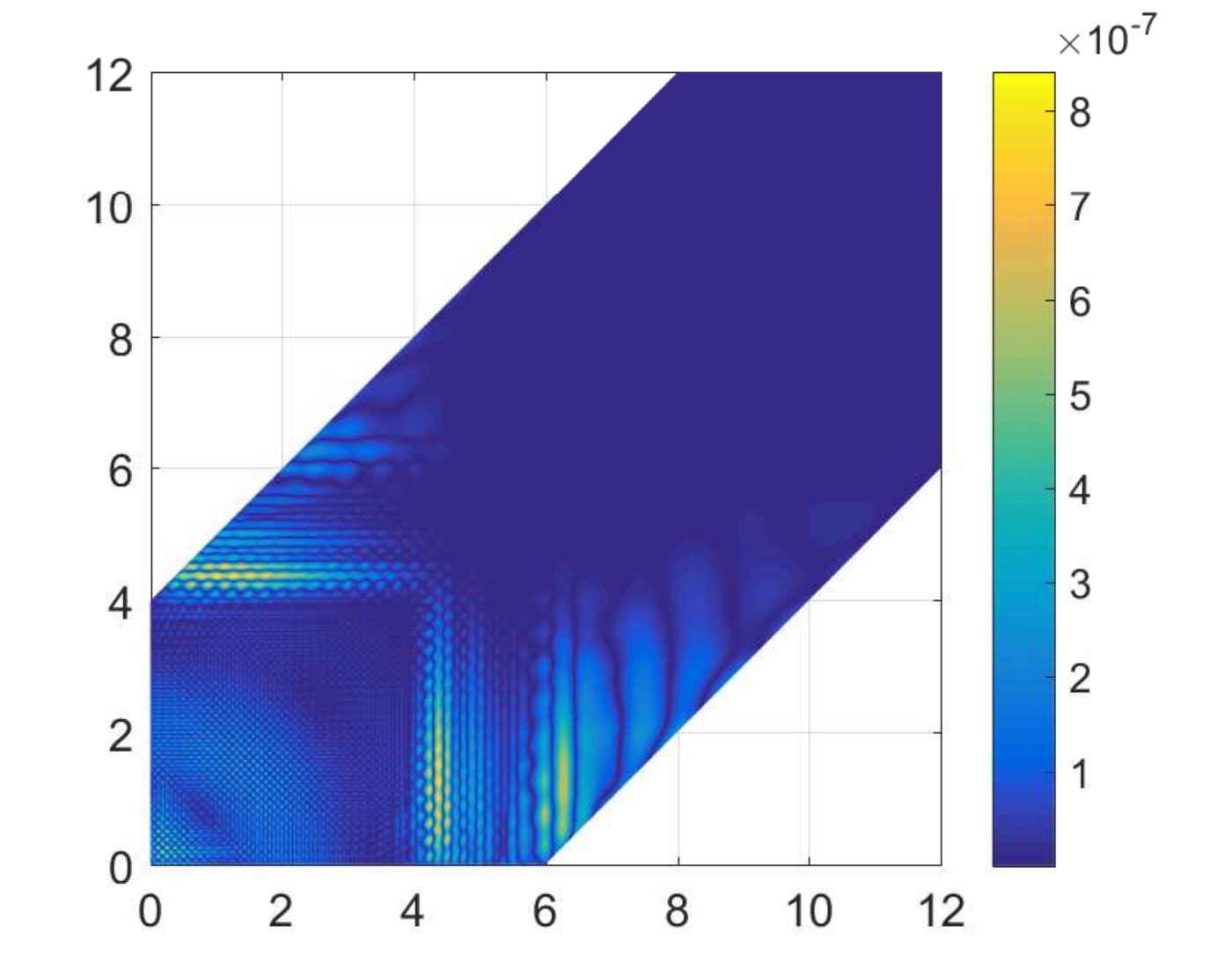}} \\
\subfigure[$\Omega$, HBox-spline]
{\includegraphics[trim = 0.5cm 0.5cm 0.5cm 0cm, clip = true, height = 4cm]
{hexagon_hierar4_omega-eps-converted-to}}
\subfigure[$\Omega_\Gamma$, HBox-spline]
{\includegraphics[trim = 0.5cm 0.5cm 0.5cm 0cm, clip = true, height = 4cm]
{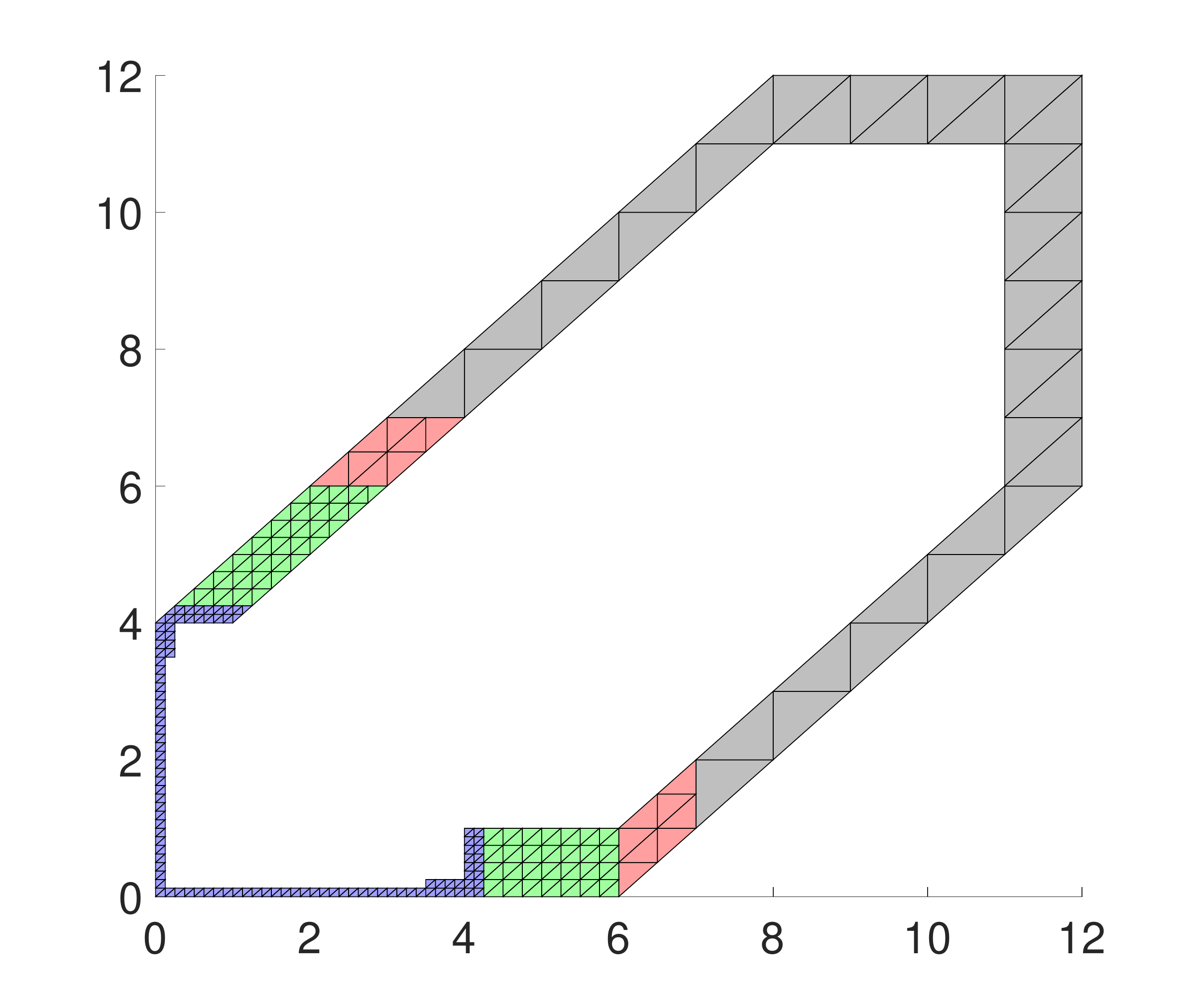}}
\subfigure[error, HBox-spline]
{\includegraphics[trim = 0.5cm 0.25cm 0.5cm 0cm, clip = true, height = 4cm]
{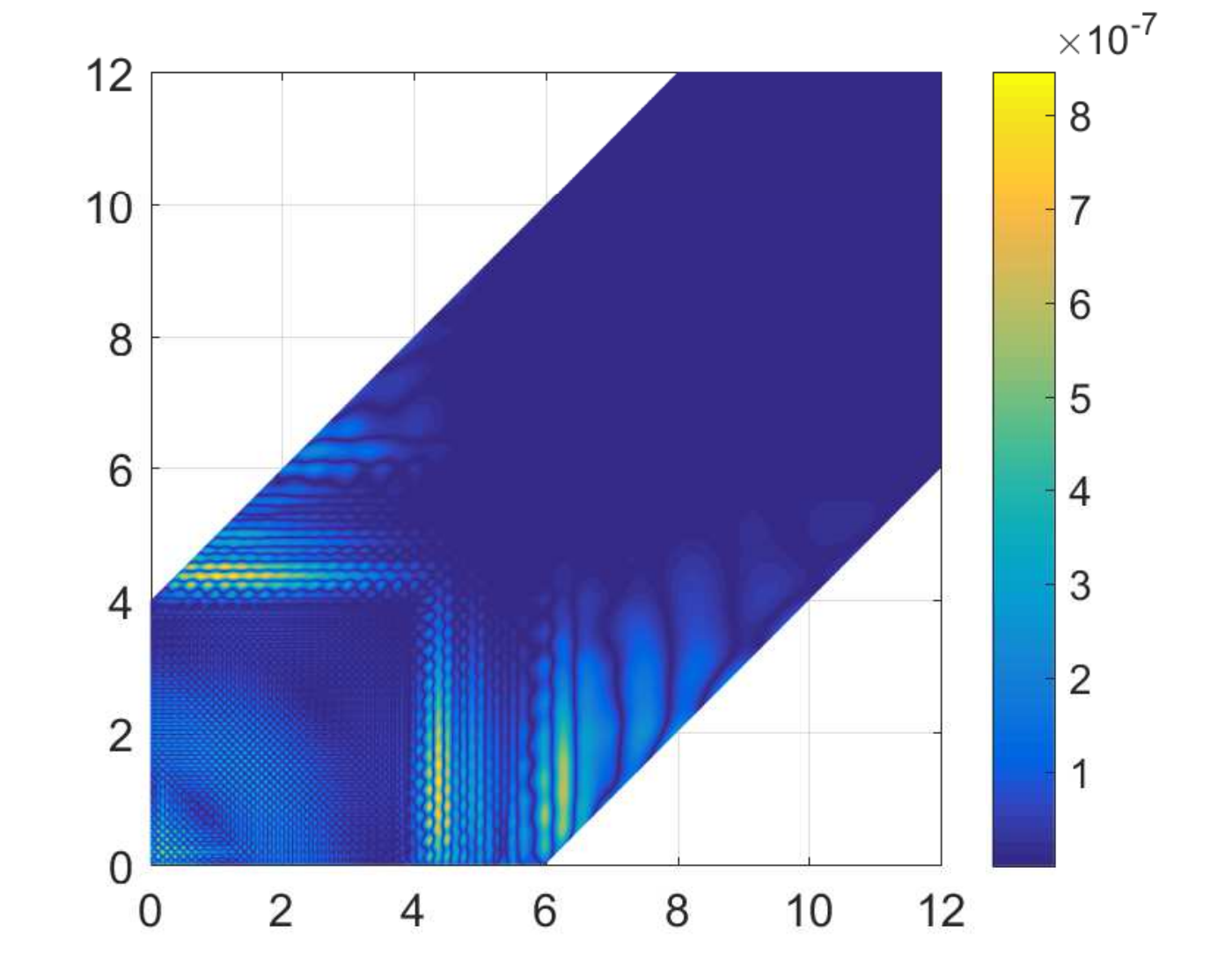}} \\
\subfigure[$\Omega$, THBox-spline]
{\includegraphics[trim = 0.5cm 0.5cm 0.5cm 0cm, clip = true, height = 4cm]
{hexagon_hierar4_omega-eps-converted-to}}
\subfigure[$\Omega_\Gamma$, THBox-spline]
{\includegraphics[trim = 0.5cm 0.5cm 0.5cm 0cm, clip = true, height = 4cm]
{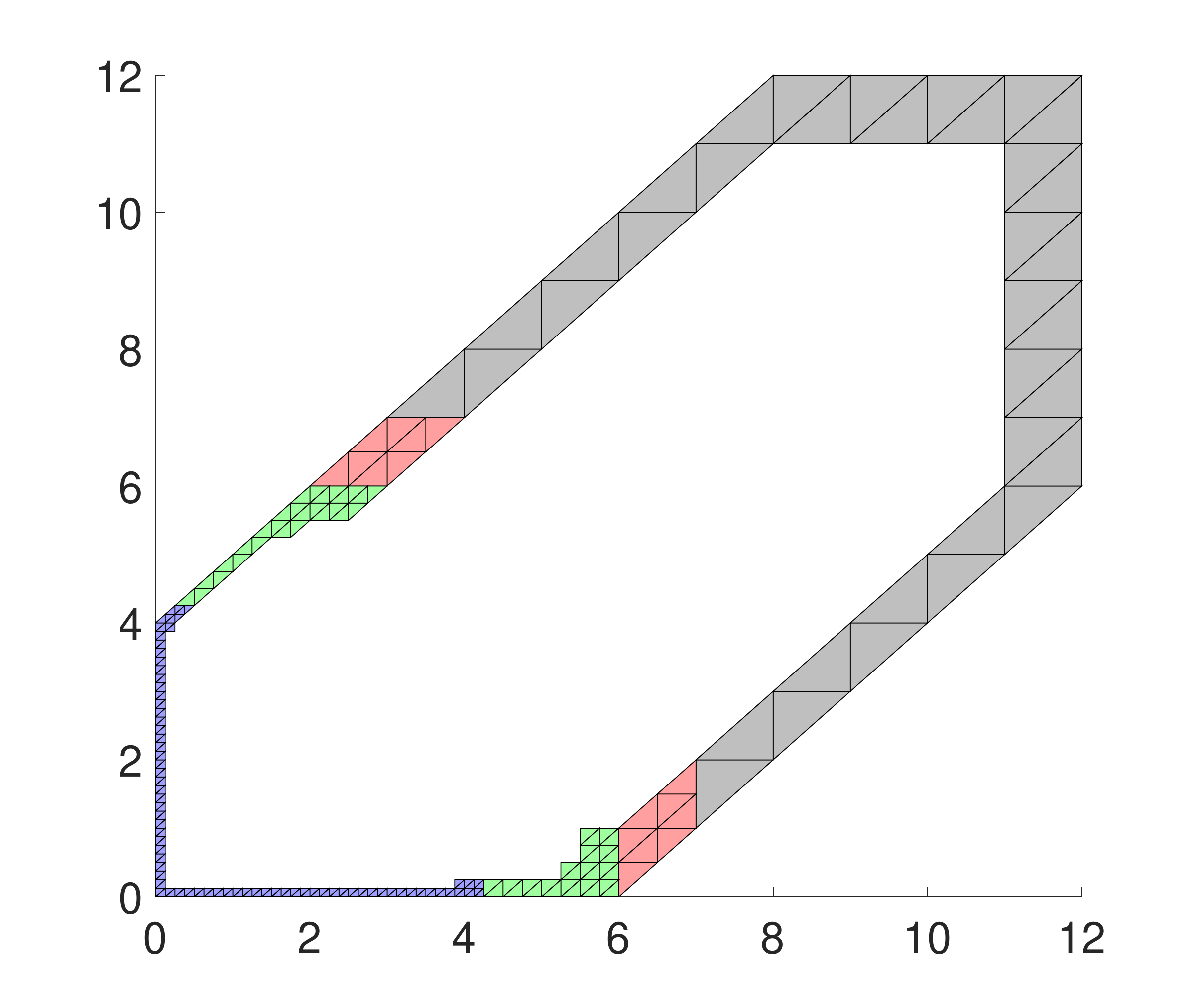}}
\subfigure[error, THBox-spline]
{\includegraphics[trim = 0.5cm 0.25cm 0.5cm 0cm, clip = true, height = 4cm]
{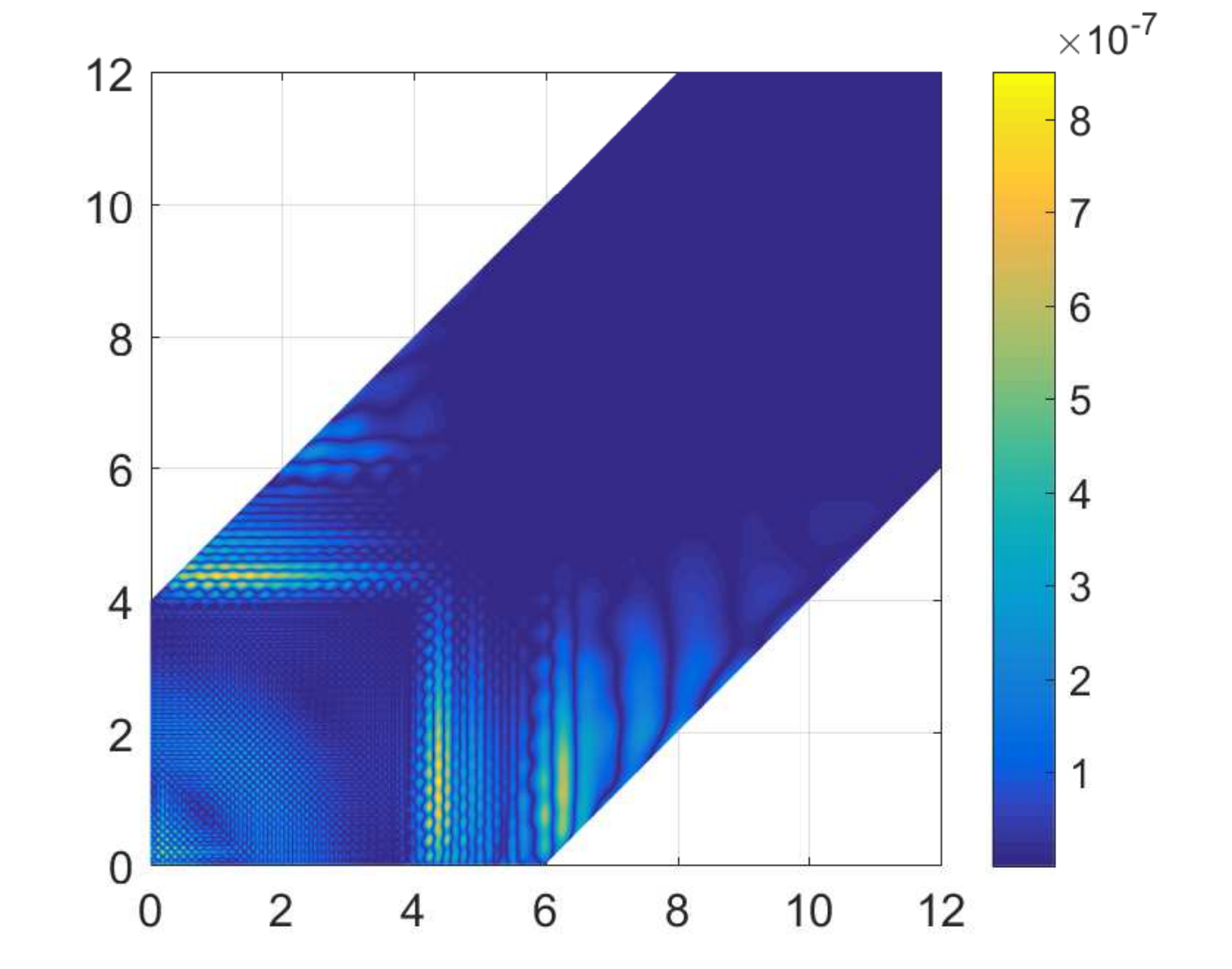}}
\caption{The hexagonal domain problem: Uniform and hierarchical meshes on $\hat \Omega = \Omega$ and $\hat \Omega_\Gamma = \Omega_\Gamma$, together with the obtained box spline errors ($N=4$).}
\label{fig:hexagon-domain}
\end{center}
\end{figure}

In this example we solve (\ref{problem1-poisson}) over the hexagonal domain $\Omega \subset [0,12]^2$ shown in Figure~\ref{fig:hexagon-sol}, where 
\begin{equation*}
g(x,y) = u(x,y) = \exp(- x^2/4 - y^2/4), \quad 
f(x,y) = -\Delta u(x,y) = - \exp(-x^2/4 - y^2/4) (x^2/4 + y^2/4 - 1).
\end{equation*}
Figure~\ref{fig:hexagon-domain} shows the uniform and hierarchical three-directional meshes used for this problem with $N=4$ (left column), the different hierarchical boundary strip constructions (center column), and the corresponding error functions $|u-u_h|$ (right column).
Table~\ref{tab:hexagon} collects the maximum values of the error with
respect to the exact solution, for different depths $N$ of uniform and
adaptive refinement. The finest mesh-size is $h_{N-1}=2^{-N+1}$.
We clearly observe that (T)HBox-splines maintain the optimal box spline convergence rate (order 4), but with a substantially smaller amount of degrees-of-freedom than in case of uniform refinement. 

\begin{table}[t!]
\begin{center}
\footnotesize
\begin{tabular}{c@{\hspace{.75em}}c{|}*{4}{c}
}
&& \multicolumn{4}{c}{error} 
 \\
\hline
$N$ & $h_{N-1}$ & uniform & fixed & HBox&  THBox 
 \\
\hline
1 & 1 & $4.44 \cdot 10^{-3}$ & $4.44 \cdot 10^{-3}$ &$4.44 \cdot 10^{-3}$ & $4.44 \cdot 10^{-3}$
\\
2 & 1/2 & $1.21 \cdot 10^{-3}$ & $1.13 \cdot 10^{-3}$ & $1.21 \cdot 10^{-3}$& $1.21 \cdot 10^{-3}$
\\
3 & 1/4 & $1.03 \cdot 10^{-5}$ &  $1.01 \cdot 10^{-5}$ &  $1.03 \cdot 10^{-5}$ & $1.03 \cdot 10^{-5}$
\\
4 & 1/8 & $7.13 \cdot 10^{-7}$ &  $8.39 \cdot 10^{-7}$&$8.45 \cdot 10^{-7}$ & $8.50 \cdot 10^{-7}$
\end{tabular}

\medskip

\begin{tabular}{c@{\hspace{.75em}}c{|}*{4}{c}{|}*{4}{c}}
&& \multicolumn{4}{c|}{dof} & \multicolumn{4}{c}{$\Omega_\Gamma$ dof}\\
\hline
$N$ & $h_{N-1}$ & uniform & fixed & HBox&  THBox & uniform  & fixed & HBox & THBox \\
\hline
1 & 1 & 152 & 152 & 152 & 152
  & 124 & 124& 124 & 124\\
2 & 1/2 & 485 & 290 & 290 &290
    & 260 & 199 & 190 &185\\
3 & 1/4 & 1715 & 704 & 704& 704
    & 532 & 391 & 326& 292\\
4 & 1/8 &  6431 & 1571 & 1571 & 1571
     & 1076 & 829 &490& 432
\end{tabular}
\caption{The hexagonal domain problem: Maximum error 
 and degrees-of-freedom (dof), for uniform and hierarchical box splines with different types of boundary strips and different levels of refinements $N=1,\ldots,4$.}
\label{tab:hexagon}
\end{center}
\end{table}


\subsubsection{Quarter of a circular domain}
\begin{figure}
\begin{center}
\includegraphics[trim = 0.5cm 0.25cm 0.5cm 0cm, clip = true, height = 4cm]
{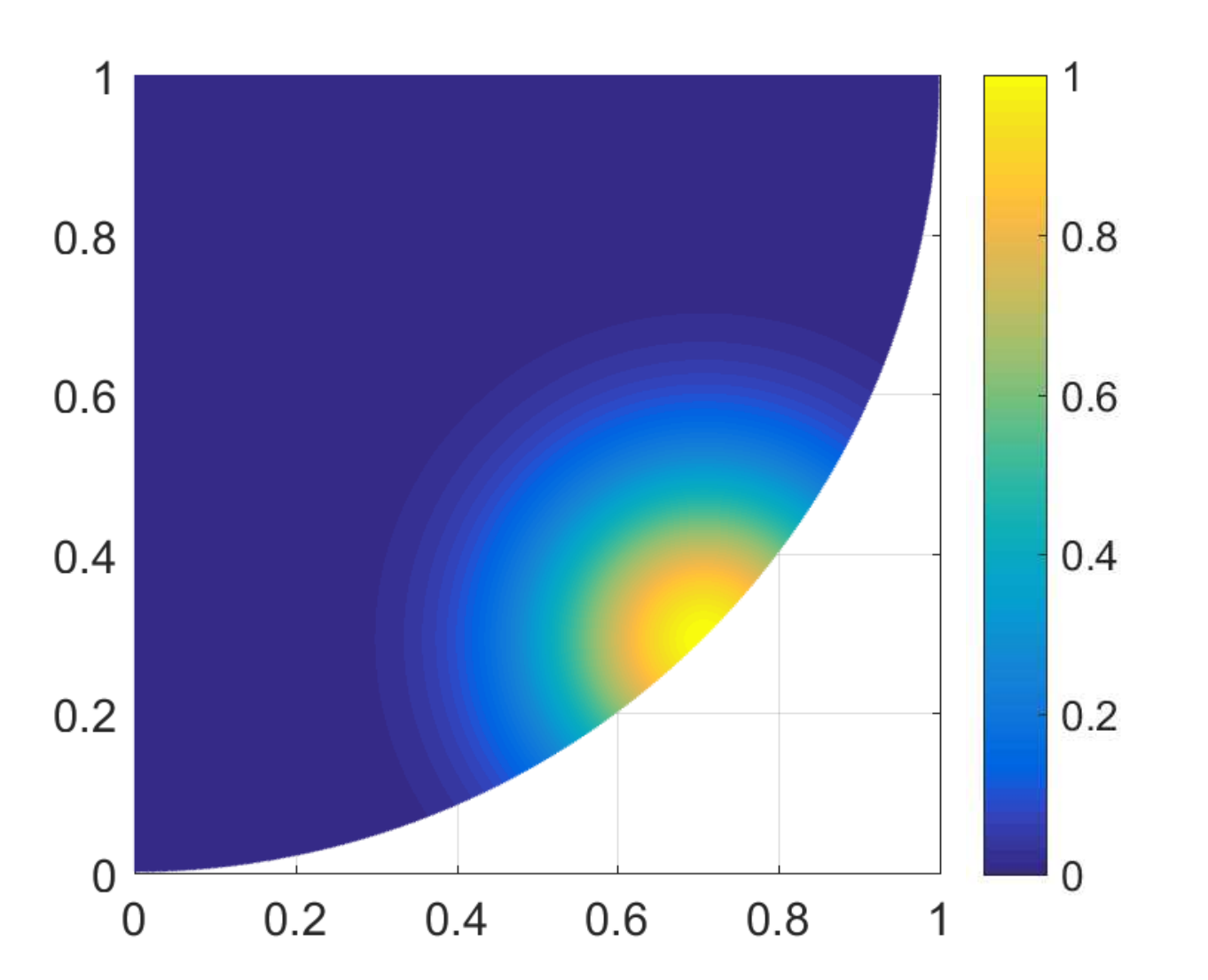}
\caption{The circular domain problem: Exact solution.}
\label{fig:circular-sol}
\end{center}
\medskip
\begin{center}
\subfigure[$\hat \Omega$, uniform]
{\includegraphics[trim = 0.5cm 0.5cm 0.5cm 0cm, clip = true, height = 4cm]
{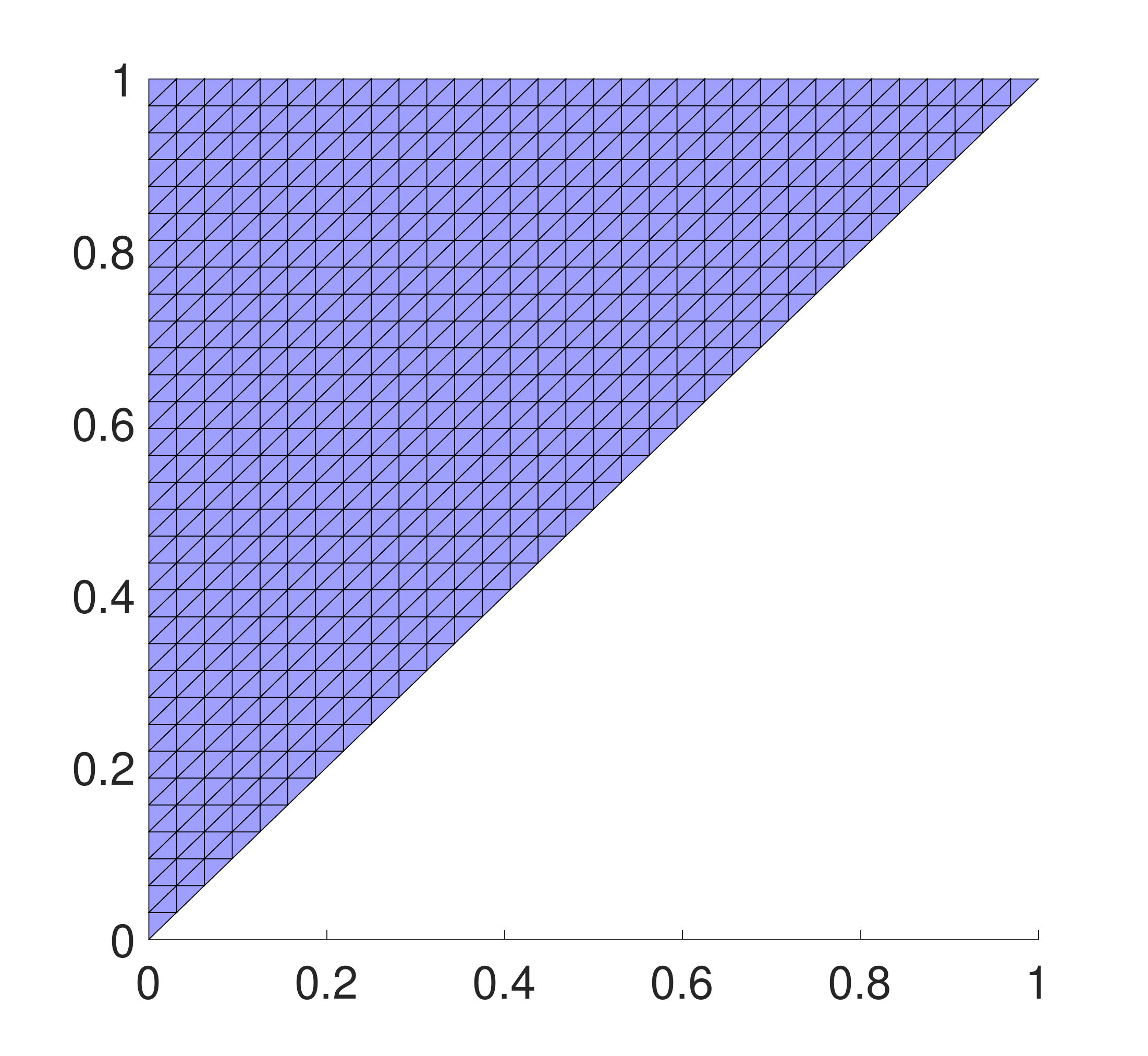}}
\subfigure[$\hat \Omega_\Gamma$, uniform]{
\includegraphics[trim = 0.5cm 0.5cm 0.5cm 0cm, clip = true, height = 4cm]
{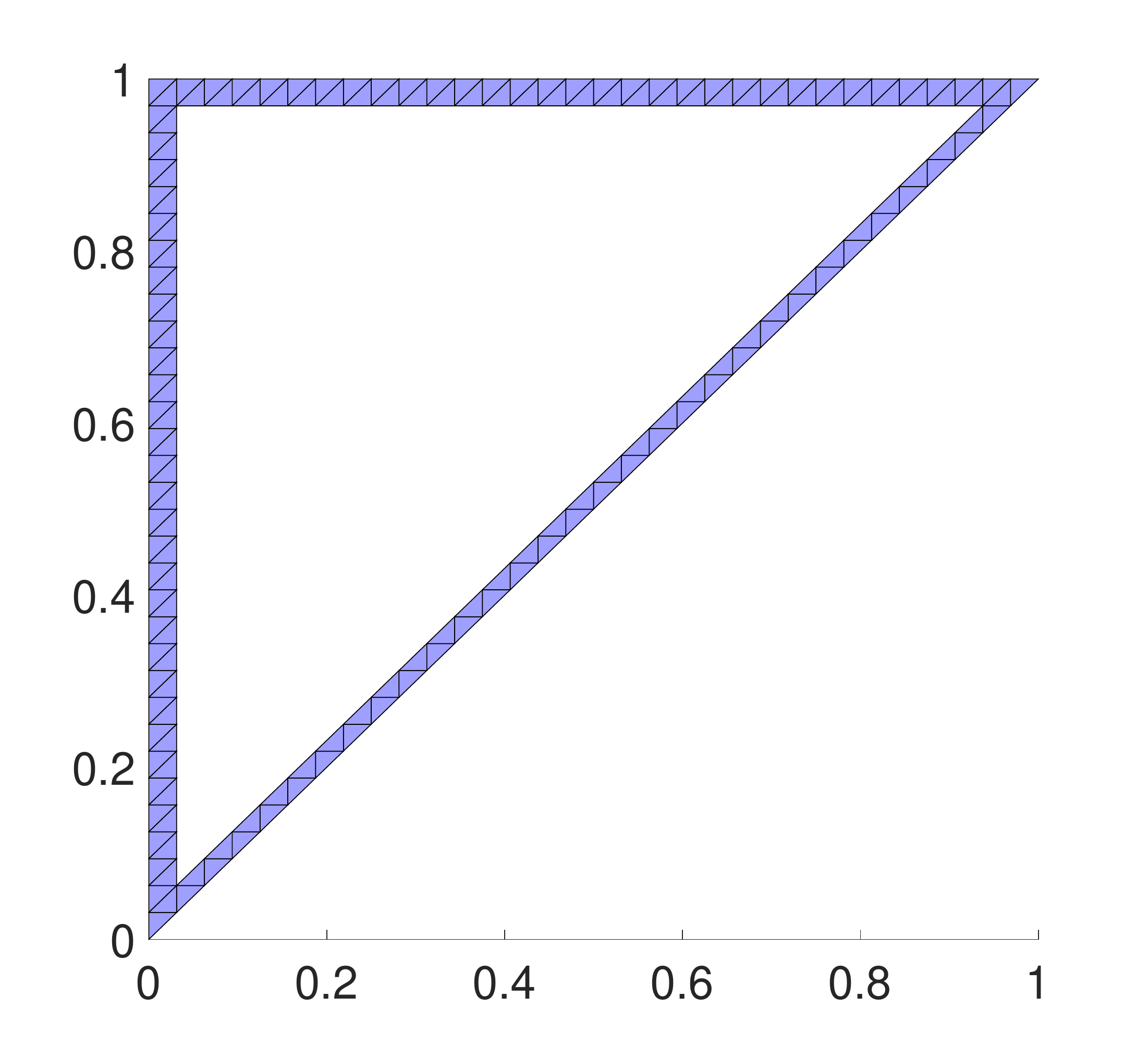}}
\subfigure[error, uniform]
{\includegraphics[trim =0.5cm 0.25cm 0.5cm 0cm, clip = true, height=4.cm]{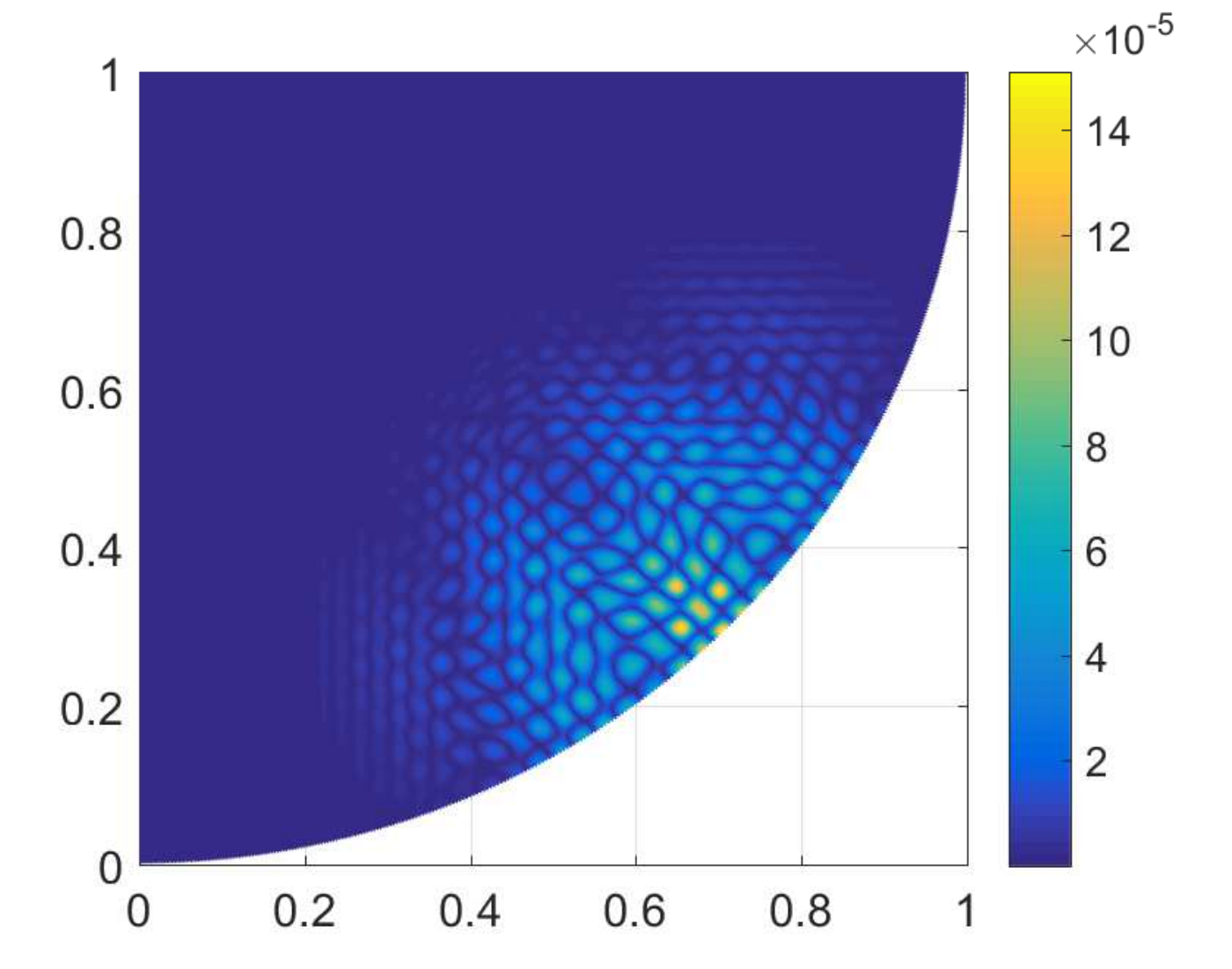}} \\
\subfigure[$\hat \Omega$, HBox-spline]
{\includegraphics[trim = 0.5cm 0.5cm 0.5cm 0.cm, clip = true, height = 4cm]
{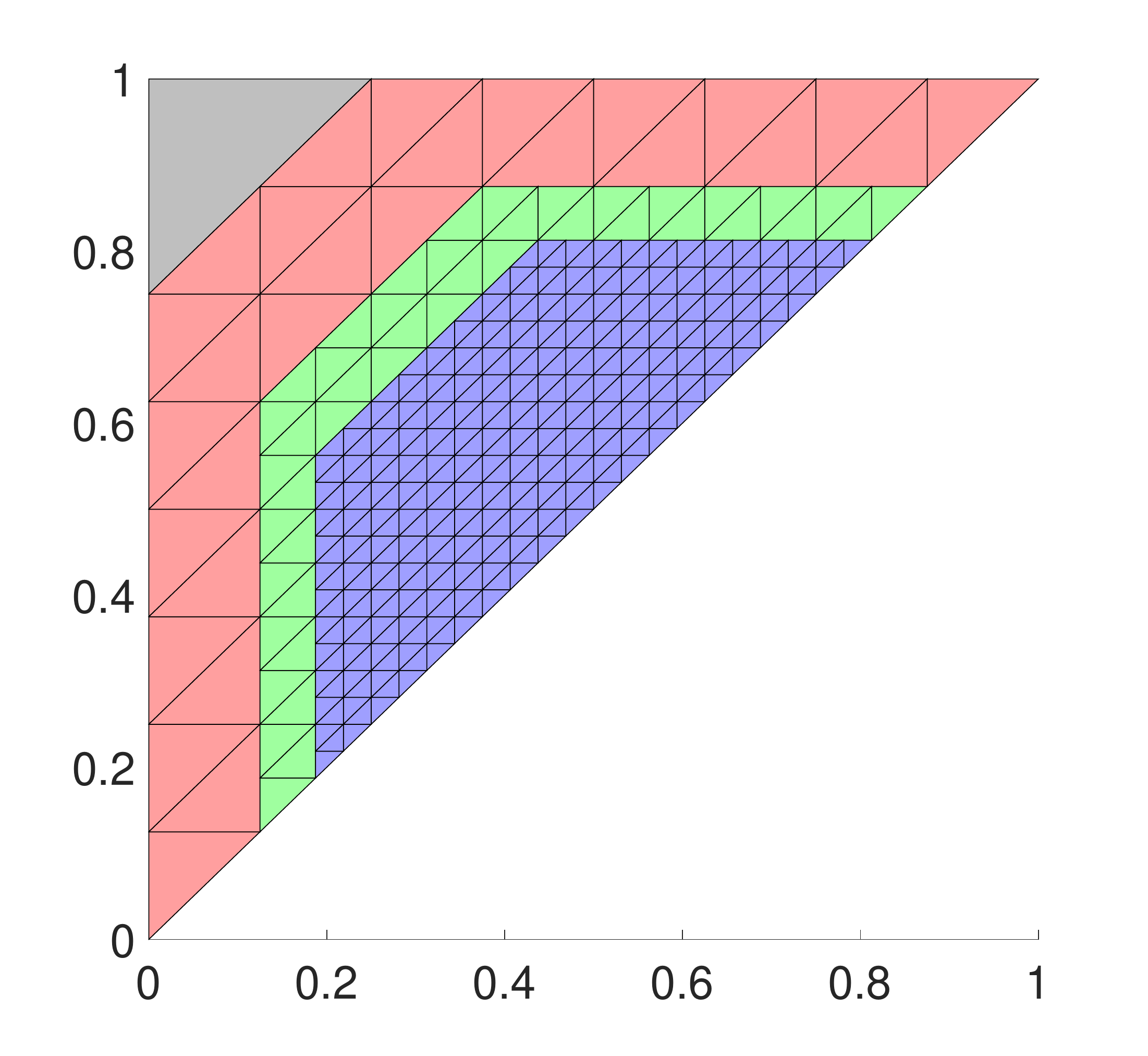}}
\subfigure[$\hat \Omega_\Gamma$, HBox-spline]
{\includegraphics[trim = 0.5cm 0.5cm 0.5cm 0.cm, clip = true, height = 4cm]
{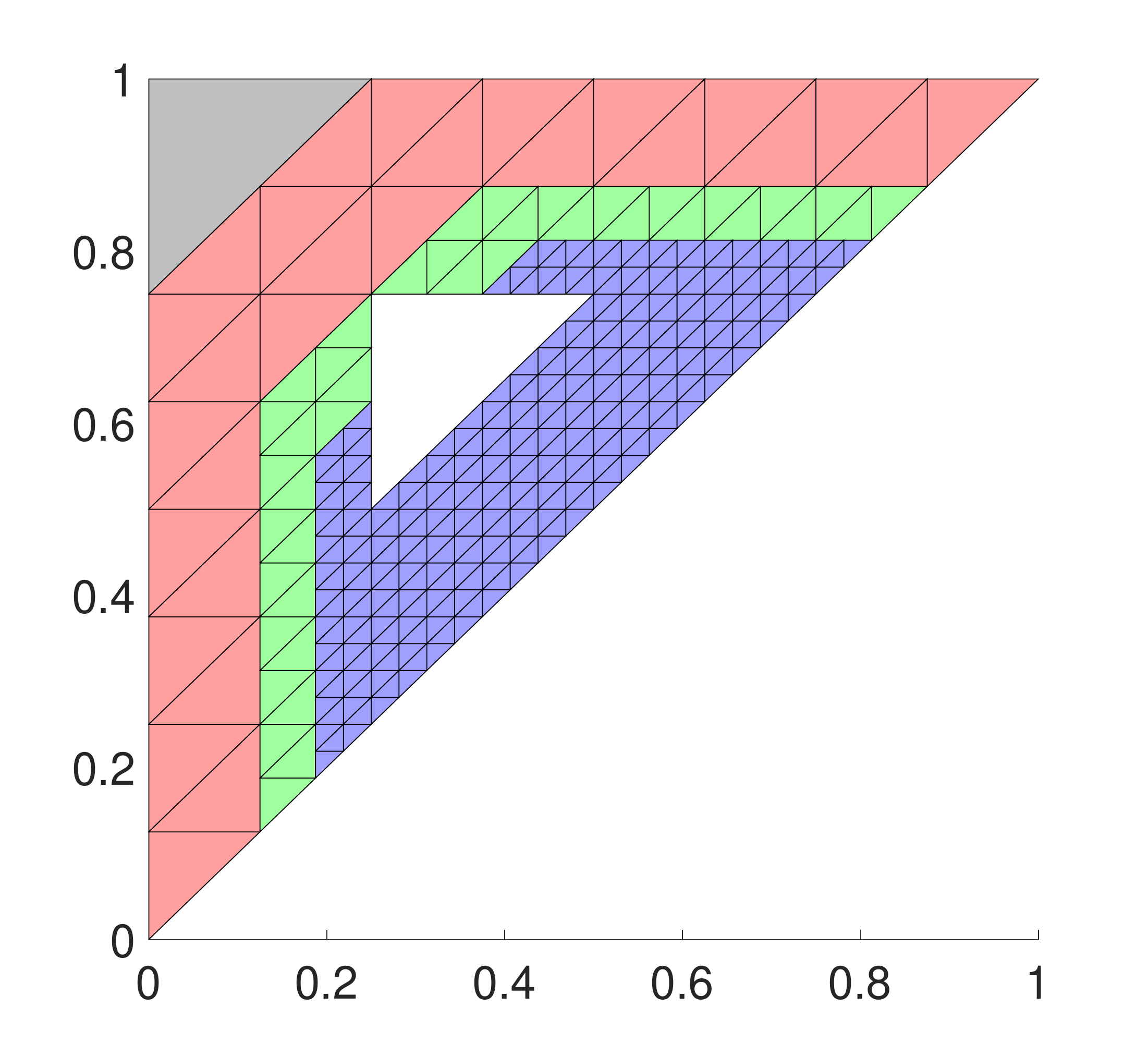}}
\subfigure[error, HBox-spline]
{\includegraphics[trim = 0.5cm 0.25cm 0.5cm 0cm, clip = true, height=4.cm]{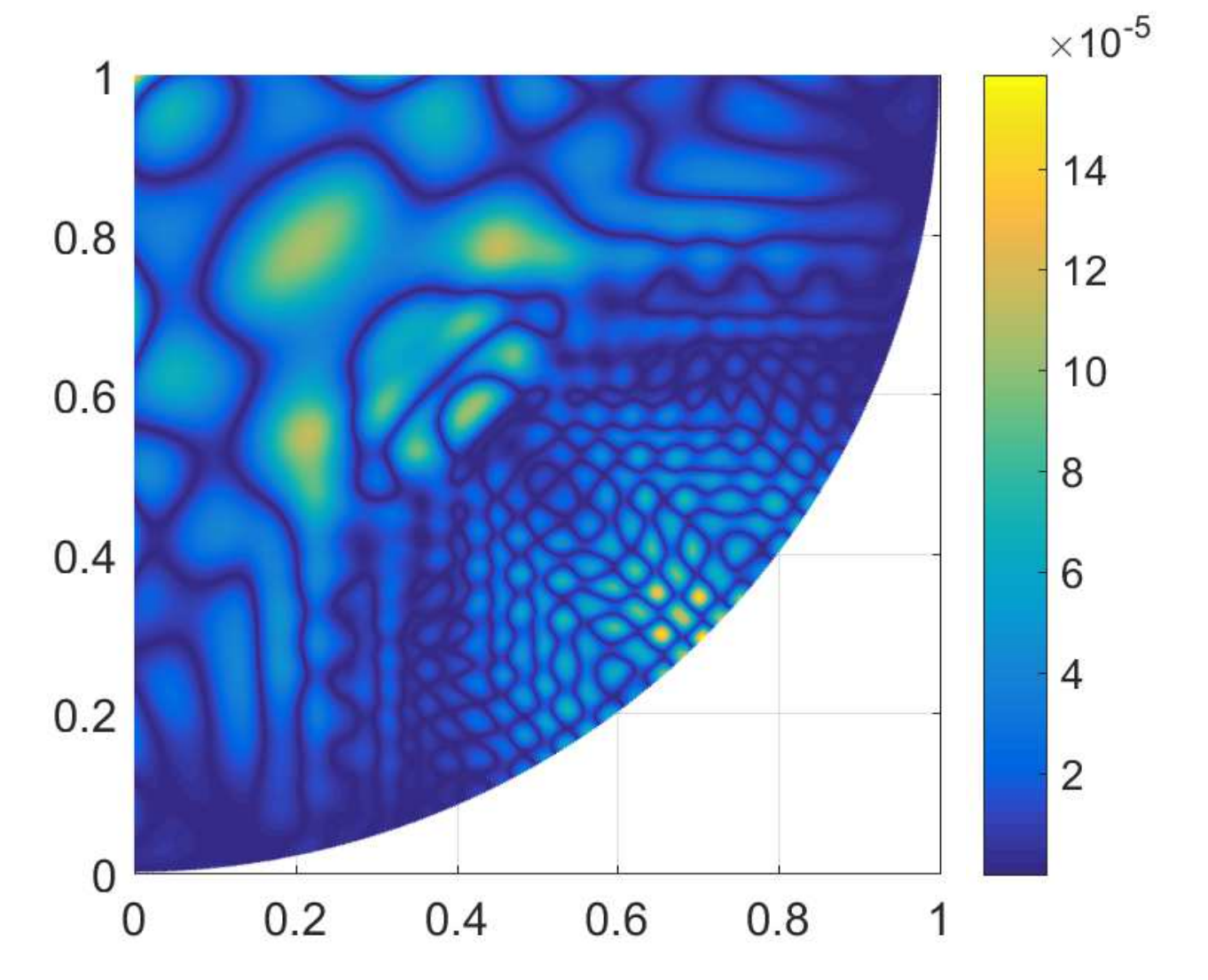}} \\
\subfigure[$\hat \Omega$, THBox-spline]
{\includegraphics[trim = 0.5cm 0.5cm 0.5cm 0.cm, clip = true, height = 4cm]
{mapped_circle_L4_omega-eps-converted-to}}
\subfigure[$\hat \Omega_\Gamma$, THBox-spline]
{\includegraphics[trim = 0.5cm 0.5cm 0.5cm 0cm, clip = true, height = 4cm]
{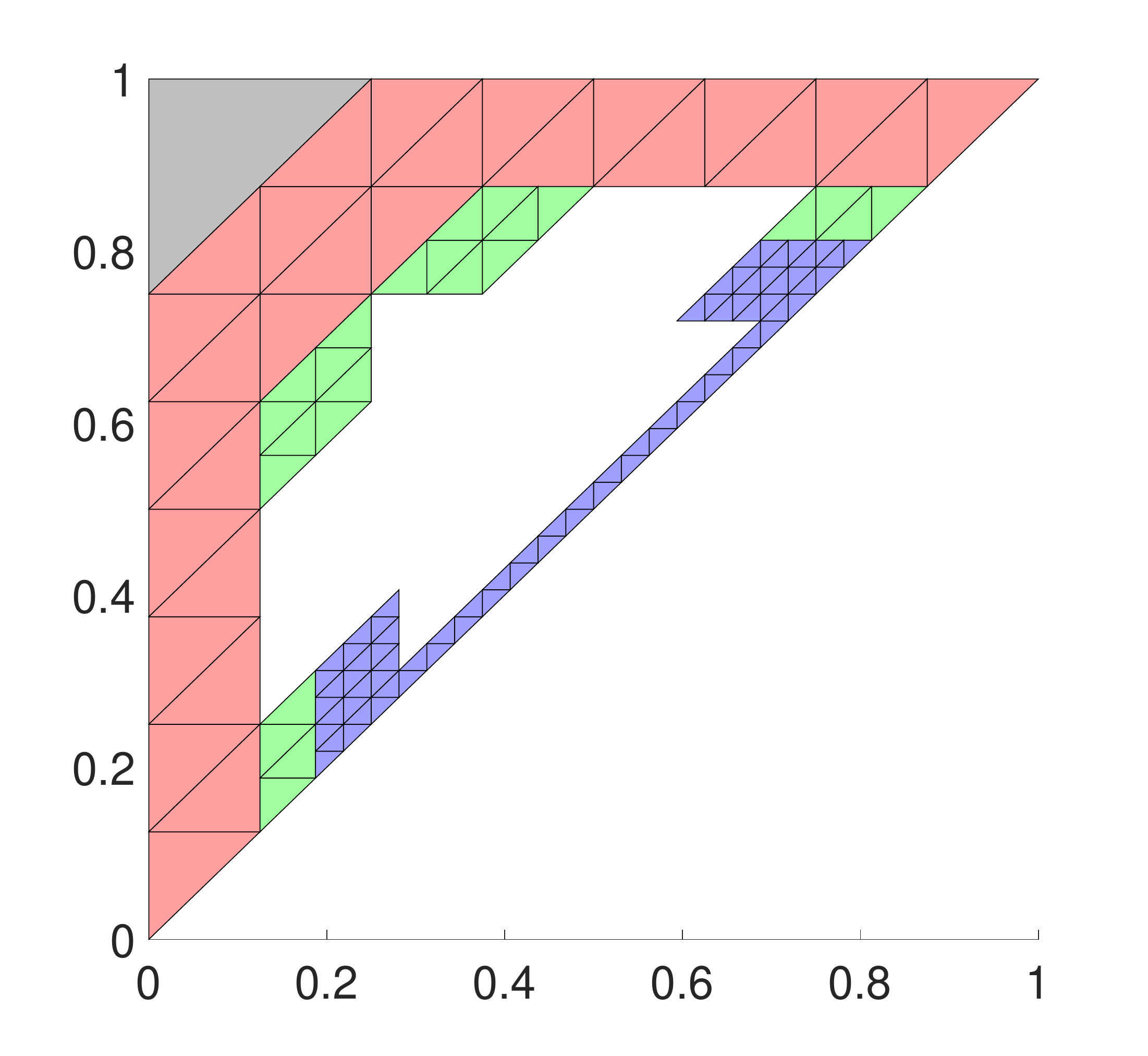}}
\subfigure[error, THBox-spline]
{\includegraphics[trim = 0.5cm 0.25cm 0.5cm 0cm, clip = true, height=4.cm]{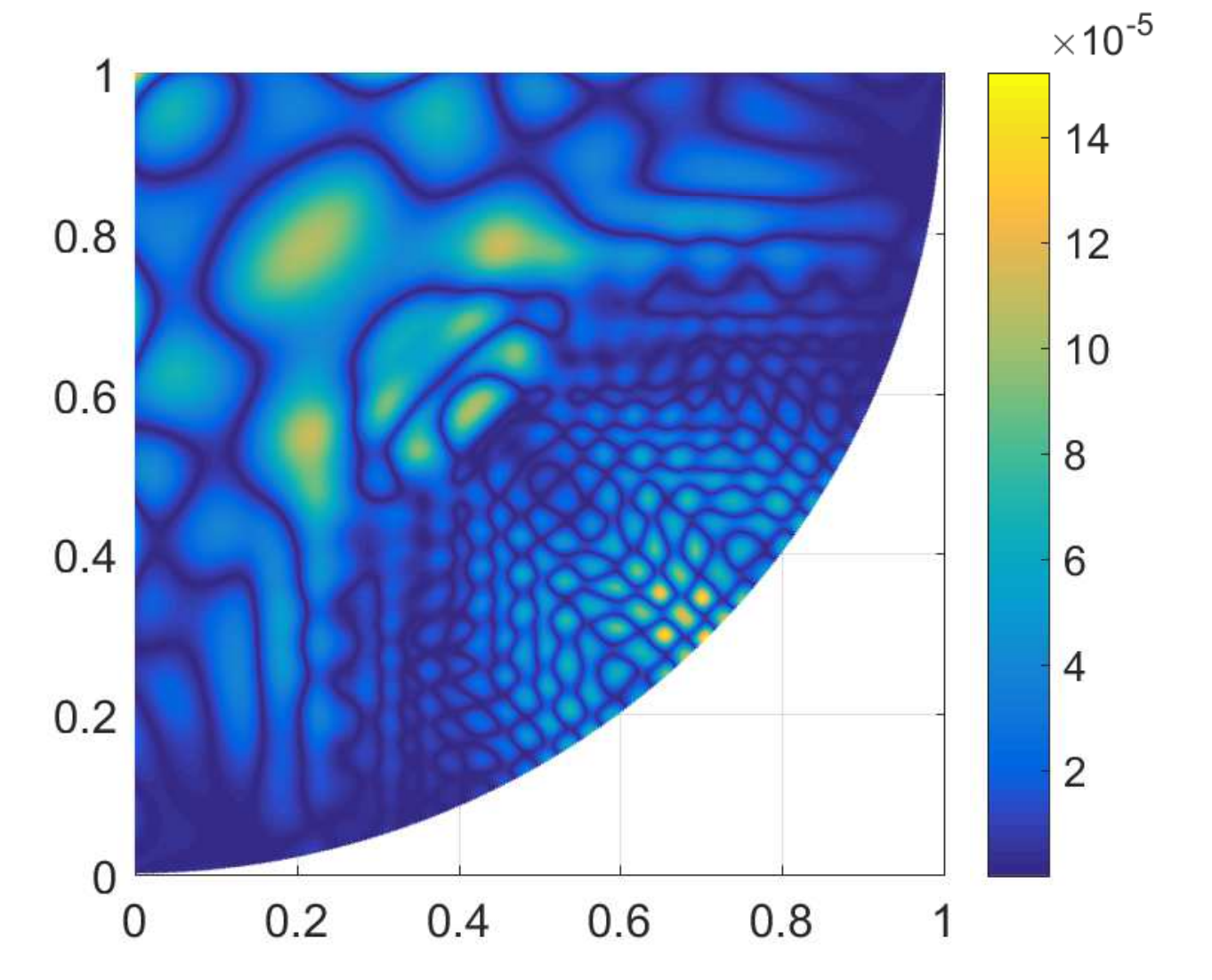}}
\caption{The circular domain problem: Uniform and hierarchical meshes on $\hat \Omega$ and $\hat \Omega_\Gamma$, together with the obtained box spline errors ($N=4$).}
\label{fig:circular-domain}
\end{center}
\end{figure}

In this example we consider the Poisson problem (\ref{problem1-poisson}) on a quarter of a circular domain. 
A unit right-angled triangular parametric domain $\hat \Omega$ (with its diagonal on the bottom right) is mapped onto a quarter of a unit circle $\Omega$ via the mapping $\bF: \hat \Omega(r,s) \to \Omega(x,y)$,
\begin{align*}
\begin{bmatrix}
x\\
y
\end{bmatrix}
=
\bF \left(\begin{bmatrix}
r\\
s
\end{bmatrix}
\right)
= \left( - r + s + \sqrt{r^2 + (s-1)^2} - \frac{2r(s-1)}{\sqrt{r^2 + (s-1)^2}} \right)
\begin{bmatrix}
r  \\
s-1
\end{bmatrix}
+
\begin{bmatrix}
0  \\
1
\end{bmatrix}.
\end{align*}
Figure~\ref{fig:circular-sol} illustrates the physical domain and the exact solution.
The exact solution $u$ (and so $g$) is chosen to have a peak on the curved boundary as
\begin{equation}\label{eq:ex-circle-g}
g(x,y) = u(x,y) = \exp \left(-25 (x - \sqrt 2/2)^2 - 25 (y - 1 + \sqrt 2/2)^2 \right),
\end{equation}
and 
\begin{align}
f(x,y) = -\Delta u(x,y) 
   &=100 \exp \left(-25 (x - \sqrt 2/2)^2 -25 (y - 1 + \sqrt 2 /2)^2 \right) \notag\\
   &\quad \cdot \left(- 25 x^2 - 25 y^2 + 25 \sqrt 2 x  + 25 (2 - \sqrt 2 )y - 49 + 25 \sqrt 2 \right).
\label{eq:ex-circle-f}
\end{align}

Figure~\ref{fig:circular-domain} shows the uniform and hierarchical three-directional meshes used for this problem with $N=4$ (left column), the different hierarchical boundary strip constructions (center column), and the corresponding error functions $|u-u_h|$ (right column).
Note that in this case the hierarchical fixed boundary strip and the HBox-spline boundary strip produce exactly the same domain $\Omega_\Gamma$. The latter strip is not shrunk because basis functions related to level 0 are active on all cells in the fixed strip.
Table~\ref{tab:circular} collects the maximum values of the error with
respect to the exact solution, for different levels of uniform and
adaptive refinement. 
As before, we observe the optimal convergence rate (order $4$) with a strong reduction of the number of degrees-of-freedom in the hierarchical cases, especially using the truncation for $\Omega_\Gamma$. 

\begin{table}[t!]
\begin{center}
\footnotesize
\begin{tabular}{c@{\hspace{.75em}}c{|}*{3}{c}
}
&& \multicolumn{3}{c}{error} 
\\
\hline
$N$ & $h_{N-1}$ & uniform & HBox &  THBox 
\\
\hline
1 & 1/4 & $1.85 \cdot 10^{-1}$  & $1.85 \cdot 10^{-1}$  & $1.85 \cdot 10^{-1}$
\\
2 & 1/8 & $4.48 \cdot 10^{-2}$  & $4.55 \cdot 10^{-2}$  & $4.49 \cdot 10^{-2}$
\\
3 & 1/16 & $4.18 \cdot 10^{-3}$  & $4.16 \cdot 10^{-3}$  & $4.09 \cdot 10^{-3}$
\\
4 & 1/32 & $1.52 \cdot 10^{-4}$  & $1.59 \cdot 10^{-4}$  & $1.52 \cdot 10^{-4}$
\end{tabular}

\medskip

\begin{tabular}{c@{\hspace{.75em}}c{|}*{3}{c}{|}*{3}{c}}
&& \multicolumn{3}{c|}{dof} & \multicolumn{3}{c}{$\Omega_\Gamma$ dof}\\
\hline
$N$ & $h_{N-1}$ & uniform & HBox &  THBox & uniform & HBox & THBox \\
\hline
1 & 1/4 & 33 & 33 &      33 & 33 & 33 & 33 \\
2 & 1/8 & 75 & 69 &      69 & 75 & 69 & 69 \\
3 & 1/16 & 207 & 111 &  111 & 171 & 111 & 104\\
4 & 1/32 & 663 & 228 &  228 & 363 & 225 & 158
\end{tabular}
\caption{The circular domain problem: Maximum error 
and degrees-of-freedom (dof), for uniform and hierarchical box splines with different types of boundary strips and different levels of refinements $N=1,\ldots,4$.}
\label{tab:circular}
\end{center}
\end{table}


\subsubsection{Mapped Z-shape domain}

\begin{figure}[!t]
\begin{center}
\subfigure[$\hat \Omega$]
{\includegraphics[trim = 0.5cm 0.5cm 0.5cm 0cm, clip = true, height = 4cm]
{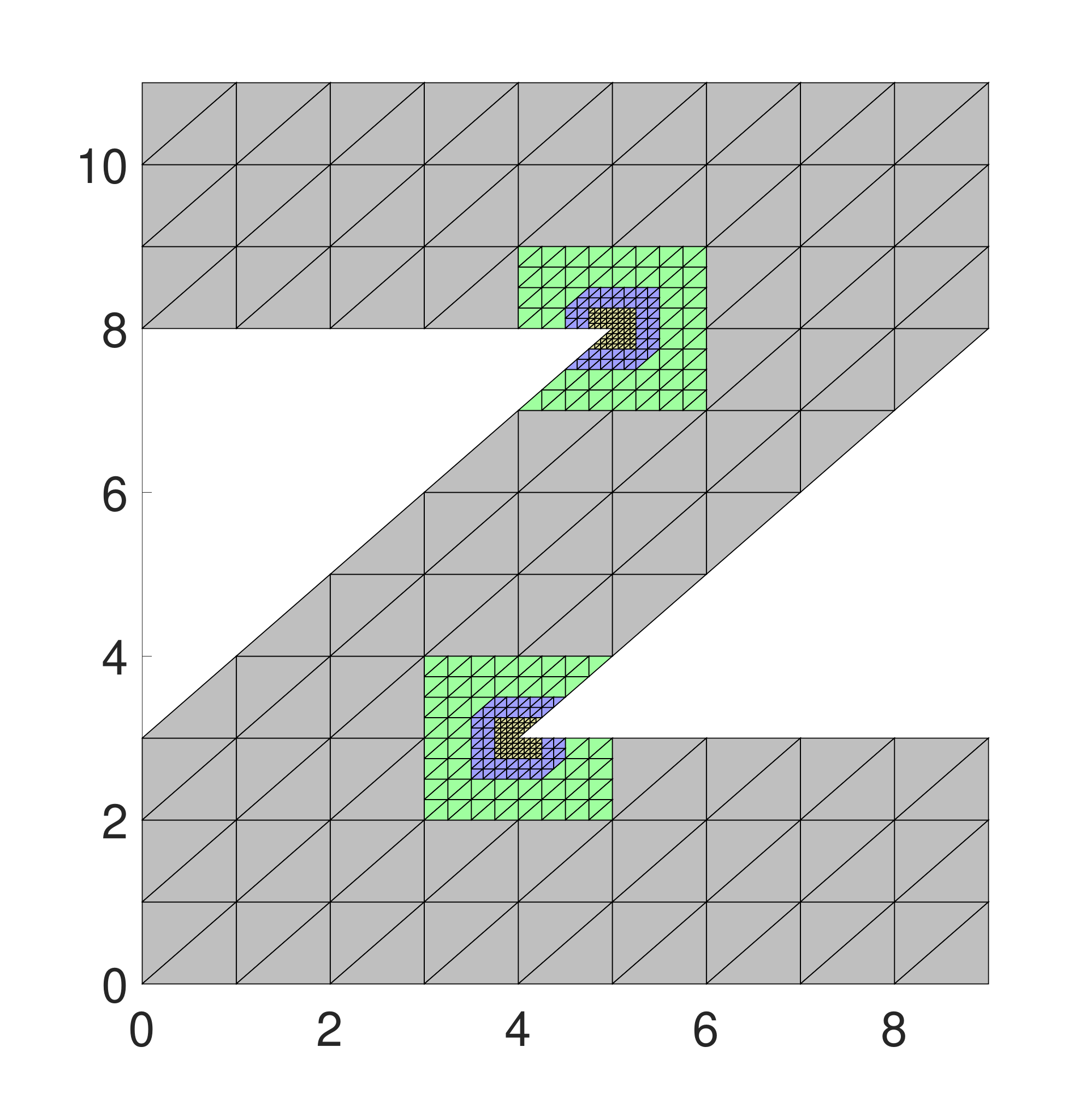}}\hspace*{1cm}
\subfigure[$\hat \Omega_\Gamma$]
{\includegraphics[trim = 0.5cm 0.5cm 0.5cm 0cm, clip = true, height = 4cm]
{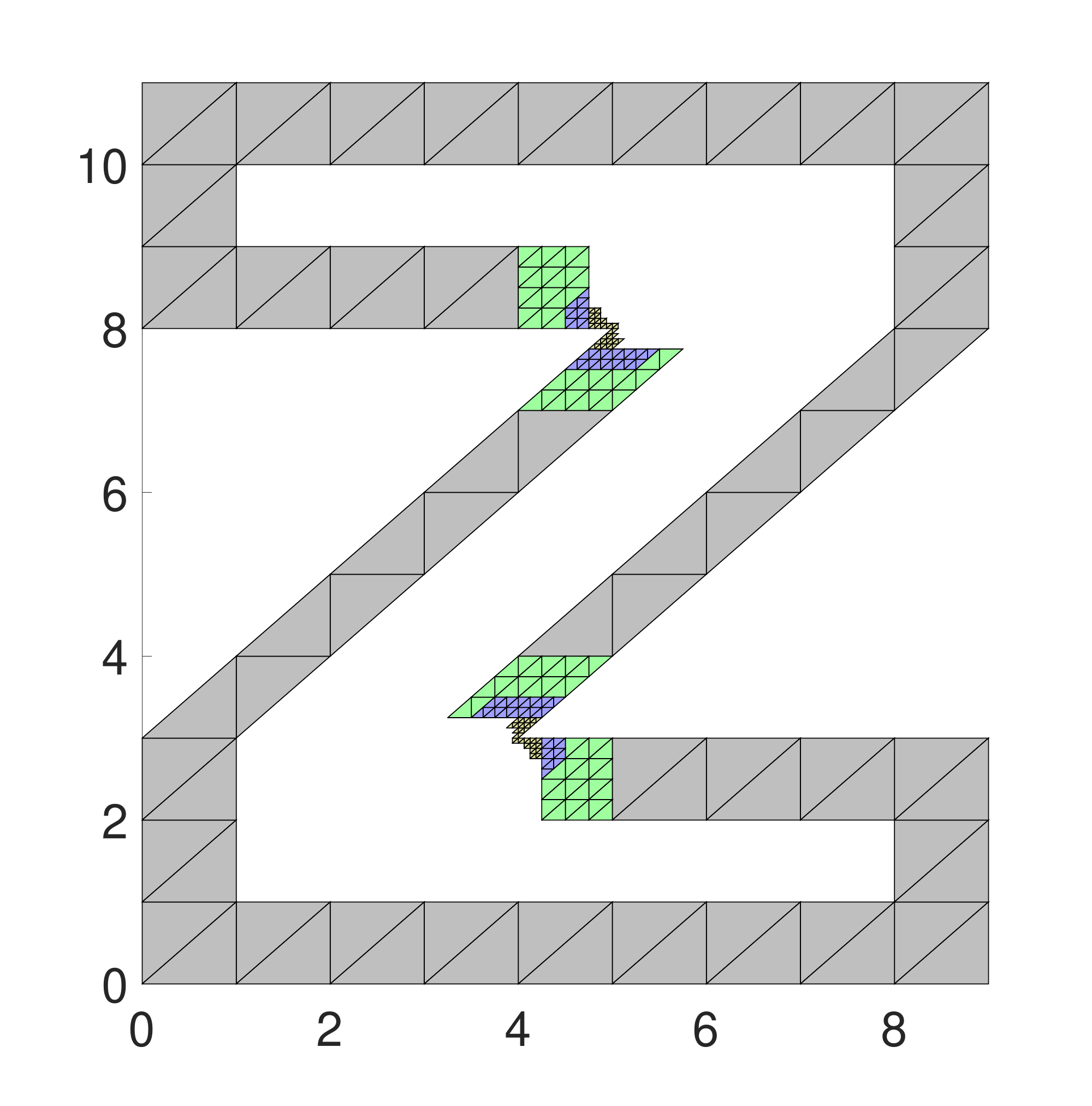}}\\
\subfigure[$\Omega$]
{\includegraphics[trim = 0.5cm 0.75cm 1cm 0.5cm, clip = true, height = 4cm]
{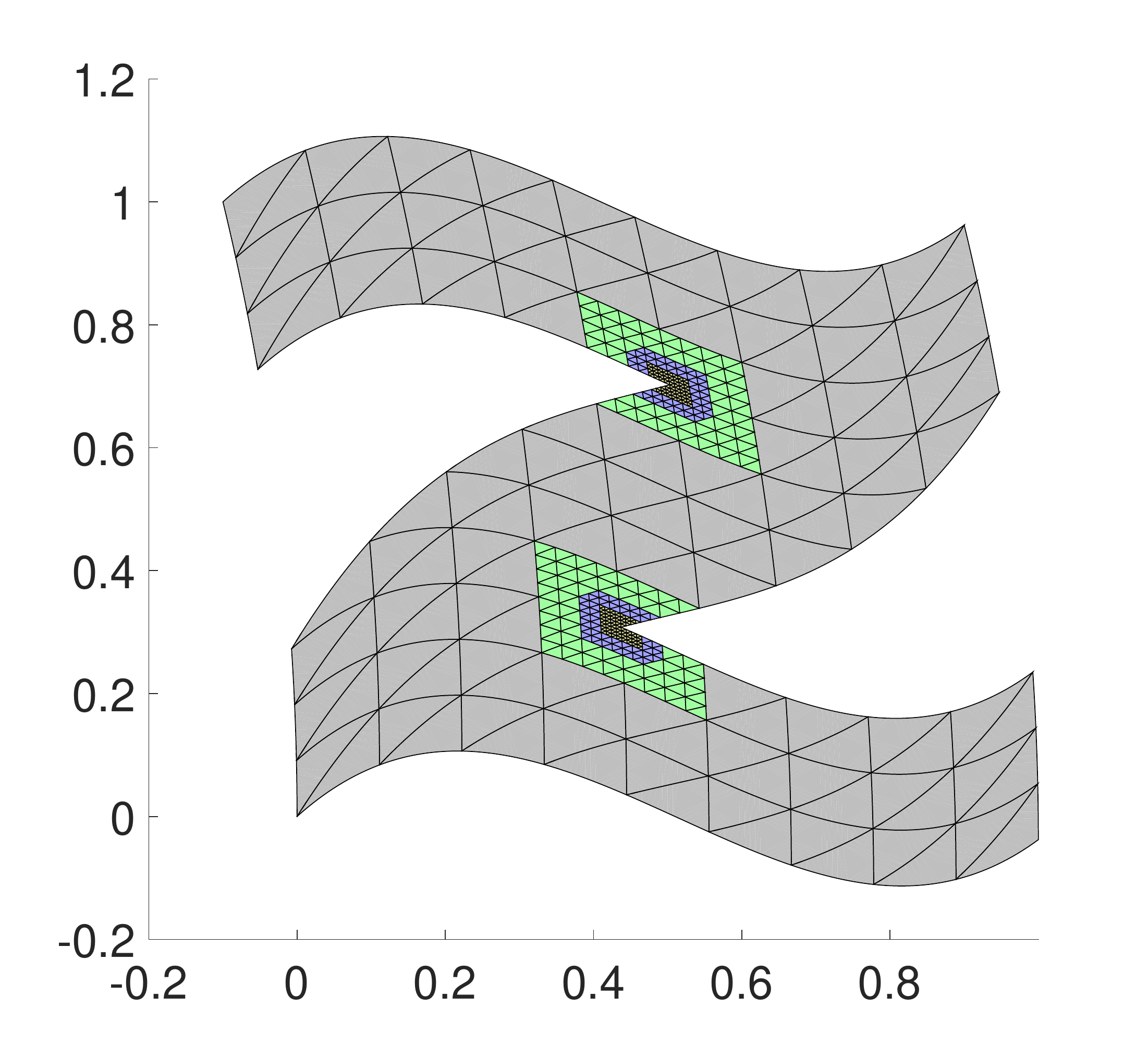}}\hspace*{.5cm}
\subfigure[solution]
{\includegraphics[trim = 0.5cm 0.5cm 0.5cm 0cm, clip = true, height = 4cm]
{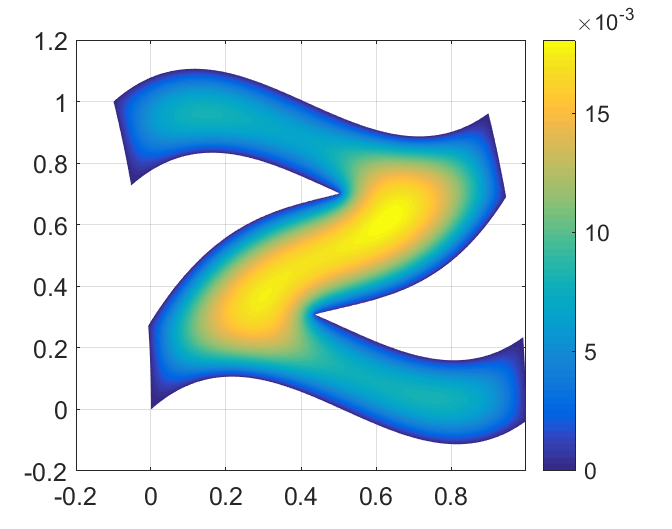}}\hspace*{.5cm}
\subfigure[error]
{\includegraphics[trim = 0.5cm 0.5cm 0.5cm 0cm, clip = true, height = 4cm]
{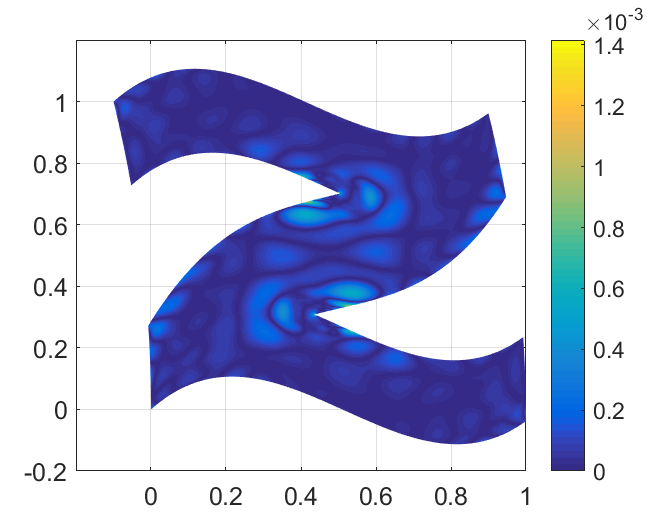}}
\caption{The Z-shaped domain problem: THBox-spline meshes on $\hat \Omega$, $\hat \Omega_\Gamma$ and $\Omega$, together with the obtained solution and its error ($N=5$).}
\label{fig:PoissonZDom}
\end{center}
\end{figure}

In this example we consider a Z-shape parametric domain $\hat \Omega \subset [0,9]\times[0,11]$ and a wave-like cubic mapping $\bF: \hat \Omega(r,s) \to \Omega(x,y)$,
\begin{equation*}
\begin{bmatrix}
x\\
y
\end{bmatrix}
=
\bF \left(
\begin{bmatrix}
r\\
s
\end{bmatrix}
\right)=
\begin{bmatrix}
r/9 - s^2/1210\\
 s/11 + (14 r^3 - 195 r^2 + 600 r)/5000
\end{bmatrix}.
\end{equation*}
The parametric and physical domains are depicted in Figure~\ref{fig:PoissonZDom}.
Then, we solve (\ref{problem1-poisson}) with 
\begin{equation*}
g(x,y) \equiv 0, \quad f(x,y) \equiv 1,
\end{equation*}
by considering uniform and hierarchical box splines up to 5 levels of refinement. The refinements were applied around the two corners of the Z-domain where the solution presents a geometric singularity (see Figure~\ref{fig:PoissonZDom}(a)).
A THBox-spline solution and its error are shown in Figure~\ref{fig:PoissonZDom}, together with the used hierarchical meshes ($N=5$). 

The exact solution is simulated by taking the uniform box spline solution computed on a very fine mesh ($N=7$).
Table~\ref{tab:PoissonZDom} reports the maximum error 
and degrees-of-freedom for all box spline solutions. 
The results of hierarchical splines using fixed boundary strips are not shown as they are identical to the ones using HBox-spline boundary strips. The latter strips are not thinner because the refinement regions are very small and basis functions related to several levels are active on each cell in the strip. Due to the nature of the considered problem, the convergence with increasing levels is slow but steady. The number of degrees-of-freedom is substantially reduced in case of hierarchical meshes compared to uniform meshes. 

\begin{table}[t!]
\begin{center}
\footnotesize
\begin{tabular}{c@{\hspace{.75em}}c{|}*{3}{c}
}
&& \multicolumn{3}{c}{error} 
 \\
\hline
$N$ & $h_{N-1}$ & uniform & HBox & THBox 
 \\
\hline
1 & 1 & $1.05 \cdot 10^{-2}$ & $1.05 \cdot 10^{-2}$ & $1.05 \cdot 10^{-2}$
\\
2 & 1/2 & $6.25 \cdot 10^{-3}$ & $8.74 \cdot 10^{-3}$ & $8.74 \cdot 10^{-3}$
\\
3 & 1/4 & $3.97 \cdot 10^{-3}$ & $4.80 \cdot 10^{-3}$ & $4.32 \cdot 10^{-3}$
\\
4 & 1/8 & $2.39 \cdot 10^{-3}$ & $2.99 \cdot 10^{-3} $ & $2.60 \cdot 10^{-3}$ 
\\
5 & 1/16 & $1.26 \cdot 10^{-3}$ & $1.66 \cdot 10^{-3} $ & $1.42 \cdot 10^{-3}$
\end{tabular}

\medskip

\begin{tabular}{c@{\hspace{.75em}}c{|}*{3}{c}{|}*{3}{c}}
&& \multicolumn{3}{c|}{dof} & \multicolumn{3}{c}{$\Omega_\Gamma$ dof}\\
\hline
$N$ & $h_{N-1}$ & uniform & HBox & THBox & uniform & HBox & THBox \\
\hline
1 & 1 & 154 & 154 & 154 & 154 & 154 & 154\\
2 & 1/2 & 451 & 156 & 156 & 380 & 180 & 180\\
3 & 1/4 & 1489 & 208 & 208 & 780 & 282 & 238\\
4 & 1/8 & 5341 & 260 & 260 & 1580 & 334 & 284\\
5 & 1/16 & 20149 & 312 
& 312 & 3180 & 386
& 336 
\end{tabular}
\caption{The Z-shaped domain problem: Maximum error 
 and degrees-of-freedom (dof), for uniform and hierarchical box splines with different types of boundary strips and different levels of refinements $N=1,\ldots,5$.}
\label{tab:PoissonZDom}
\end{center}
\end{table}

\subsubsection{Mapped triangular domain with a hole}

\begin{figure}[!t]
\begin{center}
\subfigure[$\hat \Omega$]
{\includegraphics[trim = 0.5cm 0.5cm 0.5cm 0cm, clip = true, height = 4cm]
{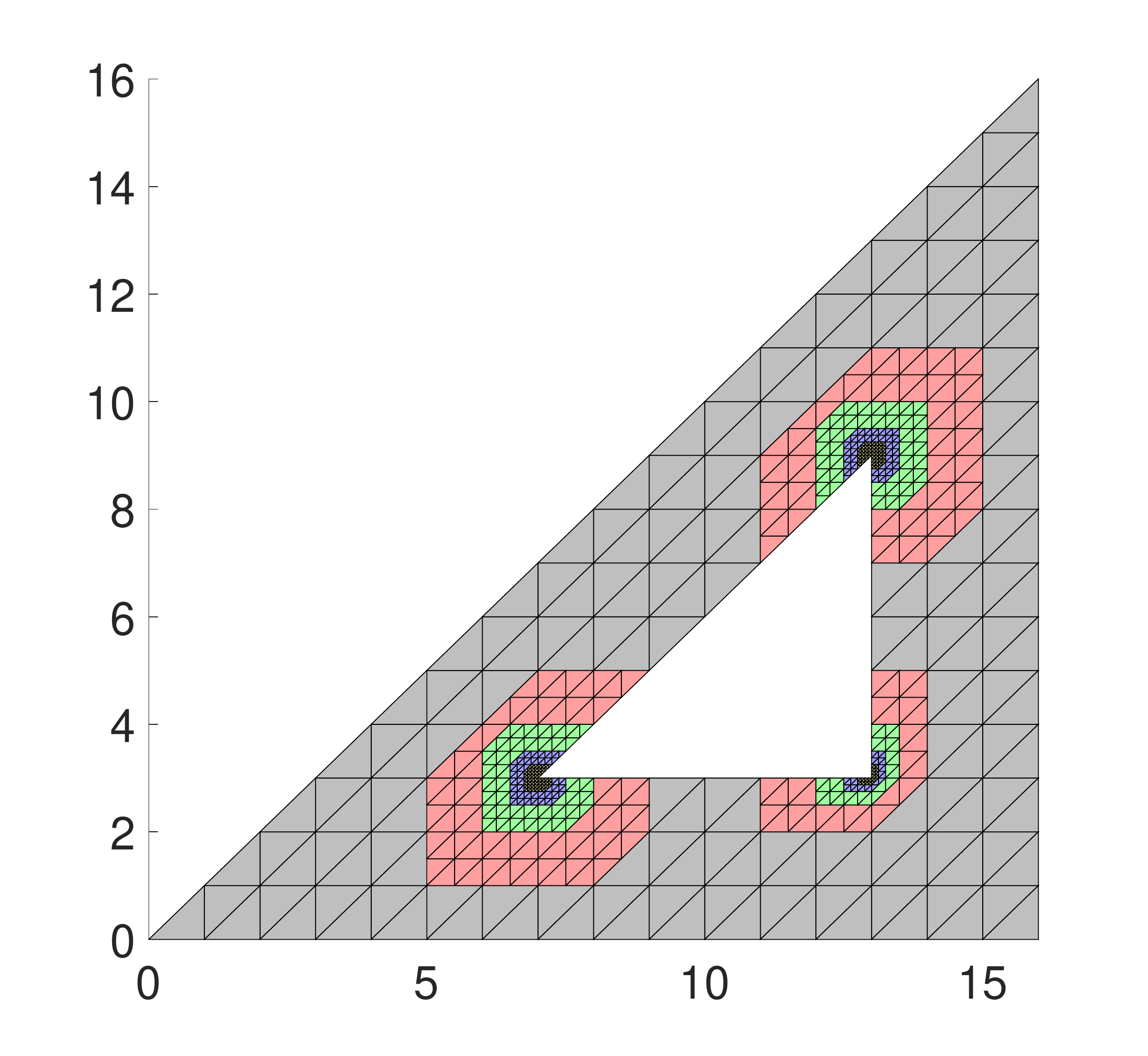}}\hspace*{1cm}
\subfigure[$\hat \Omega_\Gamma$]
{\includegraphics[trim = 0.5cm 0.5cm 0.5cm 0cm, clip = true, height = 4cm]
{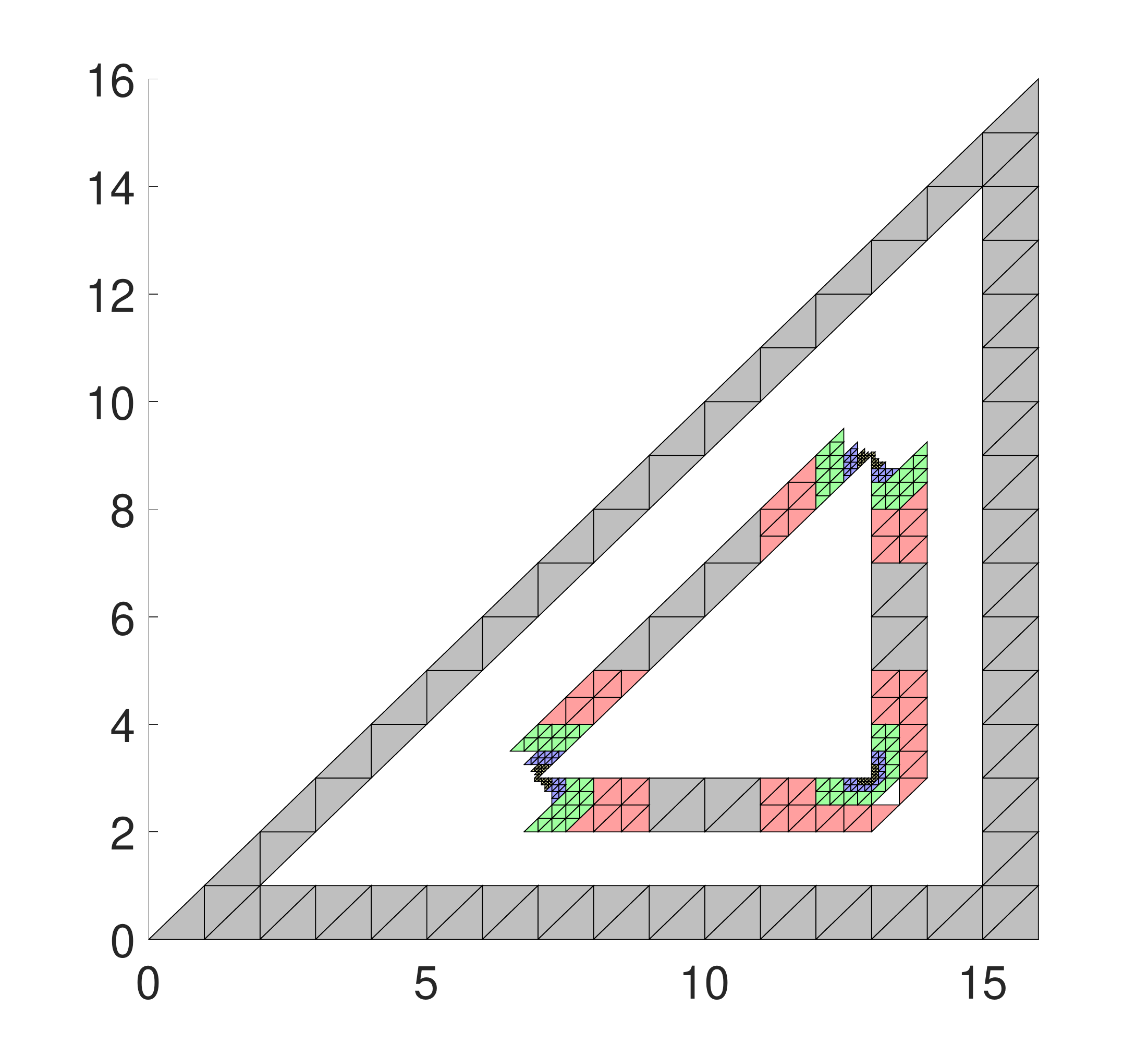}}\\
\subfigure[$\Omega$]
{\includegraphics[trim = 0.5cm 0.75cm 0.75cm 0.25cm, clip = true, height = 4cm]
{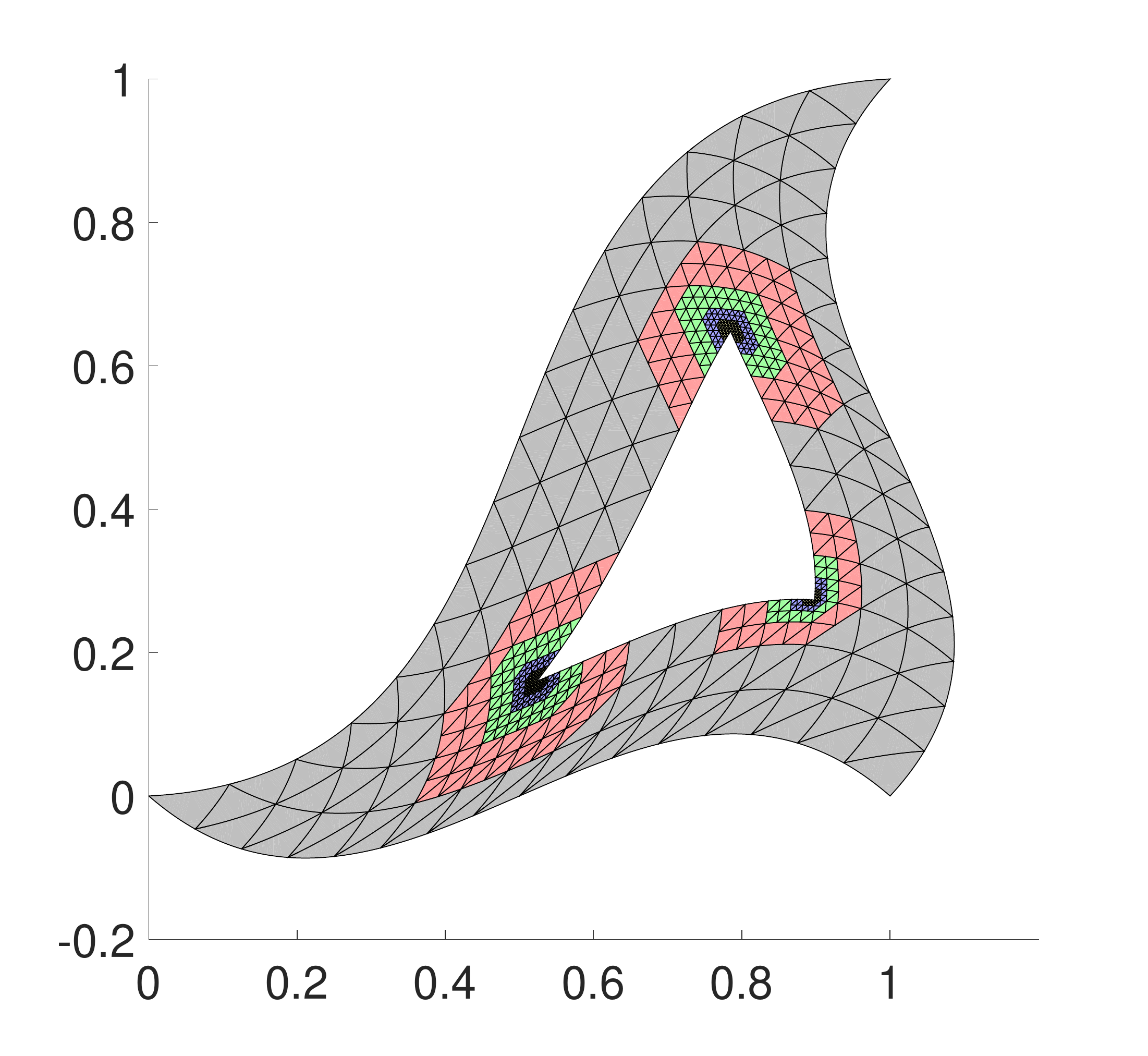}}\hspace*{.5cm}
\subfigure[solution]
{\includegraphics[trim = 0.5cm 0.5cm 0.5cm 0cm, clip = true, height = 4cm]
{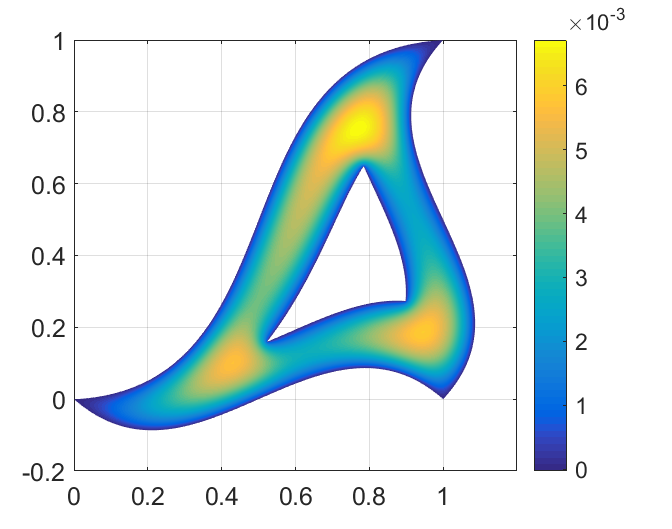}}\hspace*{.5cm}
\subfigure[error]
{\includegraphics[trim = 0.5cm 0.5cm 0.5cm 0cm, clip = true, height = 4cm]
{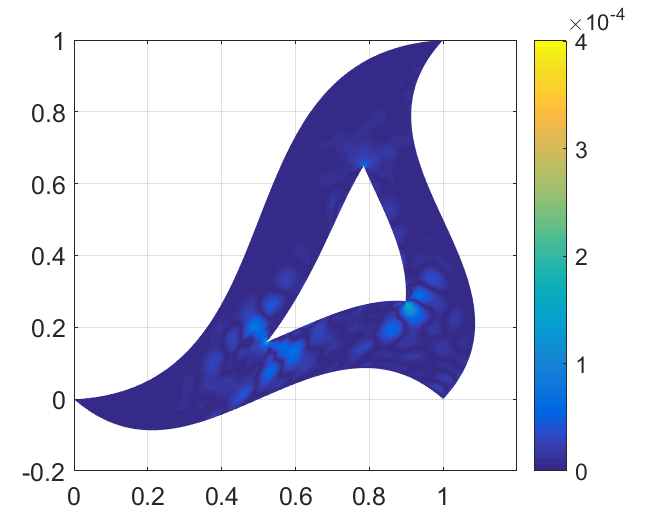}}
\caption{The triangular domain problem: THBox-spline meshes on $\hat \Omega$, $\hat \Omega_\Gamma$ and $\Omega$, together with the obtained solution and its error ($N=5$).}
\label{fig:PoissonTriangHole}
\end{center}
\end{figure}

In this example we consider a right-angled triangular domain with a hole as parametric domain $\hat \Omega \subset [0,16]^2$ and a wave-like cubic mapping $\bF: \hat \Omega(r,s) \to \Omega(x,y)$,
\begin{align*}
\begin{bmatrix}
x\\
y
\end{bmatrix}
=
\bF \left(
\begin{bmatrix}
r\\
s
\end{bmatrix}
\right)=
\begin{bmatrix}
r/16 + 9 s^3/20480 - 27 s^2/2560 + 9 s/160\\
s/16 - 9 r^3/20480 + 27 r^2/2560 - 9 r/160
\end{bmatrix}.
\end{align*}
The parametric and physical domains are depicted in Figure~\ref{fig:PoissonTriangHole}. Then, we solve the same homogeneous boundary Poisson problem as in the previous example:
\begin{equation*}
g(x,y) \equiv 0, \quad f(x,y) \equiv 1,
\end{equation*}
by considering uniform and hierarchical box splines up to 5 levels of refinement. The refinements were applied around the three interior corners of the domain (see Figure~\ref{fig:PoissonTriangHole}(a)).
A THBox-spline solution and its error are shown in Figure~\ref{fig:PoissonTriangHole}, together with the used hierarchical meshes ($N=5$). 

Similarly to the Z-shaped domain problem, the exact solution is simulated by taking the uniform box spline solution computed on a very fine mesh ($N=7$).
Table~\ref{tab:PoissonTriangHole} reports the maximum error 
 and degrees-of-freedom for all box spline solutions. 
Like in the previous example, the convergence with increasing levels is slower than the optimal rate, but the hierarchical splines require a substantially smaller number of degrees-of-freedom than the uniform counterparts. 

\begin{table}[t!]
\begin{center}
\footnotesize
\begin{tabular}{c@{\hspace{.75em}}c{|}*{3}{c}{|}*{3}{c}}
&& \multicolumn{3}{c|}{error} 
\\
\hline
$N$ & $h_{N-1}$ & uniform & HBox & THBox 
\\
\hline
1 & 1 & $3.76 \cdot 10^{-3}$ & $3.76 \cdot 10^{-3}$ & $3.76 \cdot 10^{-3}$
\\
2 & 1/2 & $2.29 \cdot 10^{-3}$ & $2.63 \cdot 10^{-3}$ & $2.58 \cdot 10^{-3}$
\\
3 & 1/4 & $1.45 \cdot 10^{-3}$ & $1.76 \cdot 10^{-3}$ & $1.56 \cdot 10^{-3}$
 \\
4 & 1/8  & $8.20 \cdot 10^{-4}$ & $1.05 \cdot 10^{-3}$ & $9.03 \cdot 10^{-4}$
\\
5 & 1/16 & $3.41 \cdot 10^{-4}$ & $4.82 \cdot 10^{-4}$ & $4.01 \cdot 10^{-4}$
\end{tabular}

\medskip

\begin{tabular}{c@{\hspace{.75em}}c{|}*{3}{c}{|}*{3}{c}}
&& \multicolumn{3}{c|}{dof} & \multicolumn{3}{c}{$\Omega_\Gamma$ dof}\\
\hline
$N$ & $h_{N-1}$ & uniform & HBox & THBox & uniform & HBox & THBox \\
\hline
1 & 1 & 206 & 206 & 206 & 206 & 206 & 206\\
2 & 1/2 & 635 & 260 & 260 & 516 & 284 & 282\\
3 & 1/4 & 2153 & 320 & 320 & 1044 & 386 & 346\\
4 & 1/8 & 7829 & 380 & 380 & 2100 & 446 & 400\\
5 & 1/16 & 29741 & 440 & 440 & 4212 & 506 & 454
\end{tabular}
\caption{The triangular domain problem: Maximum error 
 and degrees-of-freedom (dof), for uniform and hierarchical box splines with different types of boundary strips and different levels of refinements $N=1,\ldots,5$.}
\label{tab:PoissonTriangHole}
\end{center}
\end{table}

%
\subsection{Advection-diffusion problem on unit square} \label{sec:num-advection}

\begin{figure}[t!]
\begin{center}
\subfigure[domain]{
\begin{picture}(160,120)(0,-10)
\put(0,0){\includegraphics[trim = 4.5cm 9cm 4cm 9.1cm,  clip = true, height = 4cm]{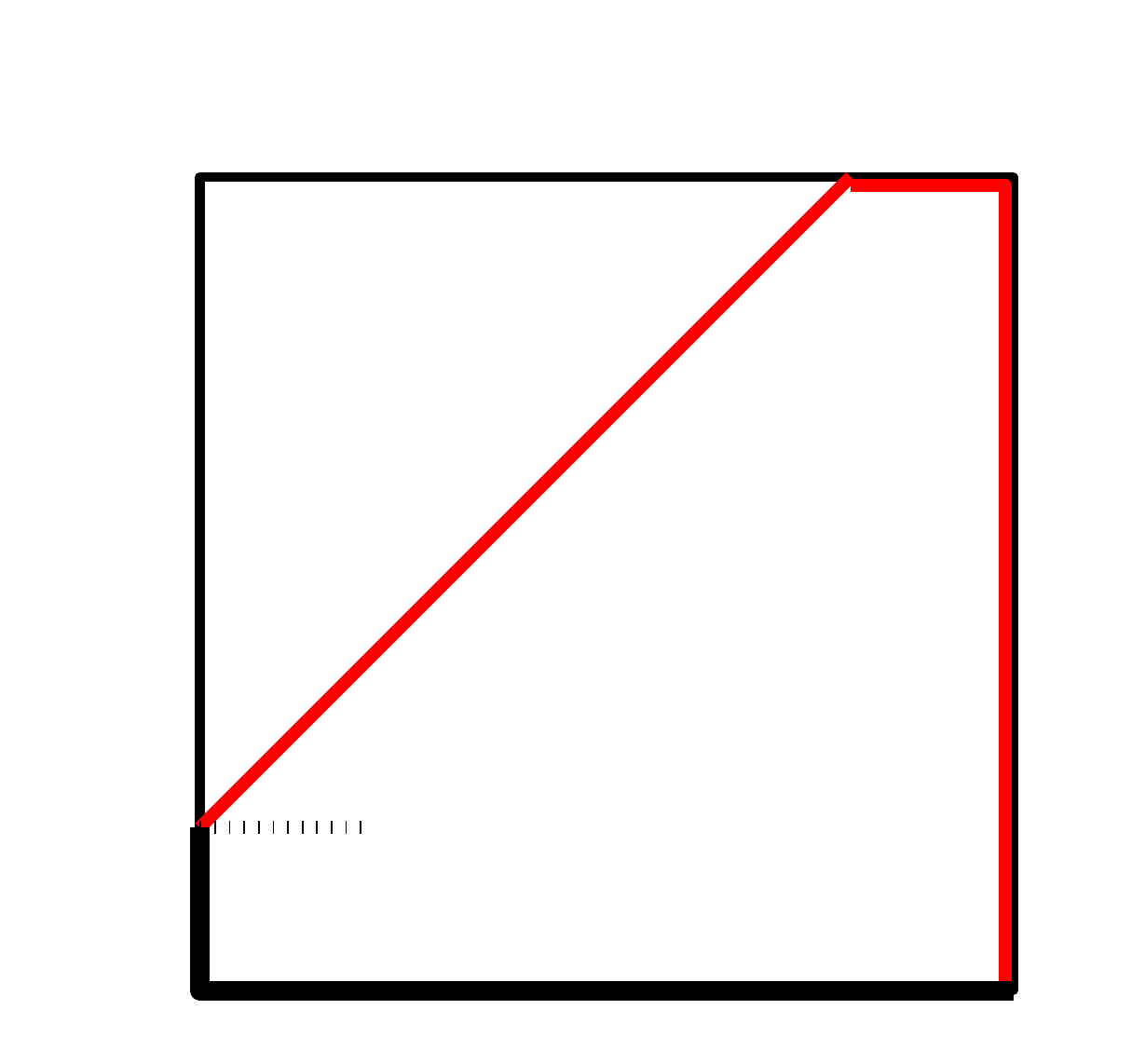}}
\put(115,50){$u=0$}%
\put(-9,60){$u=0$}%
\put(-9,15){$u=1$}%
\put(58,-5){$u=1$} %
\put(58,101){$u=0$}%
\put(35,30){$\theta$}
\end{picture}}
\subfigure[exact solution]
{\includegraphics[trim = 0.87cm 0.5cm 0.77cm 0.75cm, clip=true, height = 4cm]
{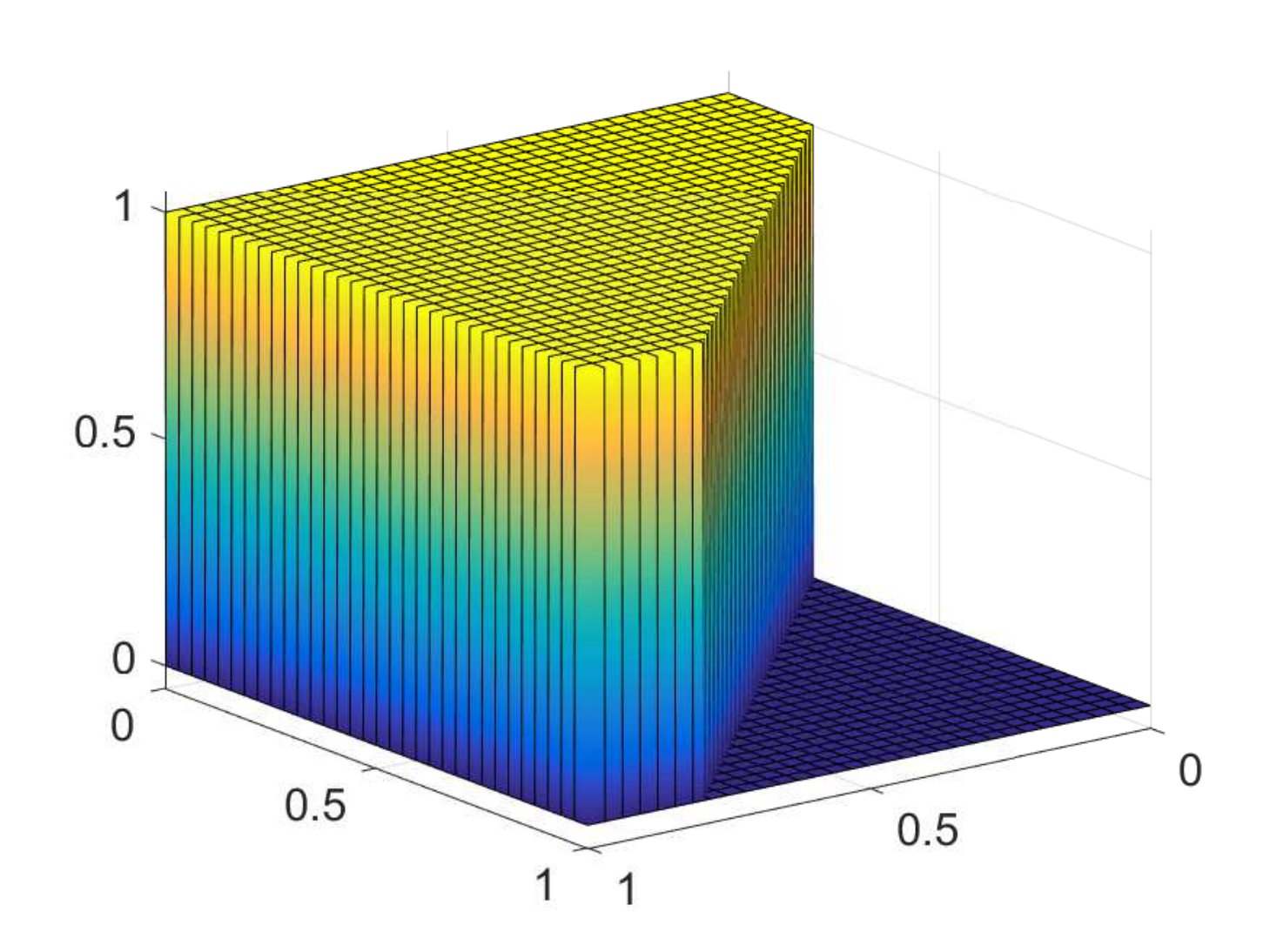}}
\caption{Advection-diffusion problem: Domain with Dirichlet boundary conditions and exact solution.}
\label{fig:advection}
\end{center}
\end{figure}

\begin{figure}[!t]
\begin{center}
\subfigure[$\Omega$, uniform]
{\includegraphics[trim = 2cm 0.75cm 1.75cm 1cm, clip = true, height = 4cm]
{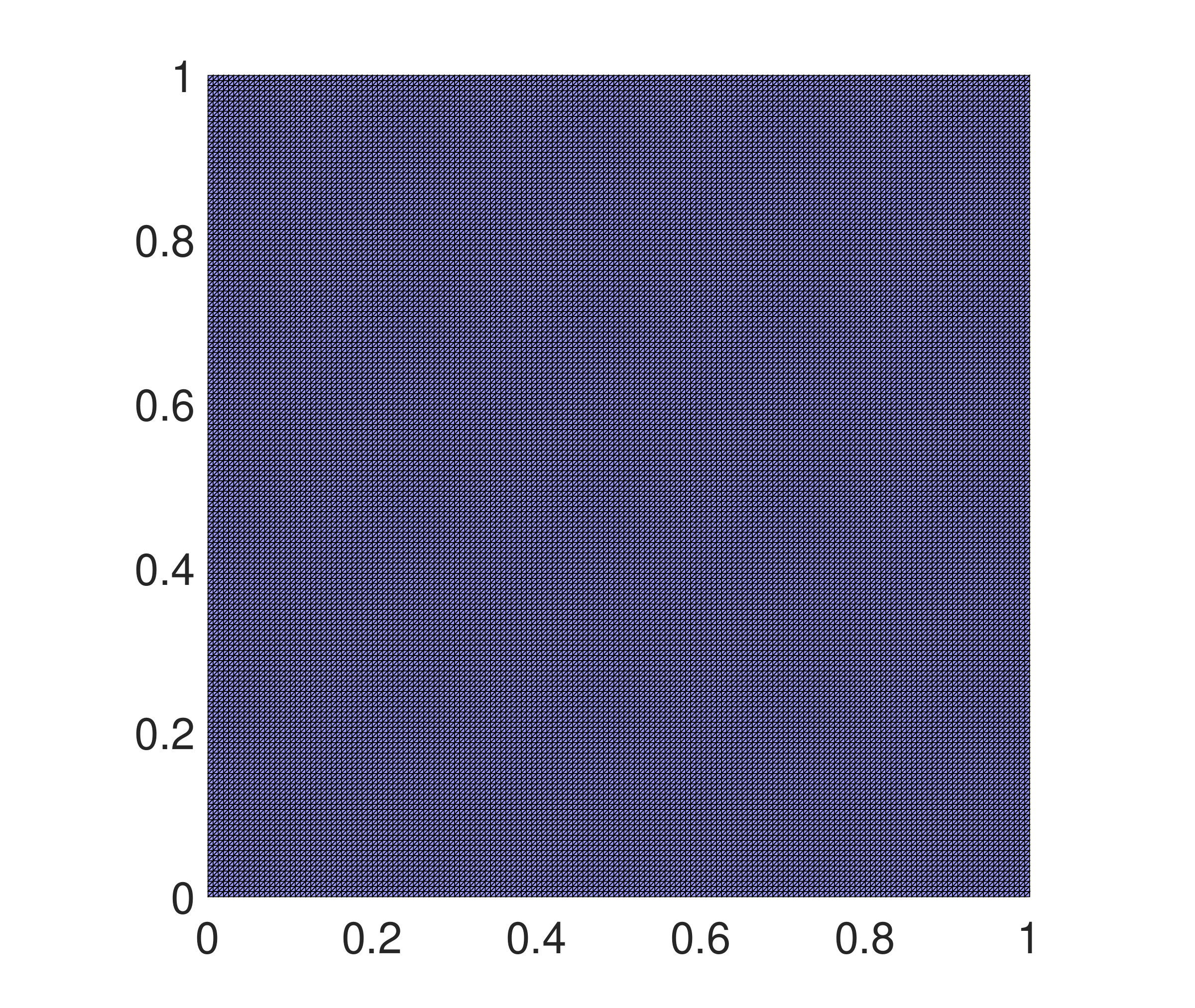}}
\subfigure[$\Omega_\Gamma$, uniform]
{\includegraphics[trim = 2cm 0.75cm 1.75cm 1cm, clip = true, height = 4cm]
{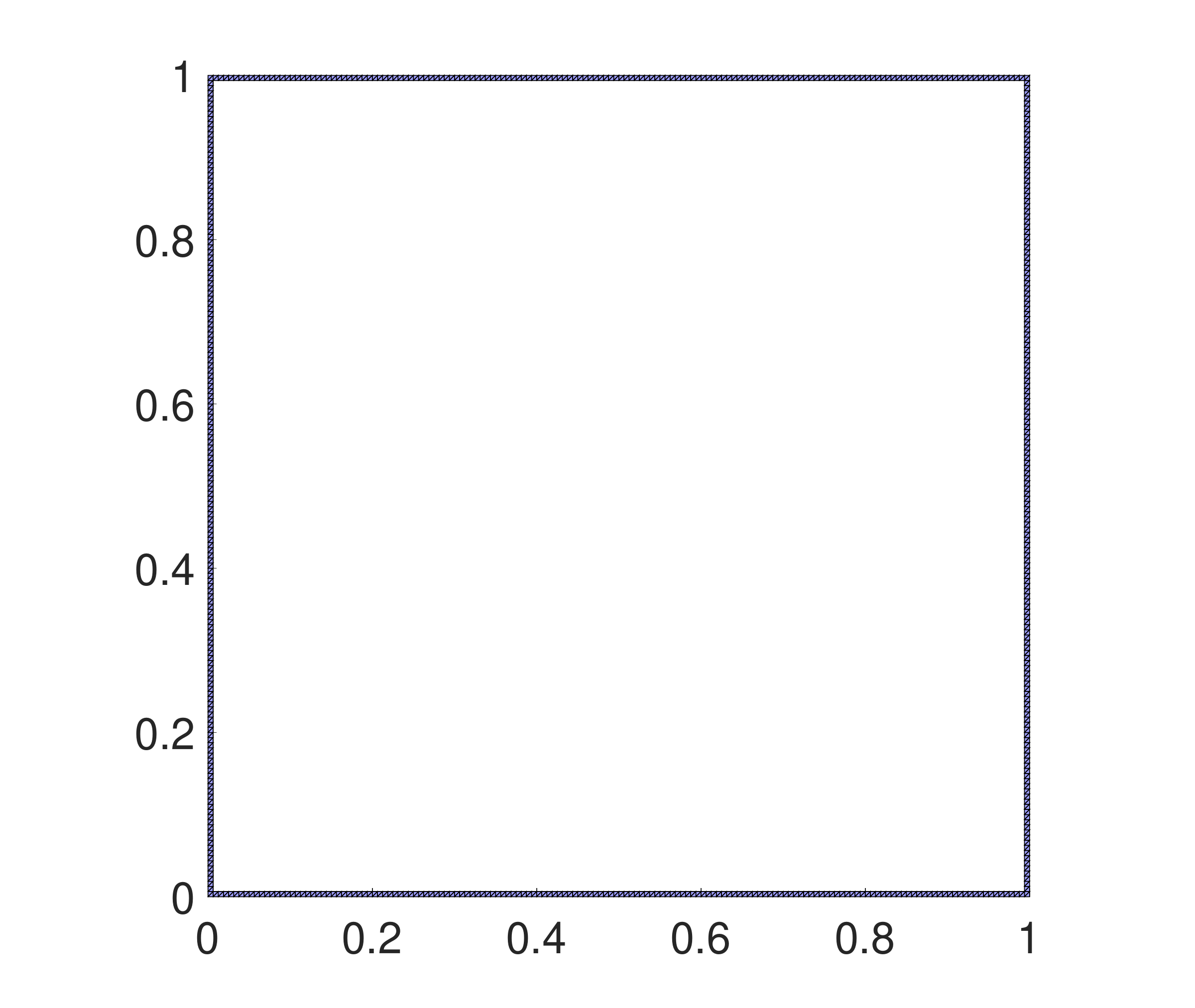}}
\subfigure[solution, uniform]
{\includegraphics[trim = 0.87cm 0.5cm 0.77cm 0.75cm, clip = true, height=3.5cm]
{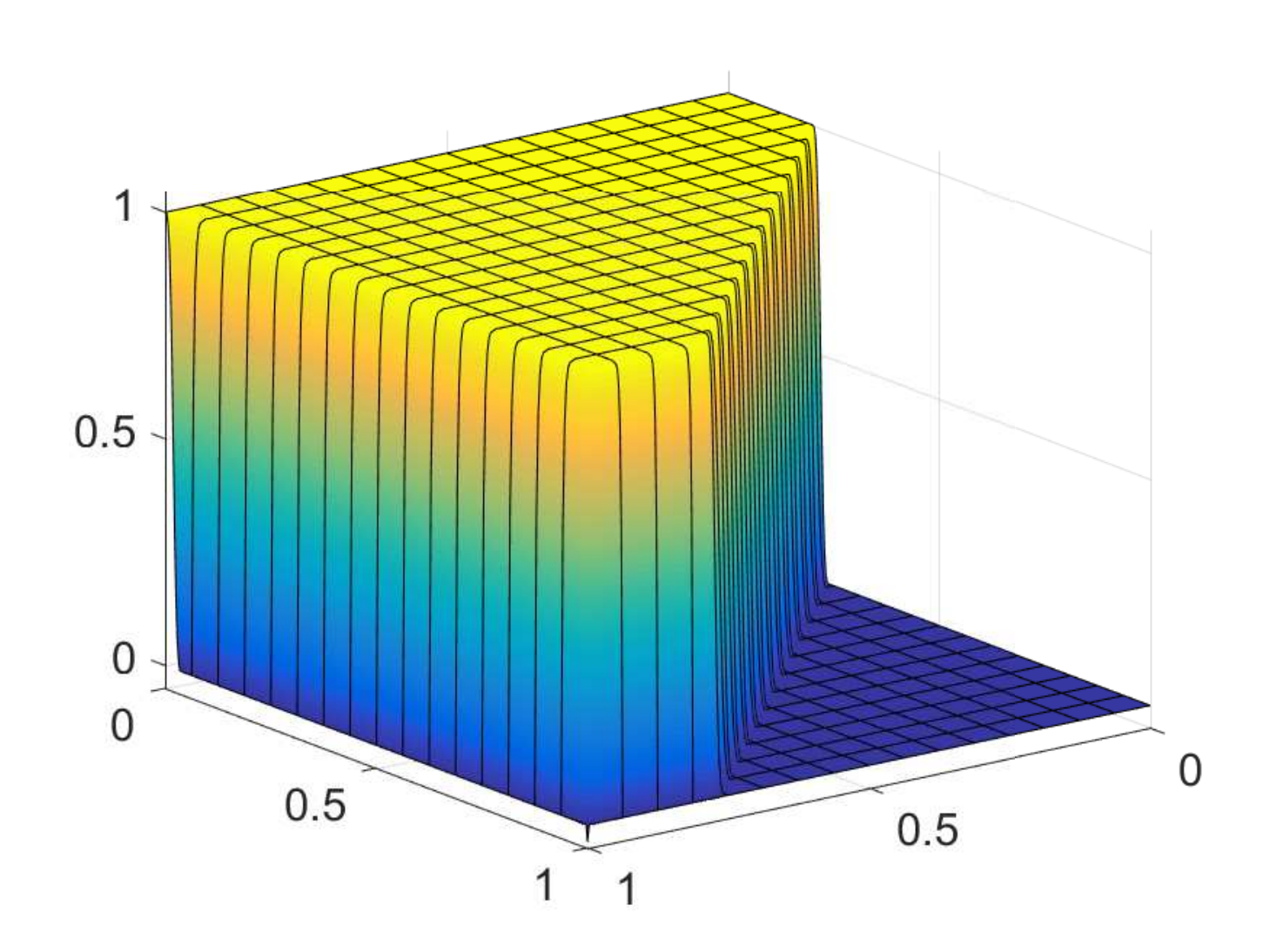}} \\
\subfigure[$\Omega$, HBox-spline]
{\includegraphics[trim = 2cm 0.75cm 1.75cm 1.cm, clip = true, height = 4cm]
{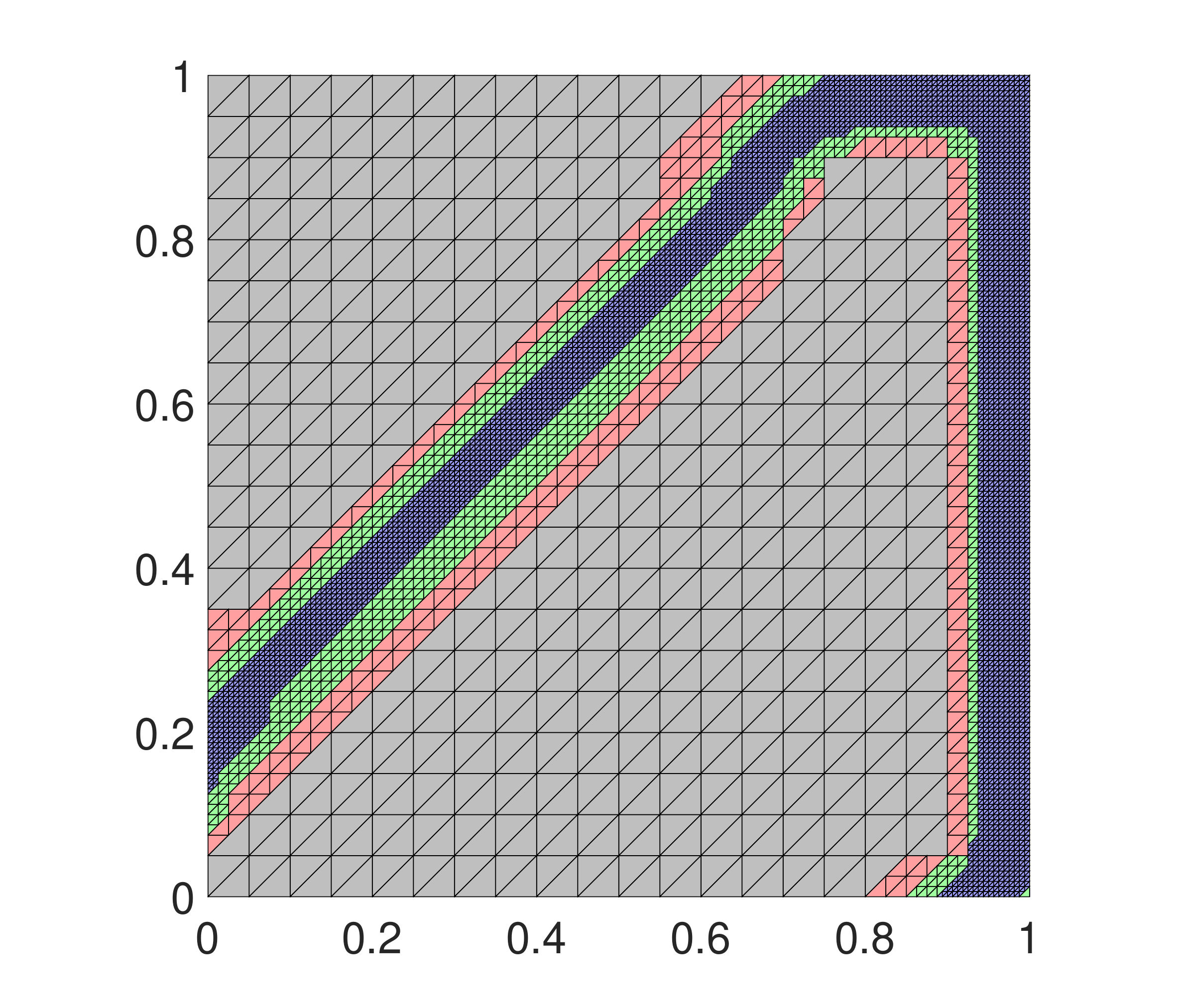}}
\subfigure[$\Omega_\Gamma$, HBox-spline]
{\includegraphics[trim = 2cm 0.75cm 1.75cm 1.cm, clip = true, height = 4cm]
{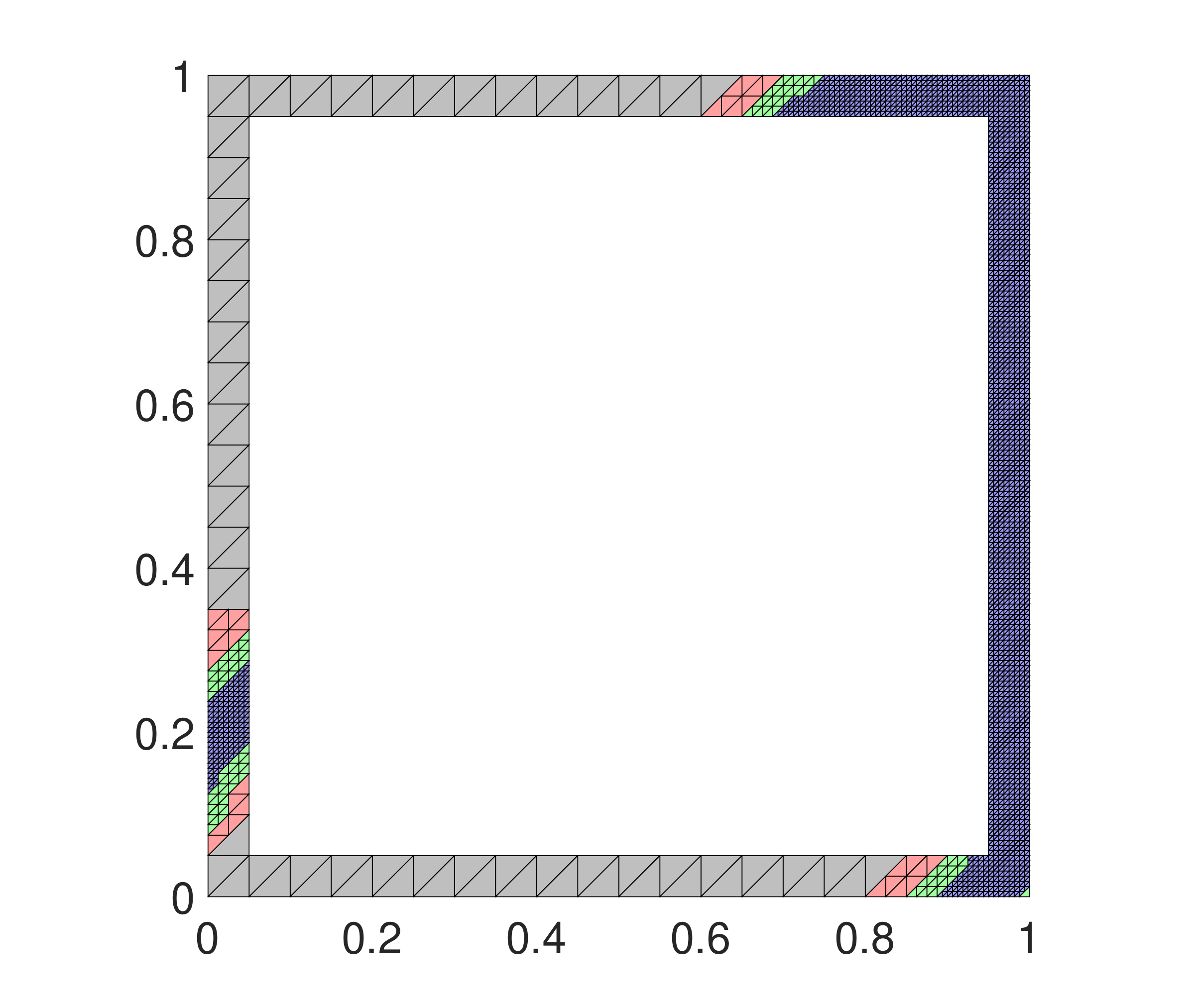}}
\subfigure[solution, HBox-spline]
{\includegraphics[trim = 0.87cm 0.5cm 0.77cm 0.75cm, clip = true, height=3.5cm]
{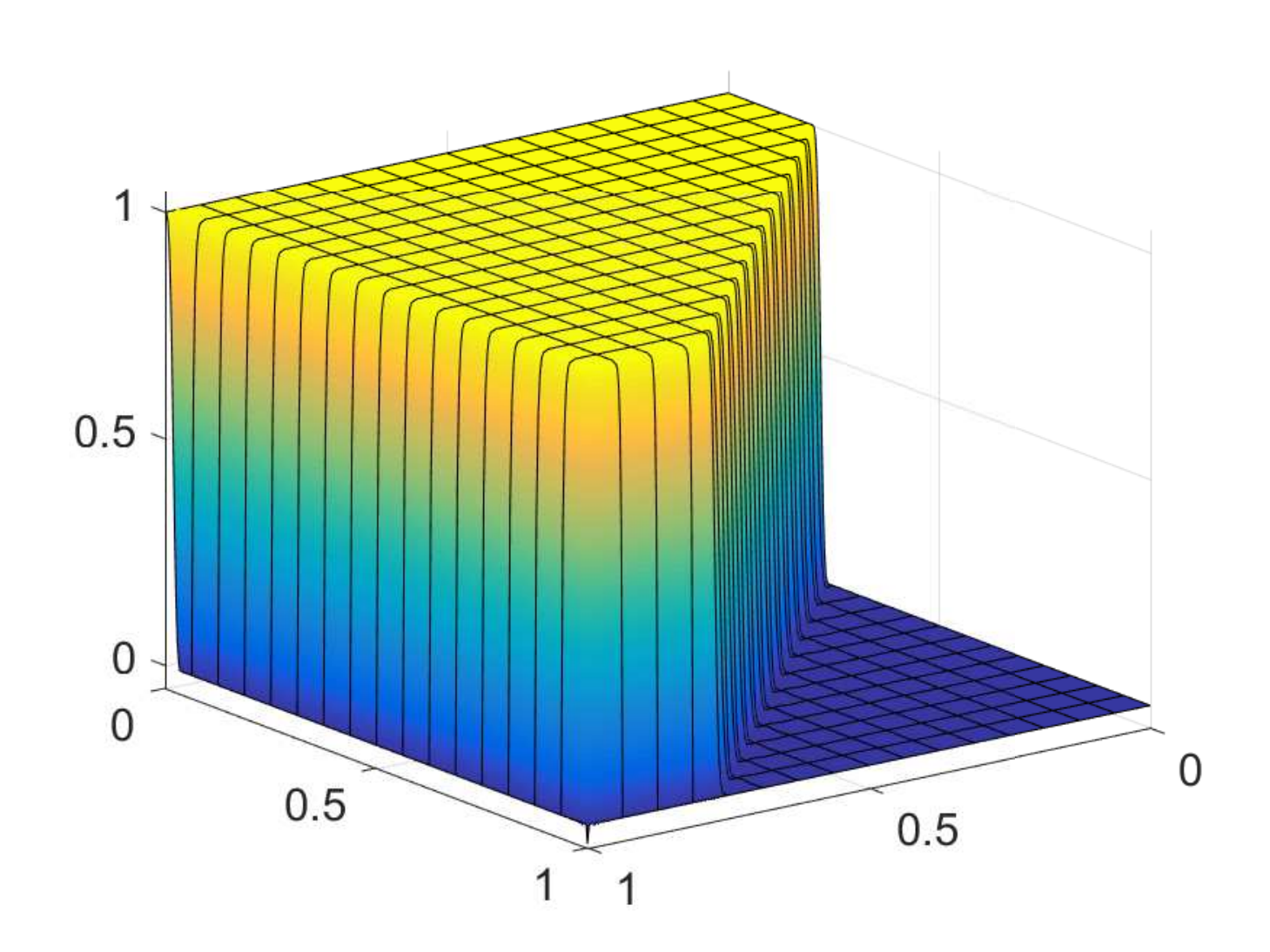}} \\
\subfigure[$\Omega$, THBox-spline]
{\includegraphics[trim = 2.5cm 0.75cm 2.25cm 1.cm, clip = true, height = 4cm]
{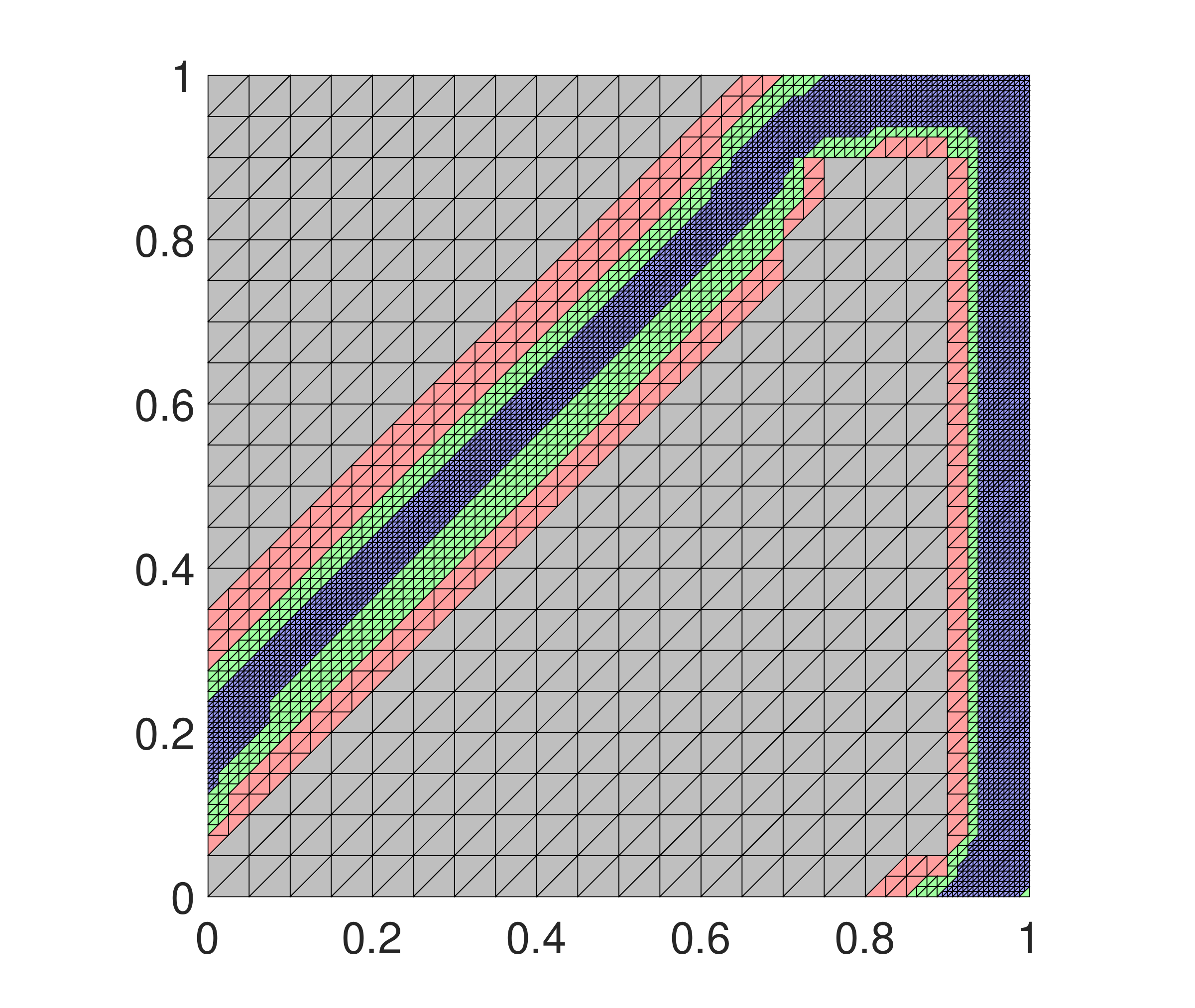}}
\subfigure[$\Omega_\Gamma$, THBox-spline]
{\includegraphics[trim = 2.5cm 0.75cm 2.25cm 1.cm, clip = true, height = 4cm]
{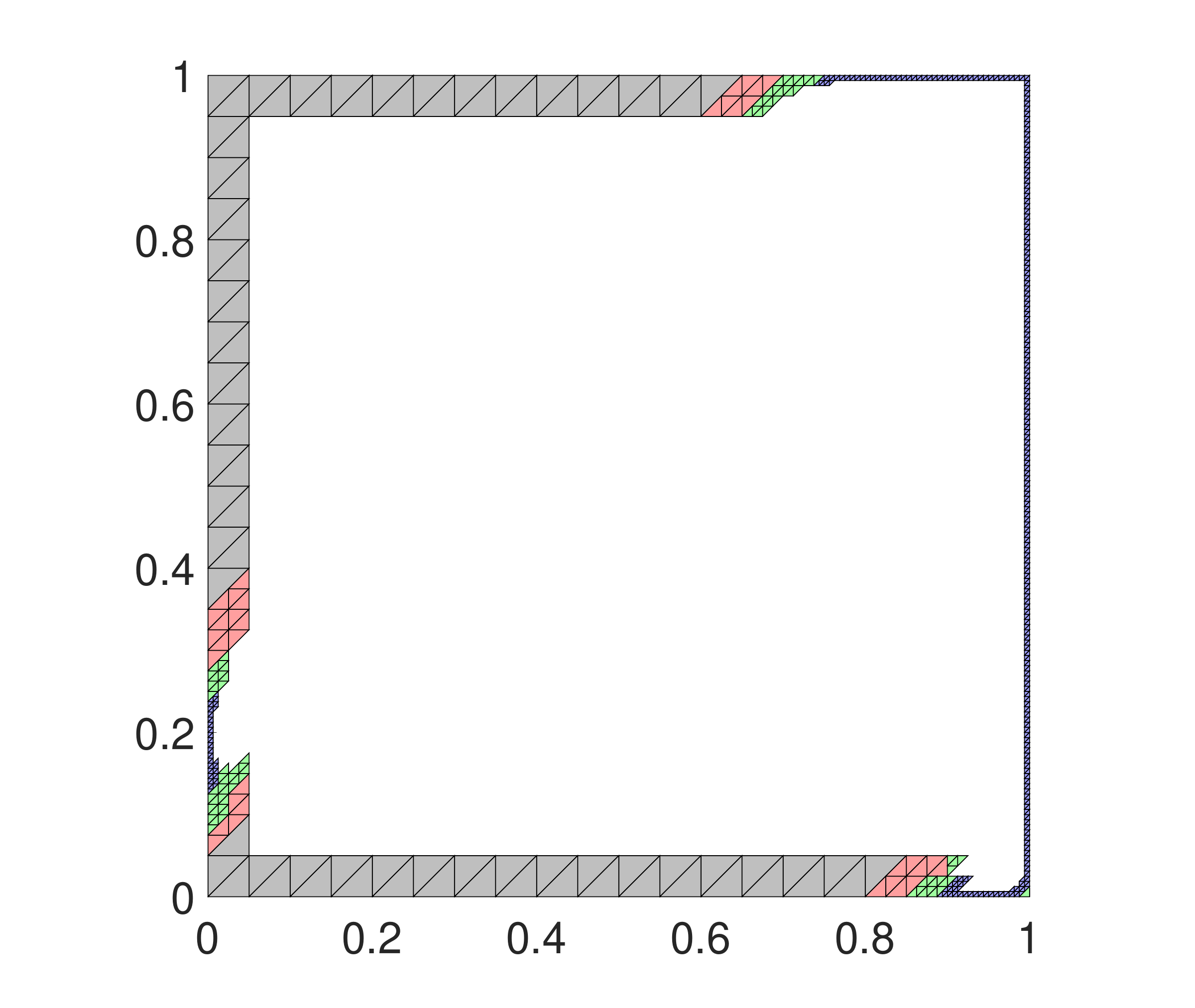}}
\subfigure[solution, THBox-spline]
{\includegraphics[trim = 0.87cm 0.5cm 0.77cm 0.75cm, clip = true, height=3.5cm]
{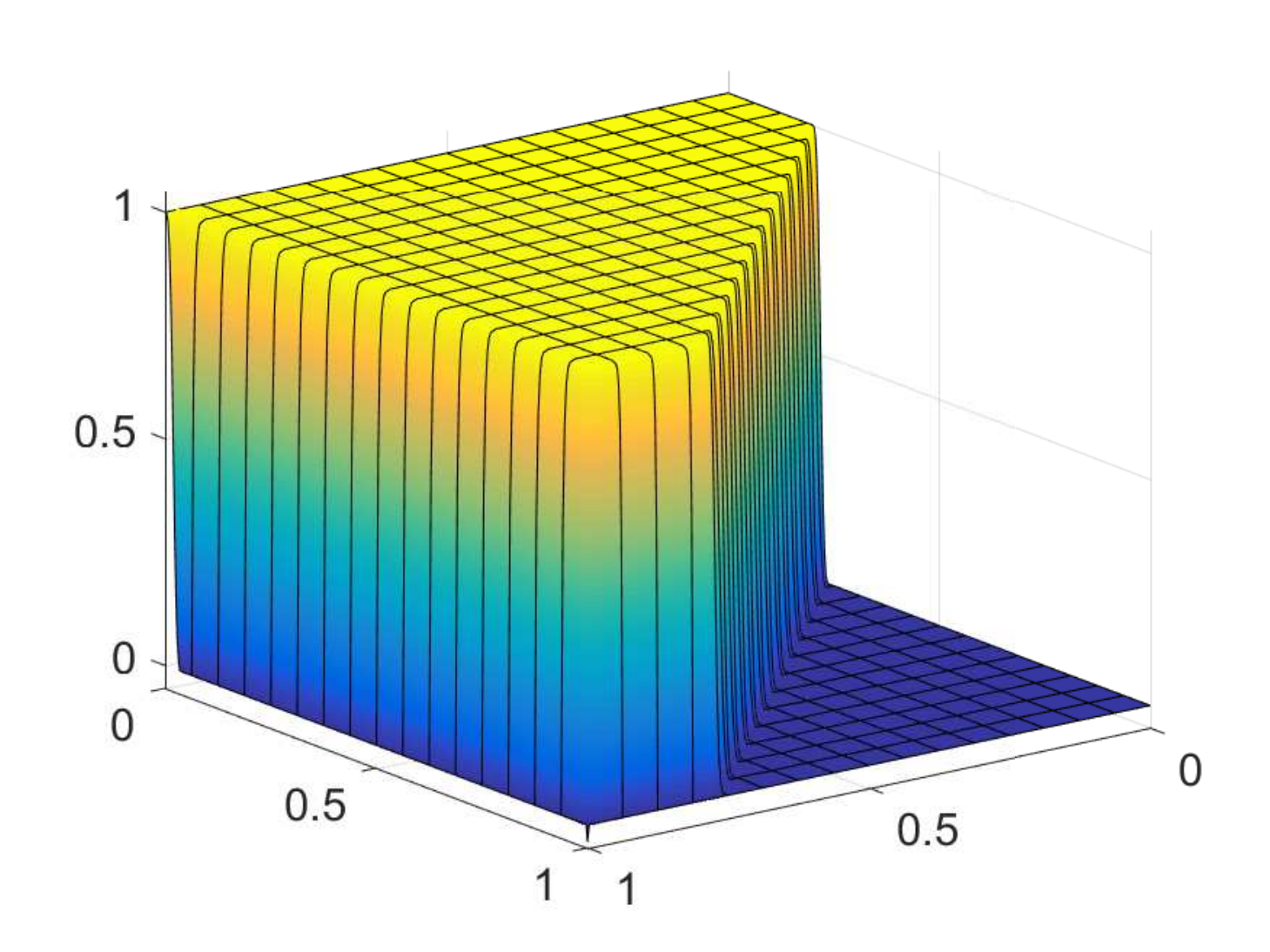}}
\caption{Advection-diffusion problem: Uniform and hierarchical meshes on $\hat \Omega = \Omega$ and $\hat \Omega_\Gamma = \Omega_\Gamma$, together with the obtained solutions and their errors ($N=4$).}
\label{fig:advection-strip}
\end{center}
\end{figure}

We now use hierarchical box splines to solve the advection-diffusion problem \eqref{problem1}
with $\kappa=10^{-6}$ and $\ba=(\cos\theta,\sin\theta)$, where $\theta=\pi/4$.
The domain with Dirichlet boundary conditions and the exact solution are shown in Figure~\ref{fig:advection}.
For system stabilization, we apply SUPG where the stabilization parameter is set to $\delta_h = h_{N-1} /\big(2 \max(\cos(\theta), \sin(\theta)) \|\ba\|\big)=h_{N-1}/\sqrt{2}$ 
(see also \cite{hughes2005}).
Moreover, to improve the quality of the solution, we replace the initial step function $g$ along the boundary with a smooth version, represented in terms of $C^2$ cubic B-splines.
{For the HBox and THBox cases, the local refinement is done automatically using a simple gradient based a posteriori error estimator. On each mesh cell $\pi$ of domain $\Omega$, the gradient indicator is computed as $\| \nabla u_h \|_{L^2(\pi)}$; see, e.g., \cite{john2000}. The gradient indicator is used in the cell marking strategy to prescribe the area of refinement: every cell that has the gradient indicator above a prescribed threshold is marked to be refined. Each area of refinement is suitably enlarged from the initial marked cell(s) so that at least one finer basis function is added on each marked cell. Starting from the initial uniform mesh, three automatic refinement steps are executed in the example.}
Figure~\ref{fig:advection-strip} shows the uniform and hierarchical three-directional meshes used for this problem with $N=4$ (left column), the different hierarchical boundary strip constructions (center column), and the corresponding box spline solutions (right column).
Again, in this case the hierarchical fixed boundary strip and the HBox-spline boundary strip produce exactly the same domain $\Omega_\Gamma$.
The quality of the approximation is usually measured in terms of the oscillations that appear in the neighborhood of the layers. Therefore, we report the minimum and maximum value of the box spline solutions in Table~\ref{tab:advection}, together with their degrees-of-freedom. 

\begin{table}[t!]
\begin{center}
\footnotesize
\begin{tabular}{c@{\hspace{.75em}}c{|}*{3}c{|}*{3}{c}}
&& \multicolumn{3}{c|}{max value} & \multicolumn{3}{c}{min value}\\
\hline
$N$ & $h_{N-1}$ & uniform & HBox & THBox & uniform & HBox & THBox\\
\hline
1 & 1/20 & $1.14\phantom{10}$ & $1.14\phantom{10}$ & $1.14\phantom{10}$ & $-0.15\phantom{1}$ & $-0.15\phantom{1}$ & $-0.15\phantom{1}$ \\
2 & 1/40 & $1.19\phantom{00}$ & $1.18$\phantom{00} & $1.18$\phantom{00} & $-0.046$ & $-0.049$ & $-0.046$\\
3 & 1/80 & $1.0072$ & $1.0068$ & $1.0072$ & $-0.042$ & $-0.046$ & $-0.042$\\
4 & 1/160 & $1.0030$ & $1.0077$ & $1.0083$ & $-0.037$ & $-0.041$ & $-0.037$
\end{tabular}

\medskip

\begin{tabular}{c@{\hspace{.75em}}
c{|}*{3}{c}{|}*{3}{c}}
& 
& \multicolumn{3}{c|}{dof} & \multicolumn{3}{c}{$\Omega_\Gamma$ dof}\\
\hline
$N$ & $h_{N-1}$ 
& uniform & HBox & THBox & uniform & HBox & THBox\\
\hline
1 & 1/20 
 & \phantom{00}527 & \phantom{0}527 & \phantom{0}527 & \phantom{0}302 & \phantom{0}302 & \phantom{0}302\\
2 & 1/40 
& \phantom{0}1847 & \phantom{0}782 & \phantom{0}782 & \phantom{0}622 & \phantom{0}485 & \phantom{0}446\\
3 & 1/80 
& \phantom{0}6887 & 1568 & 1610 & 1262 & \phantom{0}995 & \phantom{0}715\\
4 & 1/160 
 & 26567 & 4358 & 4400 & 2542 & 2660 & 1222
\end{tabular}

\caption{Advection-diffusion problem: Max/min values 
and degrees-of-freedom (dof), for uniform and hierarchical box splines with different types of boundary strips.}
\label{tab:advection}
\end{center}
\end{table}


\section{Conclusions} \label{sec:conclusions}

We have presented an adaptive isogeometric method with (truncated) hierarchical box splines using weakly imposed boundary conditions. 
By combining the locally uniform structure of the hierarchical box spline model with both the immersed boundary method and the isogeometric approach, we are not only able to effectively perform adaptive mesh refinement, but also to model non-trivial geometries.

The considered weak boundary formulation requires an appropriate selection of boundary strip where an additional flux unknown is enforced.  In order to do this, the strip must necessarily include all cells of the mesh that touch or intersect the boundary. In view of the possible interactions between different levels of resolutions in the adaptive hierarchical setting, this minimal set of boundary cells must be properly enlarged and different choices can be considered. 

Three different boundary strips have been introduced. 
The first strip simply consists of a fixed thickness $h_0$, independent of the number of refinements. 
To reduce the number of degrees-of-freedom of the overall system and the computational overhead, the strip should be thinner at regions where only finer resolutions are considered. 
%
In particular, the strip can be shrunk along the parts of the boundary where no basis functions of coarser levels are non-zero. Exploiting the supports of HBox-splines and THBox-splines results in the definition of adaptive strips that strongly reduce the number of degrees-of-freedom with respect to the fixed configuration. 

We have illustrated the capability of the presented hierarchical box spline model in isogeometric analysis with several numerical examples. Uniform and hierarchical box splines were used to solve a selection of advection-diffusion problems imposing weakly Dirichlet boundary conditions and using several non-trivial domains described by triangulated partitions and different domain mappings. These examples confirmed the optimal convergence rates of box splines on uniform and hierarchical meshes with different types of boundary strips.
Also, the number of degrees-of-freedom was substantially reduced by constructing {suitable} locally refined hierarchical meshes.

{The problem of high-quality domain parameterizations is crucial for the design of efficient isogeometric methods. It is also a non-trivial task, that becomes of remarkable interest when moving from (more common) surface parameterization algorithms to challenging volumetric configurations. Being defined as multivariate generalizations of univariate cardinal B-splines not necessarily restricted to rigid tensor-product structures, box splines offer an enhanced flexibility for the representation of computational domains.
For example, while non-singular tensor-product B-spline (and their rational extension) parameterizations may be developed in connection with planar domains that exhibit only four corners, singular solutions need to be taken into account when the number of corners is different from four. Box spline geometries may circumvent this problem but the development of optimal domain parameterizations in this context deserves a separate study and it is beyond the scope of this paper. 
For related work on domain parameterizations using splines on triangulations we refer to \cite{engvall2016,speleersM2015}. 
The aim of the current manuscript is solely to show the potential of the hierarchical box spline framework, suitably combined with simple geometric mappings, with respect to standard tensor-product counterparts.}

\section*{Acknowledgements}
This work was partially supported by the MIUR ``Futuro in Ricerca'' programme through the project DREAMS (RBFR13FBI3).

\section*{References}

\bibliographystyle{plain}

\end{document}